\documentclass[a4paper,11pt]{article}
\usepackage{amsfonts}
\usepackage{amsthm,amsfonts,epsfig,setspace}
\usepackage{colortbl,color}
\usepackage[margin=1in]{geometry}
\usepackage{epstopdf}
\usepackage{amssymb}
\usepackage{mathrsfs}
\usepackage{amsfonts,latexsym,bm}
\usepackage{amsmath}
\usepackage{amsthm}
\usepackage{float}
\usepackage{graphicx}
\usepackage{CJK}
\usepackage{graphics}
\usepackage{color}
\usepackage{subfigure}
\usepackage[english]{babel}

\usepackage{algorithm}
\usepackage{algpseudocode}

\makeatletter 
\@addtoreset{equation}{section}
\makeatother  

\newcommand {\pic}{\setlength{\unitlength}{1cm}
\begin{picture}(4,0)
\thicklines \put(-0.5,0){\line(1,0){12}}
\end{picture}
}

\title{A stochastic gradient descent approach with partitioned-truncated singular value decomposition for large-scale inverse problems of magnetic modulus data}

\author{
Wenbin Li  \thanks{\bf Corresponding Author: Wenbin Li} \thanks{School of Science, Harbin Institute of Technology, Shenzhen, Shenzhen 518055, China. Email: {\bf liwenbin@hit.edu.cn}} 
\and Kangzhi Wang \thanks{School of Science, Harbin Institute of Technology, Shenzhen, Shenzhen 518055, China. Email: {\bf 19S058004@stu.hit.edu.cn}}
\and Tingting Fan \thanks{School of Science, Harbin Institute of Technology, Shenzhen, Shenzhen 518055, China. Email: {\bf 19S058010@stu.hit.edu.cn}} 
}
\date{October 10, 2021}

\begin{document}
\thispagestyle{plain} \maketitle \vspace{0.5cm}

\section{Abstract}
We propose a stochastic gradient descent approach with partitioned-truncated singular value decomposition for large-scale inverse problems of magnetic modulus data. Motivated by a uniqueness theorem in gravity inverse problem and realizing the similarity between gravity and magnetic inverse problems, we propose to solve the level-set function modeling the volume susceptibility distribution from the nonlinear magnetic modulus data. To deal with large-scale data, we employ a mini-batch stochastic gradient descent approach with random reshuffling when solving the optimization problem of the inverse problem. We propose a stepsize rule for the stochastic gradient descent according to the Courant-Friedrichs-Lewy condition of the evolution equation. In addition, we develop a partitioned-truncated singular value decomposition algorithm for the linear part of the inverse problem in the context of stochastic gradient descent. Numerical examples illustrate the efficacy of the proposed method, which turns out to have the capability of efficiently processing large-scale measurement data for the magnetic inverse problem. A possible generalization to the inverse problem of deep neural network is discussed at the end.

\section{Introduction}
Solution of magnetic inverse problem aims to reconstruct susceptibility distributions from magnetic measurements, which is an important task in detections and surveys employing magnetic approaches. Applications include exploration of geological resources \cite{liold96}, detection of unexploded ordnance \cite{munbouulrbou07}, surveillance for naval warfare \cite{zhdcumwilpol12}, and others. A main difficulty of solving the magnetic inverse problem comes from its ill-posedness, and there is a strong non-uniqueness of solution for a given set of magnetic measurement data. To deal with the non-uniqueness, level-set approaches have been proposed for the inverse problems of magnetic intensity data \cite{liluqiali17,liqiali20} and magnetic gradient tensor data \cite{liqia20}, where it assumes a given constant susceptibility supported in an unknown domain and the level-set function is used to represent topological shape of the unknown domain. In this work, we adopt the level-set formulation and develop efficient computational method for massive datasets which were infeasible to process due to the restriction of computational resource. We will work on the large-scale inverse problem of nonlinear magnetic modulus data, and we propose a mini-batch stochastic gradient descent approach with partitioned-truncated singular value decomposition for the solution of the level-set function modeling the volume susceptibility distribution.

When meeting with large-scale dataset in realistic magnetic inverse problem, we realize that the computational resource becomes inadequate as both the memory access and the computational time are demanding. A direct and trade-off strategy is to take a subsample of the measurement dataset, but the problem is how to choose an appropriate subsample preserving the data feature under the restriction of computational resource. As a result, we look for the employment of full dataset with limited computational resource, and the stochastic approach comes to our mind. The stochastic gradient methods have achieved great success in large-scale machine learning where massive datasets are involved in the training process; e.g., see \cite{botcurnoc18} for a comprehensive review. The stochastic approaches use part of the dataset in every computation, which significantly saves memory usage and computational effort, whereas the full dataset is employed during each epoch if the without-replacement sampling strategy is adopted, or the full dataset is touched in expectation if the with-replacement sampling strategy is used \cite{safsha20}. Moreover, it is reported that with the help of noisy gradients, the stochastic approaches have the capability of escaping saddle points \cite{gehuajinyua15} and local minima \cite{kleliyua18} for non-convex optimization problems. Due to these advantages, we propose to introduce the stochastic strategy into the solution of the magnetic inverse problem. We use a mini-batch approach to reduce variances of stochastic gradients, and we formulate the iteration process by a partial differential equation and propose a stepsize rule according to its Courant-Friedrichs-Lewy condition. In addition, A partitioned-truncated singular value decomposition algorithm is developed to efficiently compute the linear part of the inverse problem in the context of stochastic gradient approach. The details will be shown in the part of main algorithm.

We consider the nonlinear modulus data in the magnetic inverse problem for the following two reasons. Firstly, the modulus data are nonlinear stacking of multiple components of the magnetic vector field, which are considered to be less sensitive to systematic measurement noises than the directionally dependent vector data. Secondly, a well-posedness theory of gravity inverse problem suggests that the modulus of gravity field can uniquely determine gravity source \cite{isa90,liqia21}, and since the magnetic inverse problem shares similarity to the gravity inverse problem, we believe that the measurement of magnetic modulus data is sufficient for the recovery of the susceptibility distribution. More details will be discussed in the formulation of the inverse problem of magnetic modulus data. In addition, we realize that the general form of the nonlinear magnetic inverse problem is similar to the architecture of neural network, and the proposed algorithm can be generalized to solve the inverse problem of deep neural network, which will be discussed in the section of conclusion and discussion.

\section{Inverse problem of magnetic modulus data}
Given an inducing  magnetic field of strength $B^0$ and direction $\hat{\mathbf{B}}^0$, the magnetic vector field resulting from a distribution of magnetic susceptibility $\kappa$ is modeled by the equation
\begin{equation} \label{eqn1}
\mathbf{B}(\mathbf{r})=\left(B_1(\mathbf{r}),\,B_2(\mathbf{r}),\,B_3(\mathbf{r})\right)^T=\frac{1}{4\pi}B^0
\int_{\Omega}\mathbf{K(r,\tilde{r})}\,\kappa(\mathbf{\tilde{r}})\,\mathrm{d}\mathbf{\tilde{
r}},
\end{equation}
with
\begin{equation} \label{eqn2}
\mathbf{K(r,\tilde{r})}=\nabla_{\mathbf{r}}\left(\hat{\mathbf{B}}^0\cdot\nabla_{\mathbf{r}}
\left(\frac{1}{|\mathbf{r}-\mathbf{\tilde{r}}|}\right)\right), 
\end{equation}
where $\mathbf{r}$ and $\mathbf{\tilde{r}}$ denote 3D spatial coordinates: $\mathbf{r}\in\Gamma$ 
(the data  measurement boundary), $\mathbf{\tilde{r}}\in\Omega$. The nonlinear magnetic modulus is given by,
\begin{equation} \label{eqn3}
d(\mathbf{r})=\left(\sum_{s=1}^3B_s(\mathbf{r})^2\right)^{\frac{1}{2}}\,;
\end{equation}
since $d(\mathbf{r})$ is nonlinear stacking of multiple components of the magnetic vector field $\mathbf{B}(\mathbf{r})$, it is presumably less sensitive to systematic measurement noises than $\mathbf{B}(\mathbf{r})$ which is directionally dependent. The nonlinear inverse problem reads as follows: given the modulus data $d(\mathbf{r})$ on the measurement boundary $\Gamma$, reconstruct the magnetic susceptibility parameter $\kappa(\mathbf{\tilde{r}})$ in the domain $\Omega$.

It is well known that the solution of magnetic inverse problem is non-unique, and there exist infinitely many susceptibility distributions which can reproduce a given set of magnetic data. To deal with the non-uniqueness, we consider the volume susceptibility distribution in the following form,
\begin{equation} \label{eqn4}
\kappa(\mathbf{\tilde{r}})=\kappa_0\chi_D(\mathbf{\tilde{r}})\,,
\end{equation}
where $\kappa_0$ is a constant that approximates the average value of the susceptibility, and $\chi_D$ is the characteristic function of a domain $D$: 
$\chi_D(\mathbf{\tilde{r}})=1$, $\mathbf{\tilde{r}}\in D$; $\chi_D(\mathbf{\tilde{r}})=0$, $\mathbf{\tilde{r}}\notin D$. We assume that the constant $\kappa_0$ is known as a priori information, and invert for the domain $D$ to recover the shape of the magnetic source. The setup is motivated by a well-posedness theory in gravity inverse problem \cite{isa90,liqia21}. Let $U$ denote the gravitational potential generated by a volume mass distribution $\mu=\rho_0\chi_D$; if $\rho_0$ is a known constant and the domain $D$ satisfies some geometric constraints, the exterior gravity modulus $|\nabla U|$ uniquely determines the shape of $D$. Realizing that the magnetic inverse problem shares similarity to the gravity inverse problem, see, e.g. \cite{liqia20}, we consider the volume susceptibility distribution as shown in equation (\ref{eqn4}), and expect that the domain $D$ can be well recovered from the magnetic modulus data $d$ as the constant value $\kappa_0$ is given.

Similar to our previous works \cite{liluqiali17,liqia20,liqiali20}, we adopt a level-set formulation for the volume susceptibility distribution,
\begin{equation} \label{eqn5}
\kappa(\mathbf{\tilde{r}})=\kappa_0 H(\phi(\mathbf{\tilde{r}}))\,,
\end{equation}
where $\phi(\mathbf{\tilde{r}})$ is the level-set function used to depict the topological shape of $D$, and $H(\cdot)$ is the Heaviside function: $H(x)=1$, $x\ge0$; $H(x)=0$, $x<0$. In level-set computation, $\phi(\mathbf{\tilde{r}})$ is maintained to be a continuous signed-distance function to the boundary of $D$,
\begin{equation} \label{eqn6}
\phi(\mathbf{\tilde{r}})=\left\{
\begin{array}{ccc}
\mathrm{dist}(\mathbf{\tilde{r}},\partial D)&,&\mathbf{\tilde{r}}\in D\\
-\mathrm{dist}(\mathbf{\tilde{r}},\partial D)&,&\mathbf{\tilde{r}}\in \bar{D}^c
\end{array}
\right.\,,
\end{equation}
where $\mathrm{dist}(\mathbf{\tilde{r}},\partial D)$ denotes the distance between $\mathbf{\tilde{r}}$ and $\partial D$. As a result, equation (\ref{eqn5}) is a natural expression for the volume susceptibility distribution $\kappa$ as shown in equation (\ref{eqn4}), and the zero level-set $\{\mathbf{\tilde{r}}\mid\phi(\mathbf{\tilde{r}})=0\}$ indicates the location of the boundary $\partial D$. The inverse problem is then transformed into the following form: given the magnetic modulus data $d(\mathbf{r})$ on the measurement boundary $\Gamma$, recover the level-set function $\phi(\mathbf{\tilde{r}})$ in $\Omega$ so that the magnetic susceptibility $\kappa(\mathbf{\tilde{r}})$ is evaluated by equation (\ref{eqn5}).

\section{Main algorithm}
\subsection{Data fitting and optimization}
The inverse problem is solved by fitting the magnetic modulus data with an optimal solution of the level-set function. Let $d^*(\mathbf{r})$ denote the measured modulus data at $\mathbf{r}\in\Gamma$, and $d(\phi,\mathbf{r})$ denote the predicting data with the value of $\phi(\mathbf{\tilde{r}})$. The usually used $l^2$ fitting function is defined as follows,
\begin{equation} \label{eqn7}
E(\phi)=\frac{1}{|\Gamma|}\int_\Gamma\frac{\left(d(\phi,\mathbf{r})-d^*(\mathbf{r})\right)^2}{2}\,\mathrm{d}\mathbf{r}
\end{equation}
where $|\Gamma|$ denotes the area of the measurement surface $\Gamma$. Equation (\ref{eqn7}) can be viewed as an expectation of the data misfit on $\Gamma$:
\begin{equation} \label{eqn8}
\mathbb{E}(f(\phi,\mathbf{r}))=\int_\Gamma f(\phi,\mathbf{r})\, p(\mathbf{r})\,\mathrm{d}\mathbf{r}
\end{equation}
with
\begin{equation} \label{eqn9}
p(\mathbf{r})=\frac{1}{|\Gamma|} \qquad \mathrm{and}\qquad f(\phi,\mathbf{r})=\frac{1}{2}\left(d(\phi,\mathbf{r})-d^*(\mathbf{r})\right)^2\,.
\end{equation}

A regularization term is added to the data fitting function, so that minimizing the total objective function leads to the admissible solution of $\phi$. We propose to use the $H^1$ semi-norm of $\phi$ as the regularization term:
\begin{equation} \label{eqn10}
E_r(\phi)=\frac{1}{2}\int_{\Omega}|\nabla\phi(\mathbf{\tilde{r}})|^2\,\mathrm{d}\mathbf{\tilde{r}}\,.
\end{equation}
As discussed in \cite{liqia21}, minimizing $E_r(\phi)$ has the effect of shrinking the measure (area or length) of the interface characterized by the zero level-set $\{\mathbf{\tilde{r}}\mid\phi(\mathbf{\tilde{r}})=0\}$, so that it avoids forming sharp oscillations on the interface.

The level-set function $\phi(\mathbf{\tilde{r}}))$ is then recovered by solving the optimization problem:
\begin{equation} \label{eqn11}
\arg\min_{\phi}\, E_t(\phi):=\mathbb{E}(f(\phi,\mathbf{r}))+\alpha E_r(\phi)\,,
\end{equation}
where $\alpha$ is a weighting parameter that controls the amount of regularization applied.

\subsection{Mini-batch stochastic gradient descent with random reshuffling}
The gradient-descent based approach is favorable for solving the optimization problem of equation (\ref{eqn11}) with large scale. The derivative of $E_t(\phi)$ is given by the formula,
\begin{equation} \label{eqn12}
\frac{\partial E_t(\phi)}{\partial\phi}=\frac{\partial \mathbb{E}(f(\phi,\mathbf{r}))}{\partial\phi}+\alpha \frac{\partial E_r(\phi)}{\partial\phi}=\frac{\partial \mathbb{E}(f(\phi,\mathbf{r}))}{\partial\phi}-\alpha \Delta\phi
\end{equation}
with
\begin{equation} \label{eqn13}
\frac{\partial \mathbb{E}(f(\phi,\mathbf{r}))}{\partial\phi}=\int_\Gamma \frac{\partial f(\phi,\mathbf{r})}{\partial\phi}\, p(\mathbf{r})\,\mathrm{d}\mathbf{r}\,;
\end{equation}
the level-set function $\phi$ is updated according to the negative gradient direction,
\begin{equation} \label{eqn_gradient}
\frac{\partial\phi}{\partial t}=-\frac{\partial E_t(\phi)}{\partial\phi}=-\frac{\partial \mathbb{E}(f(\phi,\mathbf{r}))}{\partial\phi}+\alpha \Delta\phi
\end{equation}
where $t$ denotes the pseudo time variable that enables the evolution of $\phi$. As $f(\phi,\mathbf{r})$ is defined in the form of equation (\ref{eqn9}), we have
\begin{equation} \label{eqn14}
\frac{\partial f(\phi,\mathbf{r})}{\partial\phi}=\left(d(\phi,\mathbf{r})-d^*(\mathbf{r})\right)\frac{\partial d(\phi,\mathbf{r})}{\partial\phi}\,.
\end{equation}
Considering that $d(\phi,\mathbf{r})$ is related to $\phi$ through equations (\ref{eqn1})-(\ref{eqn3}) and equation (\ref{eqn5}), formally we can compute $\frac{\partial d}{\partial\phi}$ by the chain rule,
\begin{equation} \label{eqn15}
\frac{\partial d(\phi,\mathbf{r})}{\partial\phi(\mathbf{\tilde{r}})}=\sum_{s=1}^3\frac{\partial d}{\partial B_s}\frac{\partial B_s}{\partial\kappa}\frac{\partial\kappa}{\partial\phi}=\frac{B^0\kappa_0}{4\pi}\delta(\phi(\mathbf{\tilde{r}}))\sum_{s=1}^3\frac{B_s(\mathbf{r})}{d(\mathbf{r})}K_s(\mathbf{r}, \mathbf{\tilde{r}})\,.
\end{equation}
Here, $\delta(\phi)$ denotes the Dirac delta function arising from the derivative of the Heaviside function $H(\phi)$, and $K_s(\mathbf{r}, \mathbf{\tilde{r}})$ denotes the $s$-th component of the magnetic integral kernel $\mathbf{K(r,\tilde{r})}$ so that
\begin{equation} \label{eqn16}
\big(K_1(\mathbf{r}, \mathbf{\tilde{r}}), K_2(\mathbf{r}, \mathbf{\tilde{r}}), K_3(\mathbf{r}, \mathbf{\tilde{r}})\big)=\mathbf{K(r,\tilde{r})}\,.
\end{equation}
Substituting equation (\ref{eqn15}) into equation (\ref{eqn14}), we have
\begin{equation} \label{eqn17}
\frac{\partial f(\phi,\mathbf{r})}{\partial\phi}=\frac{B^0\kappa_0}{4\pi}\,\delta(\phi(\mathbf{\tilde{r}}))\sum_{s=1}^3\frac{B_s(\mathbf{r})}{d(\mathbf{r})}\left(d(\phi,\mathbf{r})-d^*(\mathbf{r})\right) K_s(\mathbf{r}, \mathbf{\tilde{r}})\,.
\end{equation}
The derivative of the expected data misfit $\mathbb{E}(f(\phi,\mathbf{r}))$ is then obtained by putting equation (\ref{eqn17}) into equation (\ref{eqn13}).

In physical surveys, the measurement data are sampled at discrete coordinates. Denoting the set of measurement points by $\{\mathbf{r}_i\mid i=1,2,\cdots,M\}$, the expected data misfit $\mathbb{E}(f(\phi,\mathbf{r}))$ reduces to the empirical data misfit:
\begin{equation}  \label{eqn18}
f^M(\phi)=\frac{1}{M}\sum_{i=1}^M f(\phi,\mathbf{r}_i)\,.
\end{equation}
There will be difficulties in minimizing  $f^M(\phi)$ directly when the size of data becomes large. To illustrate it, we derive $\frac{\partial f^M}{\partial\phi}$ from equation (\ref{eqn18}) and equation (\ref{eqn17}),
\begin{equation} \label{eqn19}
\frac{\partial f^M(\phi(\mathbf{\tilde{r}}))}{\partial\phi(\mathbf{\tilde{r}})}=\frac{B^0\kappa_0}{4\pi M}\,\delta(\phi(\mathbf{\tilde{r}}))\sum_{s=1}^3\sum_{i=1}^M\frac{B_s(\mathbf{r}_i)}{d(\mathbf{r}_i)}\left(d(\phi,\mathbf{r}_i)-d^*(\mathbf{r}_i)\right) K_s(\mathbf{r}_i, \mathbf{\tilde{r}})\,.
\end{equation}
Here, $\mathbf{\tilde{r}}$ denotes the grid points in the computational domain $\Omega$. In a regular-size computation, the number of grid points is in the order of $10^4$, e.g. 40 grids each direction in a 3D domain; if the number of data measurements is also in this order, i.e. $M\sim 10^4$, the dimension of the integral kernel $K_s(\mathbf{r}_i, \mathbf{\tilde{r}})$ will be grater than $10^8$. It becomes demanding in memory accessing and computational time as the derivative is computed repeatedly in the gradient-descent based iterative algorithm, and the computation will be infeasible if the size of data and the number of grid points grow larger. 

We propose to utilize a mini-batch stochastic gradient descent approach \cite{botcurnoc18,lizhachesmo14} to solve the optimization problem of equation (\ref{eqn11}) with large scale. Instead of using all the measurement data, we employ a part of them to realize the expected data misfit $\mathbb{E}(f(\phi,\mathbf{r}))$ and its derivative:
\begin{eqnarray}
f^{S_k}(\phi)&=&\frac{1}{|S_k|}\sum_{i\in S_k}f(\phi,\mathbf{r}_i)\,, \label{eqn20}\\
\frac{\partial f^{S_k}(\phi)}{\partial\phi}&=&\frac{1}{|S_k|}\sum_{i\in S_k}\frac{\partial f(\phi,\mathbf{r}_i)}{\partial\phi}\,, \label{eqn21}
\end{eqnarray}
where $S_k\subset\{1,2,\cdots,M\}$ denotes a mini-batch of the index set, and $|S_k|$ denotes the batch size, i.e. the number of elements in $S_k$. Using $\frac{\partial f^{S_k}(\phi)}{\partial\phi}$ instead of $\frac{\partial f^M(\phi)}{\partial\phi}$ to evaluate $\frac{\partial \mathbb{E}(f(\phi,\mathbf{r}))}{\partial\phi}$ in equation (\ref{eqn_gradient}), the level-set function $\phi$ is then evolved according to the following formula,
\begin{equation} \label{eqn22}
\frac{\partial\phi}{\partial t}=-\frac{\partial f^{S_k}(\phi)}{\partial\phi}+\alpha \Delta\phi\,.
\end{equation}

The sampling strategies of $S_k$ can be generally divided into two categories: with-replacement sampling and without-replacement sampling \cite{bot09,bot12,safsha20}. In the with-replacement sampling, $S_k$ is uniformly drawn from the data index set $\{1,2,\cdots,M\}$ for each iteration, so that $\mathbb{E}\left(f^{S_k}(\phi)\right)=\mathbb{E}(f(\phi,\mathbf{r}))$ and $\mathbb{E}\big(\frac{\partial f^{S_k}(\phi)}{\partial\phi}\big)=\frac{\partial\mathbb{E}(f(\phi,\mathbf{r}))}{\partial\phi}$. In this case, the evolution equation (\ref{eqn22}) can be seen as a noisy version of the exact updating as shown in equation (\ref{eqn_gradient}); for instance, if one sets $|S_k|=1$, i.e. uniformly picking single sample from the data set for each iteration, the algorithm reduces to the classical stochastic gradient descent approach, which is also referred to as stochastic approximation in some context \cite{fu02,nemjudlansha09,robmon51}. The with-replacement sampling strategy is preferred by theoretical analysis since it realizes the exact updating rule in expectation. However, it suffers from practical drawbacks such as requiring truly random data access and hence longer runtime \cite{safsha20}. In practice, it is common to use the without-replacement sampling strategy, where the measurement data are employed sequentially in some random or even deterministic order. 

We propose to use the without-replacement sampling with random reshuffling in our mini-batch stochastic gradient descent algorithm. The details are summarized in Algorithm 1. In every epoch, the data set is randomly permuted and then partitioned into mini-batches, and the iterative algorithm employs the mini-batches successively in order to evaluate the stochastic gradient descent. The shuffling strategy ensures that every sampling point is reached during each epoch. The mini-batch strategy has the capability of reducing variances when estimating the stochastic gradients \cite{botcurnoc18,lizhachesmo14}. In step 7 of Algorithm 1, the level-set function $\phi$ is reinitialized to maintain the shape of the signed-distance function as shown in equation (\ref{eqn6}). This is a standard operation in level-set computation, which aims to make $\phi$ sufficiently smooth without changing the location of the zero level set. For the details, we refer readers to our previous works \cite{liqia21,liluqiali17} and textbooks on level-set methods \cite{oshfed06}.

\pic\\
\noindent { \textsf{ Algorithm 1. Mini-batch stochastic gradient descent algorithm with random reshuffling for the inverse problem of magnetic modulus data.
\begin{algorithmic}[1]
\State Initialize the level-set function $\phi$\,; set iteration number $\overline{k}=0$\,.
\For {epoch $e=1,2,\cdots,N_e$}
   \State Sample a random permutation $\{\sigma(1),\sigma(2),\cdots,\sigma(M)\}$ of the index set $\{1,2,\cdots,M\}$.
   \State With a given batch size $b$, partition the shuffled index set into $n=\frac{M}{b}$ subsets: $S_k=\{\sigma((k-1)b+1),\cdots,\sigma(kb)\}$, $k=1,2,\cdots,n$.
   \For {$k=1$ to $n$}
       \State Update $\phi$ according to equation (\ref{eqn22}): $\phi:=\phi+\Delta t\big(\frac{\partial\phi}{\partial t}\big)$.
       \State Reinitialize $\phi$ to maintain the signed distance property.
       \State Iteration number $\overline{k}:=\overline{k}+1$.
   \EndFor
\EndFor
\end{algorithmic}
} } \pic

\subsection{Stepsize rule}
The choice of stepsize is an important aspect in the stochastic gradient descent algorithm \cite{nemjudlansha09,tanmadaiqia16}. An inappropriate choice of stepsize can have a disastrous effect on the performance of convergence. In this work, we propose a stepsize rule according to the Courant-Friedrichs-Lewy (CFL) condition of the evolution equation (\ref{eqn22}).

Similar to equation (\ref{eqn19}), we have
\begin{equation} \label{eqn23}
\frac{\partial f^{S_k}(\phi)}{\partial\phi}=\frac{B^0\kappa_0}{4\pi |S_k|}\,\delta(\phi(\mathbf{\tilde{r}}))\sum_{s=1}^3\sum_{i\in{S_k}}\frac{B_s(\mathbf{r}_i)}{d(\mathbf{r}_i)}\left(d(\phi,\mathbf{r}_i)-d^*(\mathbf{r}_i)\right) K_s(\mathbf{r}_i, \mathbf{\tilde{r}})\,.
\end{equation}
We use a numerical version of the Dirac delta function \cite{zhachamerosh96},
\begin{equation} \label{eqn24}
\delta(\phi)\doteq\delta_{\epsilon}(\phi)=\chi_{T_\epsilon}|\nabla\phi|\,,
\end{equation}
where $\chi_{T_\epsilon}$ denotes the characteristic function of $T_{\epsilon}$, and $T_{\epsilon}=\{\tilde{\mathbf{r}}: |\phi(\tilde{\mathbf{r}})|<\epsilon\}$ is a neighborhood of the zero level set. Substituting (\ref{eqn24}) into (\ref{eqn23}), the stochastic gradient can be evaluated in the following way,
\begin{equation} \label{eqn25}
\frac{\partial f^{S_k}(\phi)}{\partial\phi}=V_n|\nabla\phi|
\end{equation}
with
\begin{equation} \label{eqn26}
V_n:=\frac{B^0\kappa_0}{4\pi |S_k|}\,\chi_{T_\epsilon}\sum_{s=1}^3\sum_{i\in{S_k}}\frac{B_s(\mathbf{r}_i)}{d(\mathbf{r}_i)}\left(d(\phi,\mathbf{r}_i)-d^*(\mathbf{r}_i)\right) K_s(\mathbf{r}_i, \mathbf{\tilde{r}})\,.
\end{equation}
The evolution equation (\ref{eqn22}) then reduces to 
\begin{equation} \label{eqn27}
\frac{\partial\phi}{\partial t}=-V_n|\nabla\phi|+\alpha \Delta\phi\,,
\end{equation}
which can be viewed as a Hamilton-Jacobi equation with an artificial viscosity term. Discretizing $\frac{\partial\phi}{\partial t}$ using a direct forward Euler scheme: $\frac{\partial\phi}{\partial t}\approx\frac{\phi^{n+1}-\phi^n}{\Delta t}$, the CFL condition for stability is
\begin{equation} \label{eqn28}
\Delta t\left(\frac{\max|V_n|}{\min\{\Delta x,\Delta y,\Delta z\}}+\frac{2\alpha}{\Delta x^2}+\frac{2\alpha}{\Delta y^2}+\frac{2\alpha}{\Delta z^2}\right)<1\,.
\end{equation}
The stepsize of the mini-batch stochastic gradient descent algorithm is taken as
\begin{equation} \label{eqn29}
\Delta t=C \left(\frac{\max|V_n|}{\min\{\Delta x,\Delta y,\Delta z\}}+\frac{2\alpha}{\Delta x^2}+\frac{2\alpha}{\Delta y^2}+\frac{2\alpha}{\Delta z^2}\right)^{-1}\,,
\end{equation}
where $C$ denotes a constant with $C\in(0,1)$. In practice, as the regularization parameter $\alpha$ can be very small, $\frac{\alpha}{\min\{\Delta x,\Delta y,\Delta z\}}\ll \max|V_n|$, one can drop the terms including $\alpha$ in equation (\ref{eqn29}), and the stepsize is taken as
\begin{equation} \label{eqn29_1}
\Delta t=C \left(\frac{\max|V_n|}{\min\{\Delta x,\Delta y,\Delta z\}}\right)^{-1}\,.
\end{equation}


\subsection{Partitioned-truncated SVD for matrix multiplications in the stochastic gradient descent}
Large-scale matrix multiplications arise from the evaluation of gradient direction (as shown in equation (\ref{eqn26})) and the discretization of the magnetic integral (\ref{eqn1}). Considering formula (\ref{eqn2}) and formula (\ref{eqn16}), the magnetic integral kernel has the following form,
\begin{equation} \label{eqn30}
\big(K_1(\mathbf{r}, \mathbf{\tilde{r}}), K_2(\mathbf{r}, \mathbf{\tilde{r}}), K_3(\mathbf{r}, \mathbf{\tilde{r}})\big)=\frac{1}{|\mathbf{r}-\mathbf{\tilde{r}}|^3}\left[\frac{3\left(\hat{\mathbf{B}}^0\cdot(\mathbf{r}-\mathbf{\tilde{r}})\right)(\mathbf{r}-\mathbf{\tilde{r}})}{|\mathbf{r}-\mathbf{\tilde{r}}|^2}-\hat{\mathbf{B}}^0\right]\,,
\end{equation}
which decays rapidly as the distance $|\mathbf{r}-\tilde{\mathbf{r}}|$ is increasing. This decaying property is common to the integral kernel of potential-field data, i.e. magnetic and gravity data. We have proposed a partitioned-truncated singular value decomposition (SVD) algorithm to develop low rank approximation to the decaying integral kernel in several different setups and for various types of potential-field data \cite{luleuqia15,liluqiali17,liqia21}. Here, we propose a variation of the partitioned-truncated SVD so that it is compatible with the mini-batch stochastic gradient descent approach.

\subsubsection{Matrix multiplications involving the decaying kernel}
To evaluate the stochastic gradient as shown in equation (\ref{eqn26}), we consider the following formulation,
\begin{equation} \label{eqn31}
G^{S_k}_s(\mathbf{\tilde{r}}):=\chi_{T_\epsilon}\sum_{i\in{S_k}}\frac{B_s(\mathbf{r}_i)}{d(\mathbf{r}_i)}\left(d(\phi,\mathbf{r}_i)-d^*(\mathbf{r}_i)\right) K_s(\mathbf{r}_i, \mathbf{\tilde{r}})\,, \qquad S_k\subset\{1,2,\cdots,M\},\ s\in\{1,2,3\}\,.
\end{equation}
To evaluate $G^{S_k}_s$ we need to compute the magnetic components $B_s(\mathbf{r}_i)$: $i\in{S_k},\,s\in\{1,2,3\}$, and we extract the following integral,
\begin{equation} \label{eqn32}
I^{S_k}_s(\mathbf{r}_i):=\int_{\Omega}K_s(\mathbf{r}_i, \mathbf{\tilde{r}})\,\kappa(\mathbf{\tilde{r}})\,\mathrm{d}\mathbf{\tilde{r}}\,, \qquad i\in S_k\subset\{1,2,\cdots,M\},\ s\in\{1,2,3\}\,.
\end{equation}
Let the computational domain $\Omega\subset\mathbf{R}^3$ be a rectangular domain, $\Omega=[0,X]\times[0,Y]\times[-Z,0]$, and let the measurement boundary $\Gamma$ satisfies $\Gamma\subset\{\mathbf{r}=(x,y,z)\mid z>0\}$. This setup corresponds to a common situation where the magnetic measurement is above the survey domain. Suppose that the domain $\Omega$ is uniformly discretized into $N=n_x n_y n_z$ grid points: $\mathbf{\tilde{r}}_1,\mathbf{\tilde{r}}_2,\cdots,\mathbf{\tilde{r}}_N$, and recall that there are $M$ measurement points on $\Gamma$ and the mini-batch size $|S_k|=b$. After discretization, the evaluation of $I^{S_k}_s(\mathbf{r}_i)$ and $G^{S_k}_s(\mathbf{\tilde{r}})$ reduces to the evaluation of matrix multiplications,
\begin{eqnarray}
\mathbf{I}^{S_k}_s&=&\mathbf{K}^{S_k}_s\hspace{0.5pt}\mathbf{\kappa}^T\,,\label{eqn33} \\ 
\mathbf{G}^{S_k}_s&=&\mathbf{d}^{S_k}_s  \mathbf{K}^{S_k}_s \mathbf{D}\,, \label{eqn34}
\end{eqnarray}
where
\begin{eqnarray*}
\mathbf{K}^{S_k}_s&=&\left[K_s(\mathbf{r}_i, \mathbf{\tilde{r}}_j)\right]_{b\times N}\,, \quad i\in S_k,\ j=1,\cdots, N\,,\\
\mathbf{d}^{S_k}_s&=&\left[\frac{B_s(\mathbf{r}_i)}{d(\mathbf{r}_i)}\left(d(\phi,\mathbf{r}_i)-d^*(\mathbf{r}_i)\right)\right]_{1\times b}\,,\quad i\in S_k\,,\\
\mathbf{\kappa}&=&\left[\kappa(\mathbf{\tilde{r}}_1),\cdots,\kappa(\mathbf{\tilde{r}}_N) \right]\,,\\
\mathbf{D}&=&\mathrm{diag}\left\{\chi_{T_\epsilon}(\mathbf{\tilde{r}}_1),\cdots,\chi_{T_\epsilon}(\mathbf{\tilde{r}}_N)\right\}\,.
\end{eqnarray*}

\subsubsection{Preconditioning on the full kernel matrix}
The multiplication matrix $\mathbf{K}^{S_k}_s$ appearing in (\ref{eqn33}) and (\ref{eqn34}) is a sub-matrix of the full kernel matrix $\mathbf{K}_s=\left[K_s(\mathbf{r}_i, \mathbf{\tilde{r}}_j)\right]_{M\times N}\,, \ i=1,\cdots, M,\ j=1,\cdots, N$. Considering that the full kernel matrix $\mathbf{K}_s$ is always constant once the measurement points $\mathbf{r}_i$ and the grid points $\mathbf{\tilde{r}}_j$ are fixed, we perform preconditioning including partitioning and truncating on $\mathbf{K}_s$.

The component $K_s(\mathbf{r}_i, \mathbf{\tilde{r}}_j)$ decays rapidly as $|\mathbf{r}_i-\mathbf{\tilde{r}}_j|$ increases. Since $\mathbf{r}_i\in\Gamma\subset\{\mathbf{r}=(x,y,z)\mid z>0\}$, we order the grid points $\{\mathbf{\tilde{r}}_1,\mathbf{\tilde{r}}_2,\cdots,\mathbf{\tilde{r}}_N\}$ according to their z-coordinates so that $\tilde{z}_1\le\tilde{z}_2\le\cdots\le\tilde{z}_N$, where we denote $\mathbf{\tilde{r}}_j=(\tilde{x}_j,\tilde{y}_j,\tilde{z}_j)$. The set of grid points in the computational domain is then partitioned into $n_z$ subsets,
\begin{equation} \label{eqn35}
\{\mathbf{\tilde{r}}_1,\mathbf{\tilde{r}}_2,\cdots,\mathbf{\tilde{r}}_N\}=\bigcup_{h=1}^{n_z}\{\tilde{\mathbf{r}}^h_j \mid 1\le j\le \tilde{M}:=n_x n_y\}
\end{equation}
where $\tilde{\mathbf{r}}^h_j=\tilde{\mathbf{r}}_{j+(h-1)\tilde{M}}$ and so $\tilde{z}_j^h=(h-1)\Delta z$. The full kernel matrix $\mathbf{K}_s=\left[K_s(\mathbf{r}_i, \mathbf{\tilde{r}}_j)\right]_{M\times N}$ is accordingly partitioned into $n_z$ sub-matrices along columns,
\begin{equation} \label{eqn36}
\mathbf{K}_s=\left[\mathbf{K}_s^1, \mathbf{K}_s^2, \cdots, \mathbf{K}_s^{n_z} \right],\quad \mathrm{where}\quad  \mathbf{K}_s^h=\left[K_s(\mathbf{r}_i, \mathbf{\tilde{r}}_j^h)\right]_{M\times \tilde{M}}\,,\ 1\le h\le n_z\,.
\end{equation}

In (\ref{eqn36}), every sub-matrix $\mathbf{K}_s^h$ corresponds to the grid points in the same depth. We perform the singular value decomposition for $\mathbf{K}_s^h$,
\begin{equation} \label{eqn37}
\mathbf{K}_s^h=\mathbf{U}_s^h\mathbf{S}_s^h(\mathbf{V}_s^h)^T,\quad 1\le h\le n_z,\ s\in\{1,2,3\}\,,
\end{equation}
where $\mathbf{U}_s^h\in\mathbf{R}^{M\times M},\,\mathbf{V}_s^h\in\mathbf{R}^{\tilde{M}\times \tilde{M}}$ are unitary matrices, and
\begin{equation} \label{eqn38}
\mathbf{S}_s^h=\mathrm{diag}\left\{\sigma_1^{(s,h)},\sigma_2^{(s,h)},\cdots,\sigma_{\min\{M,\tilde{M}\}}^{(s,h)}\right\} \in \mathbf{R}^{M\times\tilde{M}}
\end{equation}
with the singular values $\sigma_1^{(s,h)}\ge\sigma_2^{(s,h)}\ge\cdots\ge\sigma_{\min\{M,\tilde{M}\}}^{(s,h)}\ge 0$. Then we truncate $\mathbf{S}_s^h$ by choosing a thresholding parameter $\epsilon_{svd}$ and setting
\begin{equation} \label{eqn39}
\tilde{\sigma}_i^{(s,h)}=\left\{
\begin{array}{lcr}
\sigma_i^{(s,h)}&,&\sigma_i^{(s,h)}\ge\epsilon_{svd}\\
0&,&\sigma_i^{(s,h)}<\epsilon_{svd}
\end{array}
\right. ,\qquad 1\le i\le \min\{M,\tilde{M}\}\,.
\end{equation}
Denoting $\tau_{sh}$ the number of non-zero $\tilde{\sigma}_i^{(s,h)}$, we have the truncated singular value decomposition for $\mathbf{K}_s^h$ in the following way,
\begin{equation} \label{eqn40}
\mathbf{K}_s^h\approx\tilde{\mathbf{U}}_s^h\tilde{\mathbf{S}}_s^h(\tilde{\mathbf{V}}_s^h)^T,\quad 1\le h\le n_z,\ s\in\{1,2,3\}\,,
\end{equation}
where $\tilde{\mathbf{S}}_s^h=\mathrm{diag}\left\{\tilde{\sigma}_1^{(s,h)}, \tilde{\sigma}_i^{(s,h)},\cdots,\tilde{\sigma}_{\tau_{sh}}^{(s,h)} \right\} \in \mathbf{R}^{\tau_{sh}\times\tau_{sh}}$, $\tilde{\mathbf{U}}_s^h\in\mathbf{R}^{M\times\tau_{sh}}$ is composed of the first $\tau_{sh}$ columns of $\mathbf{U}_s^h$, and $\tilde{\mathbf{V}}_s^h\in\mathbf{R}^{\tilde{M}\times\tau_{sh}}$ is composed of the first $\tau_{sh}$ columns of $\mathbf{V}_s^h$. The truncating process reduces the dimension of the matrices.

Finally we compute
\begin{equation} \label{eqn41}
\tilde{\mathbf{Q}}_s^h=\tilde{\mathbf{S}}_s^h(\tilde{\mathbf{V}}_s^h)^T\,,
\end{equation}
so that
\begin{equation} \label{eqn42}
\mathbf{K}_s^h\approx\tilde{\mathbf{U}}_s^h\tilde{\mathbf{Q}}_s^h,\quad 1\le h\le n_z,\ s\in\{1,2,3\}\,.
\end{equation}
The preconditioning is completed and the matrices $\tilde{\mathbf{U}}_s^h,\,\tilde{\mathbf{Q}}_s^h$ are stored for later use.

\subsubsection{Evaluation of the stochastic gradient}
Equation (\ref{eqn36}) and equation (\ref{eqn42}) imply the partitioned-truncated SVD for $\mathbf{K}^{S_k}_s$:
\begin{eqnarray}
\mathbf{K}_s^{S_k}&=&\left[\mathbf{K}_s^{1,S_k}, \mathbf{K}_s^{2,S_k}, \cdots, \mathbf{K}_s^{n_z, S_k} \right]\,, \label{eqn43} \\
\mathbf{K}_s^{h,S_k}&\approx&\tilde{\mathbf{U}}_s^{h,S_k}\tilde{\mathbf{Q}}_s^h,\qquad 1\le h\le n_z,\ s\in\{1,2,3\}\,, \label{eqn44}
\end{eqnarray}
where $\tilde{\mathbf{U}}_s^{h,S_k}$ is constructed by taking the rows of $\tilde{\mathbf{U}}_s^h$ with indices $i\in S_k$.

To evaluate the matrix multiplications (\ref{eqn33}) and (\ref{eqn34}), we perform partitions on $\mathbf{\kappa}$ and $\mathbf{D}$ in the same way as the partition of $\mathbf{K}_s^{S_k}$,
\begin{eqnarray}
\mathbf{\kappa}&=&[\mathbf{\kappa}_1, \mathbf{\kappa}_2, \cdots, \mathbf{\kappa}_{n_z}]\,, \label{eqn45} \\
\mathbf{D}&=&\mathrm{diag}\left\{\mathbf{D}_1,\mathbf{D}_2,\cdots,\mathbf{D}_{n_z}\right\}\,. \label{eqn46}
\end{eqnarray}
Finally, we propose the formulas of matrix multiplications for the stochastic gradient:
\begin{eqnarray}
\mathbf{I}_s^{S_k}&=&\sum_{h=1}^{n_z}\mathbf{K}_s^{h,S_k}\left(\mathbf{\kappa}_h\right)^T\approx \sum_{h=1}^{n_z}\tilde{\mathbf{U}}_s^{h,S_k}\left(\tilde{\mathbf{Q}}_s^h\left(\mathbf{\kappa}_h\right)^T\right) \,, \label{eqn47} \\
\mathbf{G}_s^{S_k}&=&\left[\mathbf{G}_s^{1,S_k}, \mathbf{G}_s^{2,S_k},\cdots, \mathbf{G}_s^{n_z,S_k}\right]\,,\quad \mathrm{with}\quad \mathbf{G}_s^{h,S_k}\approx\left(\mathbf{d}_s^{S_k}\tilde{\mathbf{U}}_s^{h,S_k}\right)\left(\tilde{\mathbf{Q}}_s^h\mathbf{D}_h\right)\,. \label{eqn48}
\end{eqnarray}
Realizing that $\mathbf{D}_h$ is a diagonal matrix with either $0$ or $1$ entries, $\tilde{\mathbf{Q}}_s^h\mathbf{D}_h$ in (\ref{eqn48}) can be simply obtained by taking the columns of $\tilde{\mathbf{Q}}_s^h$ with the same indices of the non-zero diagonal entries in $\mathbf{D}_h$ and setting the other columns to be $\mathbf{0}$.

\section{Results}
We include some test examples to illustrate the efficacy of the proposed algorithm. The measurement magnetic modulus data are generated using equations (\ref{eqn1})-(\ref{eqn3}). We discretize the magnetic integral (\ref{eqn1}) directly and we are not using the partitioned-truncated SVD when generating the measurement data. $5\%$ Gaussian noises are added to each component $B_s$,
\[
B_s^*=B_s\left(1+\eta\cdot\mathcal{N}(0,1)\right)\,,
\]
where $\eta=5\%$ and $\mathcal{N}(0,1)$ denotes the Gaussian noises with $0$ mean and standard deviation $1$. Then the synthetic nonlinear modulus data are generated according to equation (\ref{eqn3}),
\[
d^*(\mathbf{r})=\left(\sum_{s=1}^3B_s^*(\mathbf{r})^2\right)^{\frac{1}{2}}\,.
\]
In the test examples, the inducing magnetic field has the strength of $B^0=5.95\times 10^4$\,nT, and the direction $\hat{\mathbf{B}}^0$ is prescribed by inclination $I^0$ and declination $D^0$ so that
\[
\hat{\mathbf{B}}^0=(\cos I^0\sin D^0, \cos I^0 \cos D^0,-\sin I^0)\,,
\]
where we designate the direction of $x$-axis as east, the direction of $y$-axis as north, and the direction of $z$-axis as negative depth. The average susceptibility value in the level-set formulation (\ref{eqn5}) is given as $\kappa_0=0.05$. In the inversion algorithm, we use the partitioned-truncated SVD when evaluating the predicted data and the stochastic gradient, where the thresholding parameter for truncating is taken as $\epsilon_{svd}=10^{-5}$.

\subsection{Examples in a shallow region}
\begin{figure}
\centering
\subfigure[]{\scalebox{0.4}[0.4]{\includegraphics{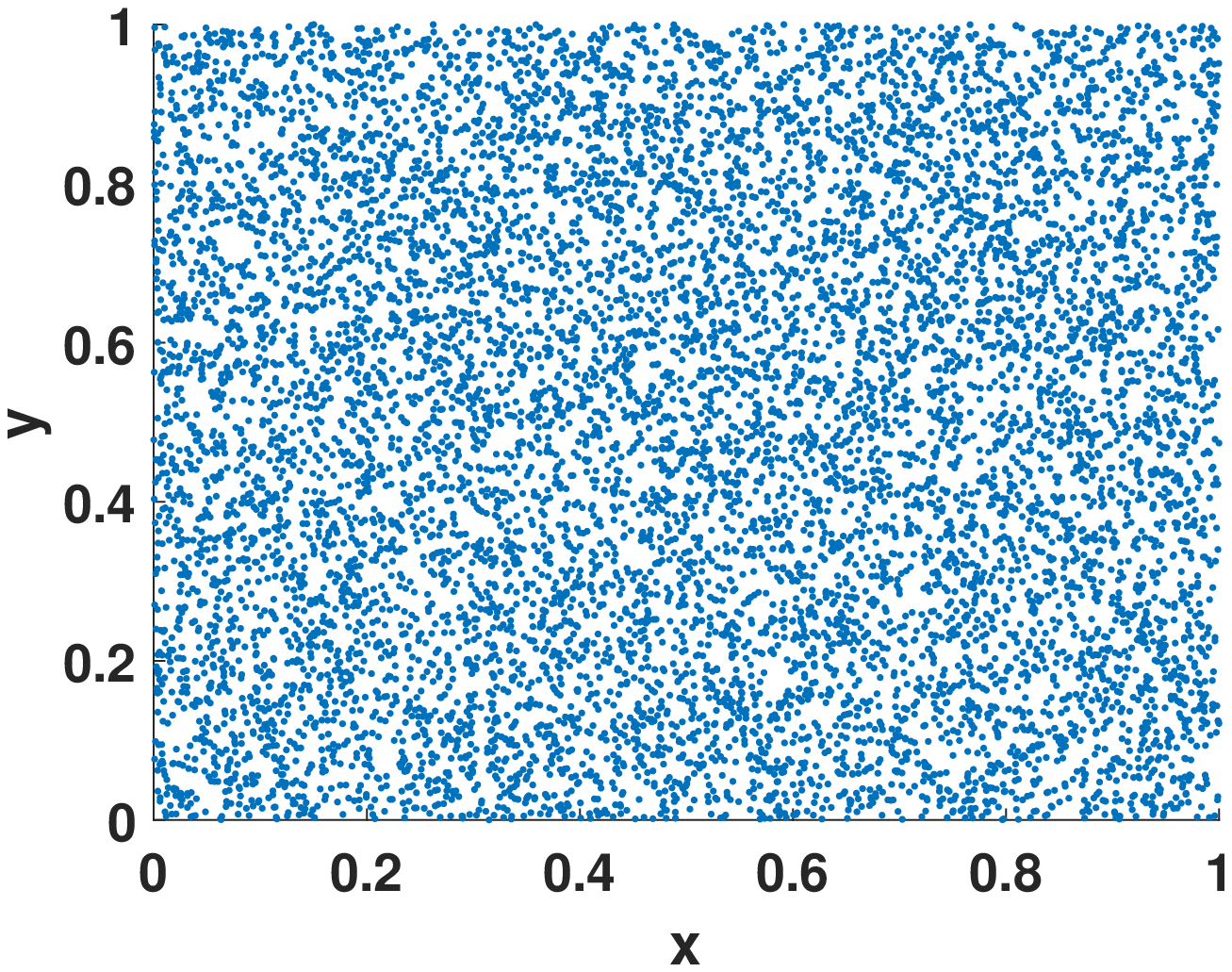}}}
\subfigure[]{\scalebox{0.4}[0.4]{\includegraphics{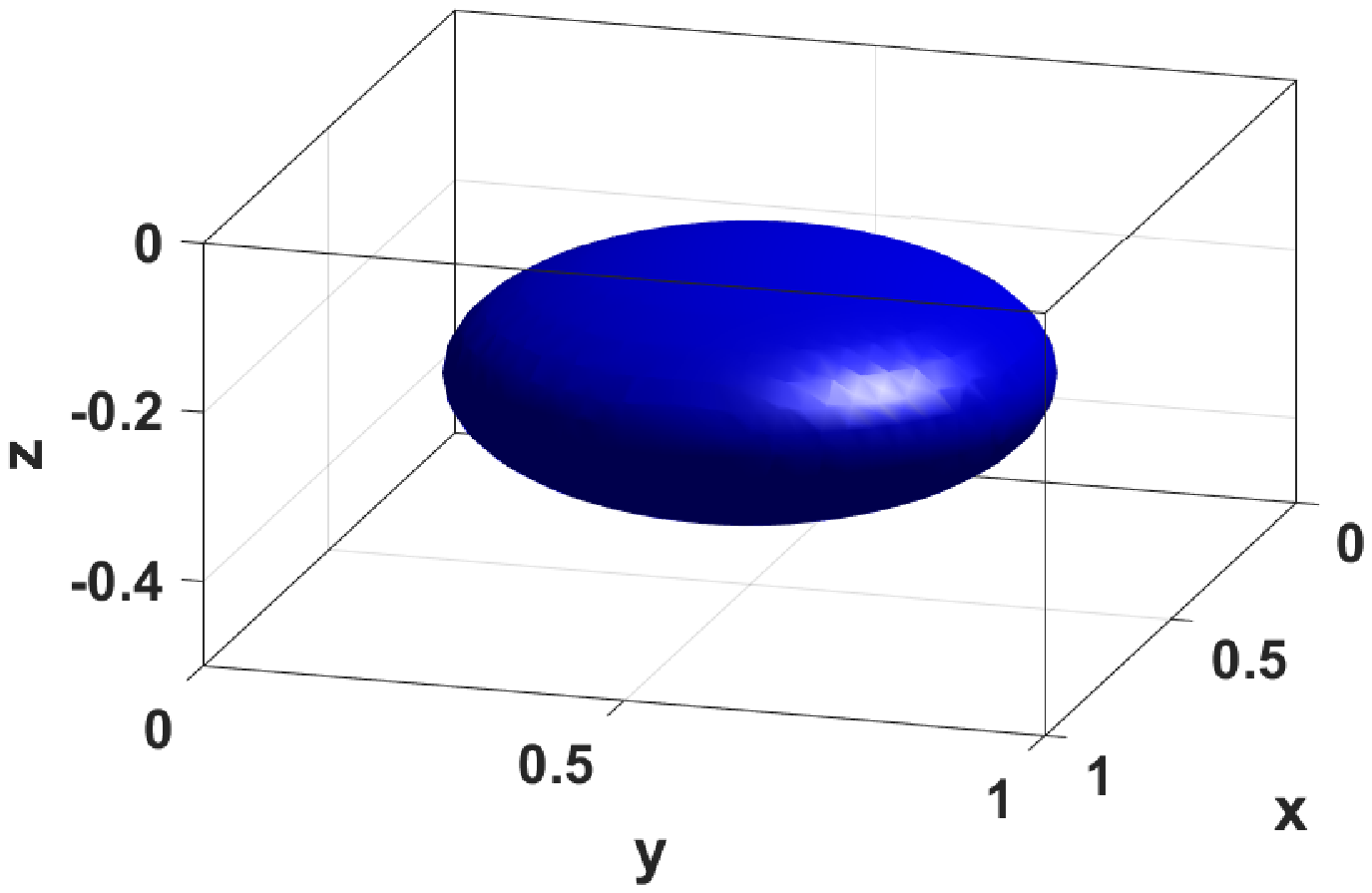}}}
\caption{Setup of examples in the shallow region. (a) Measurement points along $\Gamma=[0\ 1]\times[0\ 1]\times\{z=0.1\}\,km$; (b) initial guess.}
\label{Fig1}
\end{figure}

Magnetic susceptibilities are reconstructed in the computational domain $\Omega=[0\ 1]\times[0\ 1]\times[-0.5\ 0]\,km$, and we discretize $\Omega$ uniformly into $41\times41\times21$ mesh grids. The measurement boundary is taken as $\Gamma=[0\ 1]\times[0\ 1]\times\{z=0.1\}\,km$, along which 10,000 measurement points are randomly distributed; Figure \ref{Fig1}\,(a) shows the picture of measurement points along $\Gamma$. In the inversion algorithm, we set the mini-batch size to be $b=200$ for the stochastic gradient descent. The initial guess of the level-set function $\phi$ is set as follows,
\[
\phi_{initial}=1-\sqrt{\frac{(x-0.5)^2}{0.35^2}+\frac{(y-0.5)^2}{0.35^2}+\frac{(z-0.25)^2}{0.15^2}}\,;
\]
the zero level-set of the initial structure is an ellipsoid as shown in Figure \ref{Fig1}\,(b). 

\subsubsection{Example 1}
\begin{figure}
\centering
\subfigure[]{\scalebox{0.4}[0.4]{\includegraphics{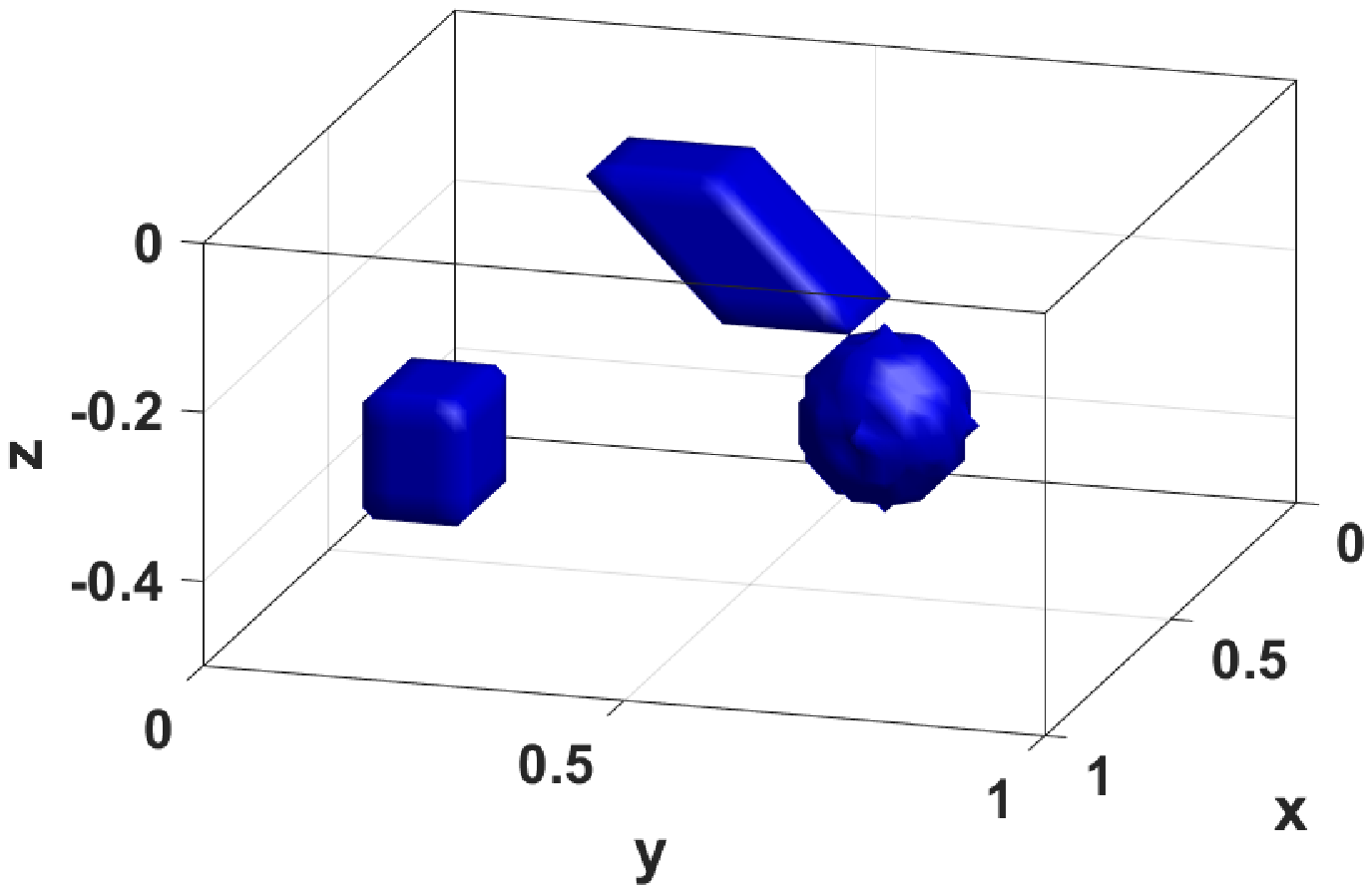}}}
\subfigure[]{\scalebox{0.4}[0.4]{\includegraphics{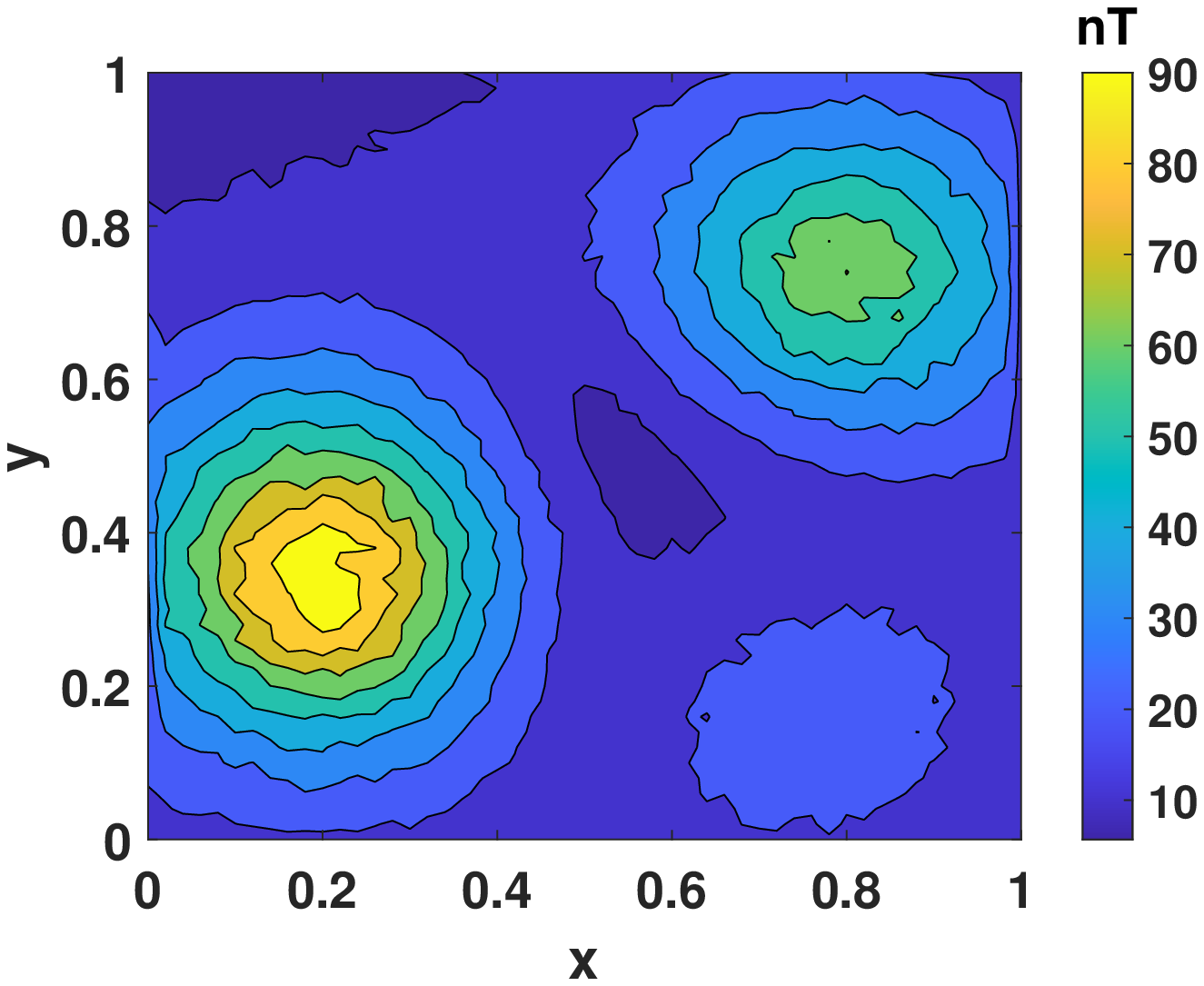}}}
\caption{Example 1. (a) True model; (b) magnetic modulus data with Gaussian noises.}
\label{Fig2}
\end{figure}

\begin{figure}
\centering
\subfigure[]{\scalebox{0.4}[0.4]{\includegraphics{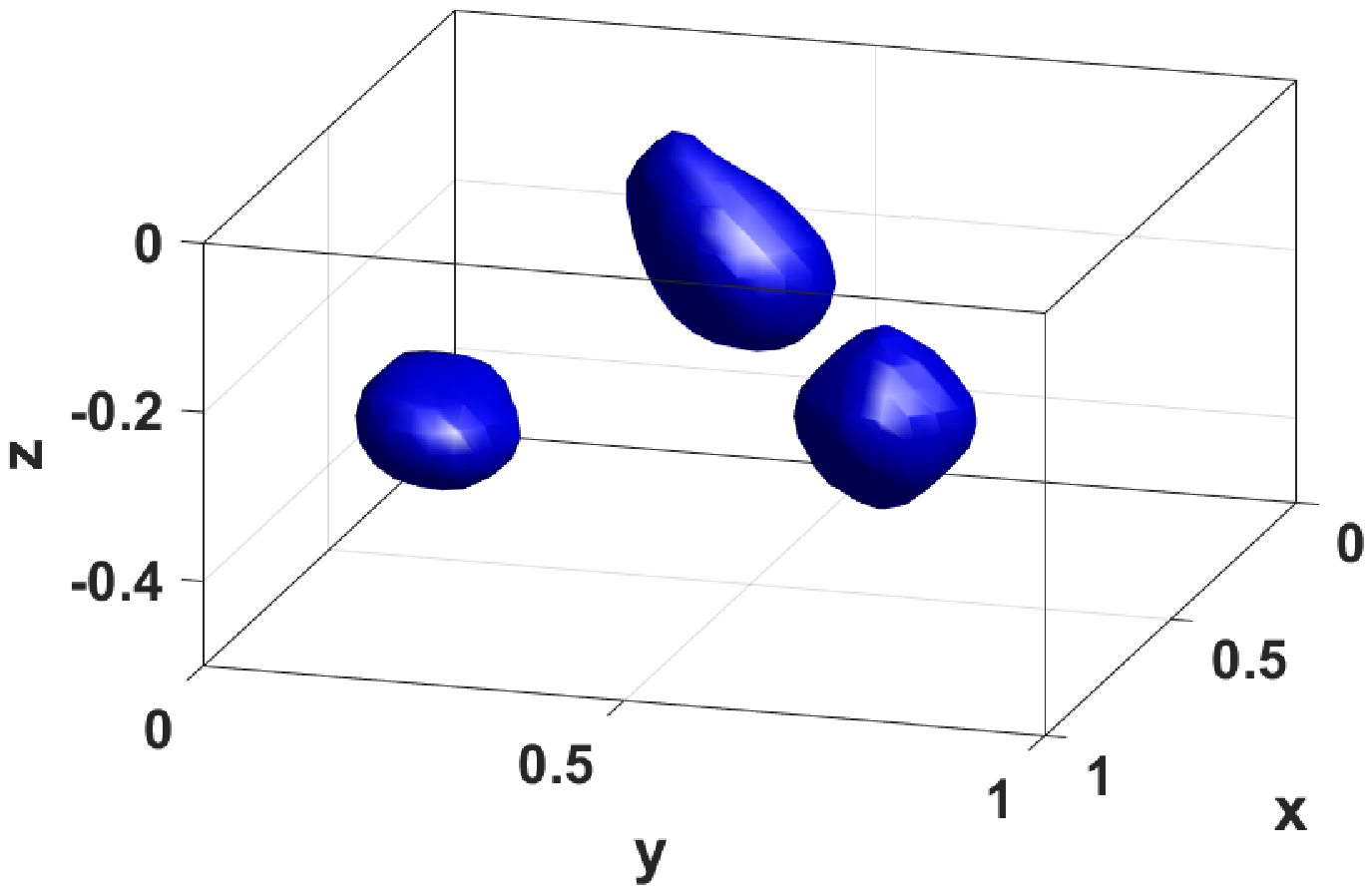}}}
\subfigure[]{\scalebox{0.4}[0.4]{\includegraphics{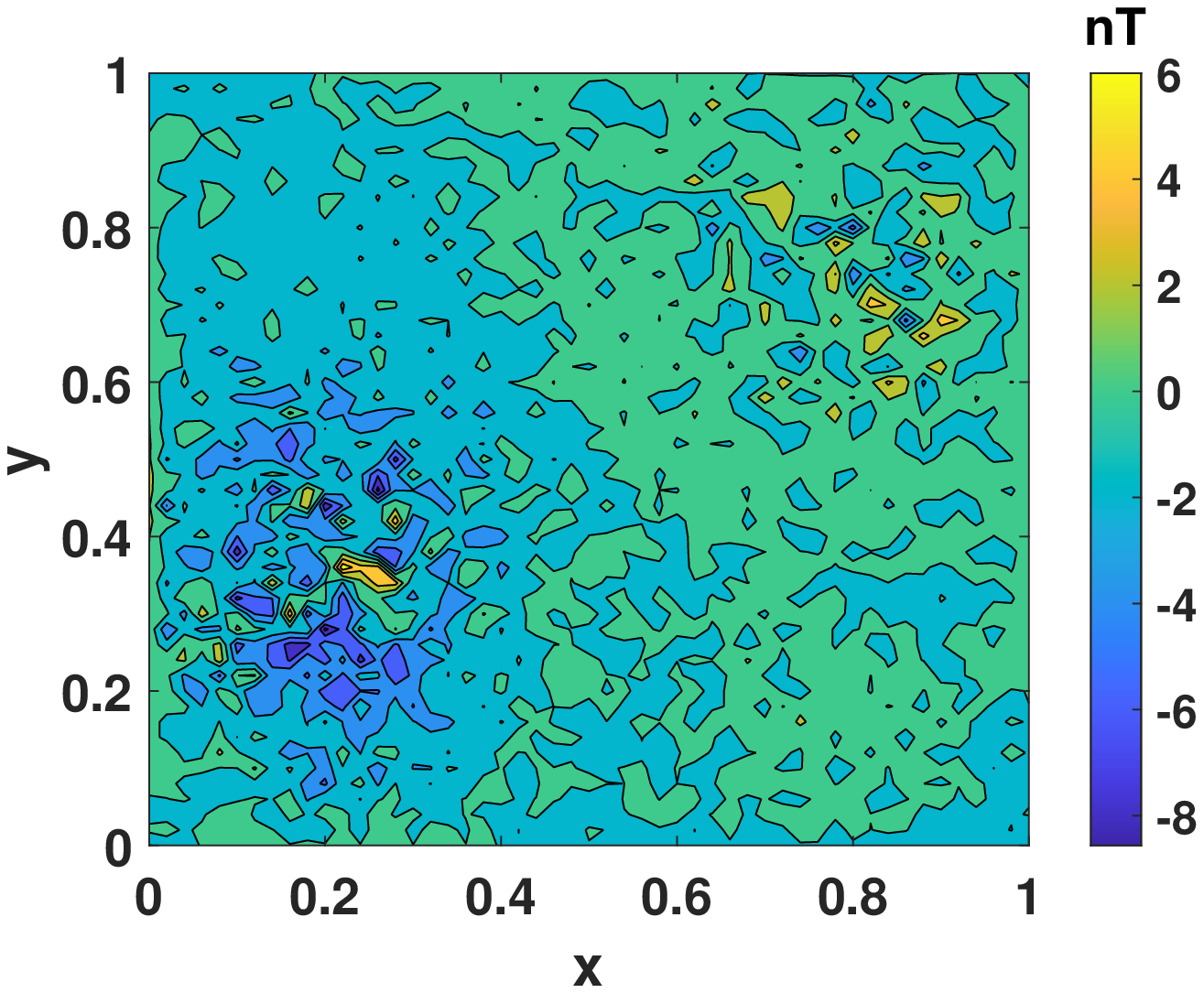}}}
\subfigure[]{\scalebox{0.4}[0.4]{\includegraphics{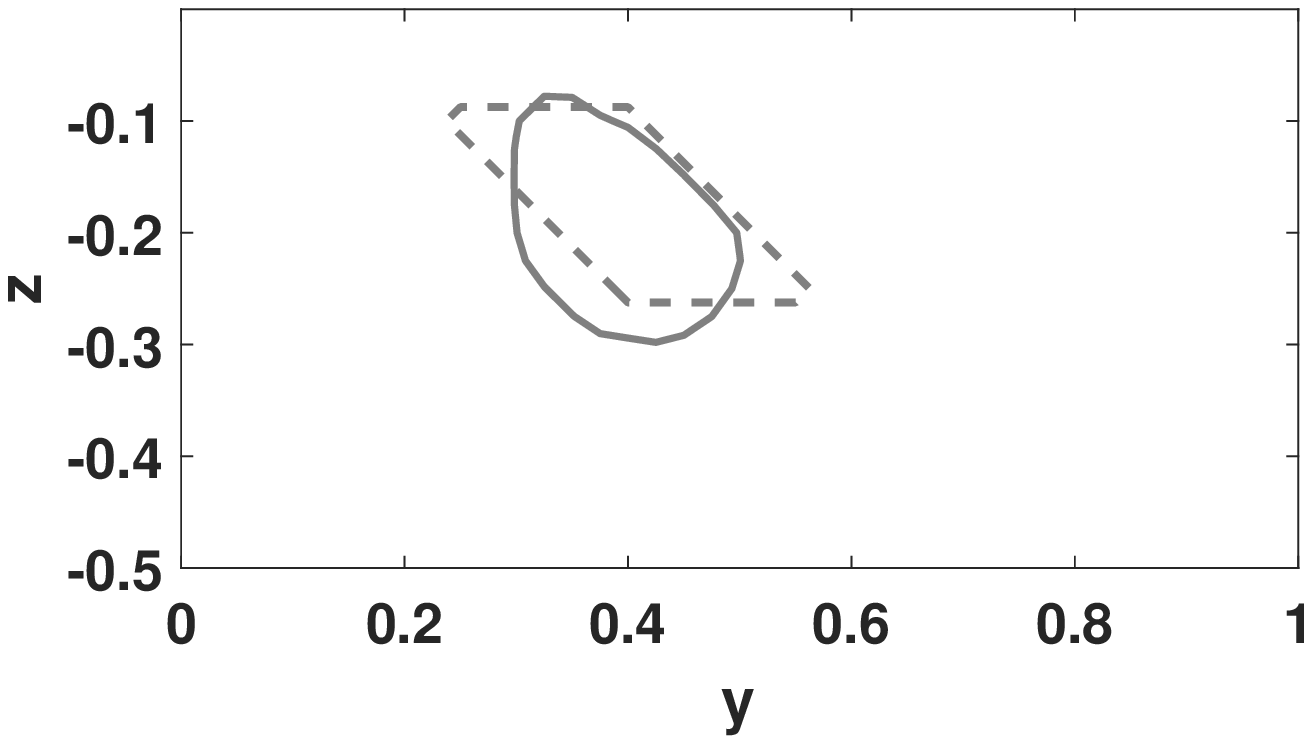}}}
\subfigure[]{\scalebox{0.4}[0.4]{\includegraphics{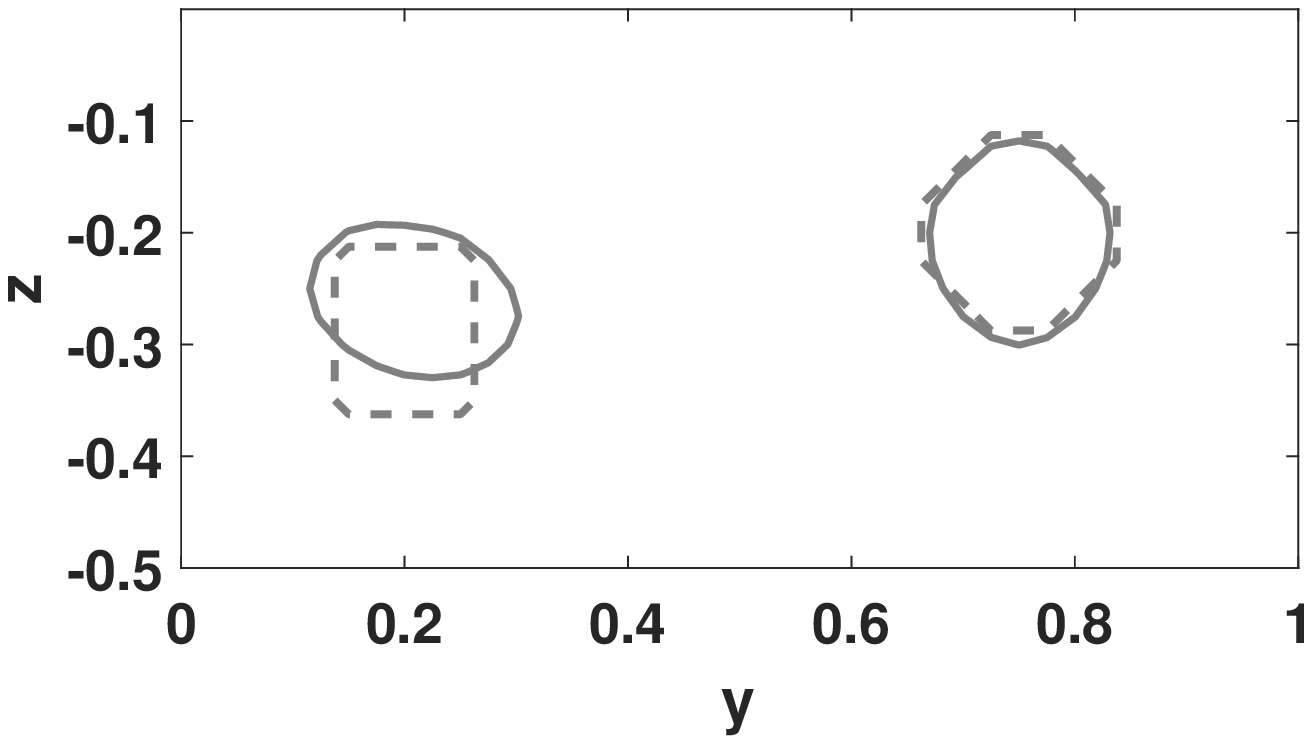}}}
\caption{Example 1. Inversion results. (a) recovered solution; (b) data discrepancy $d-d^*$; (c)-(d) cross-sections along $x=0.25$ and $x=0.75$, respectively, where the dashed line indicates the true model and the solid line indicates the recovered solution.}
\label{Fig3}
\end{figure}

\begin{figure}
\centering
\subfigure[]{\scalebox{0.4}[0.4]{\includegraphics{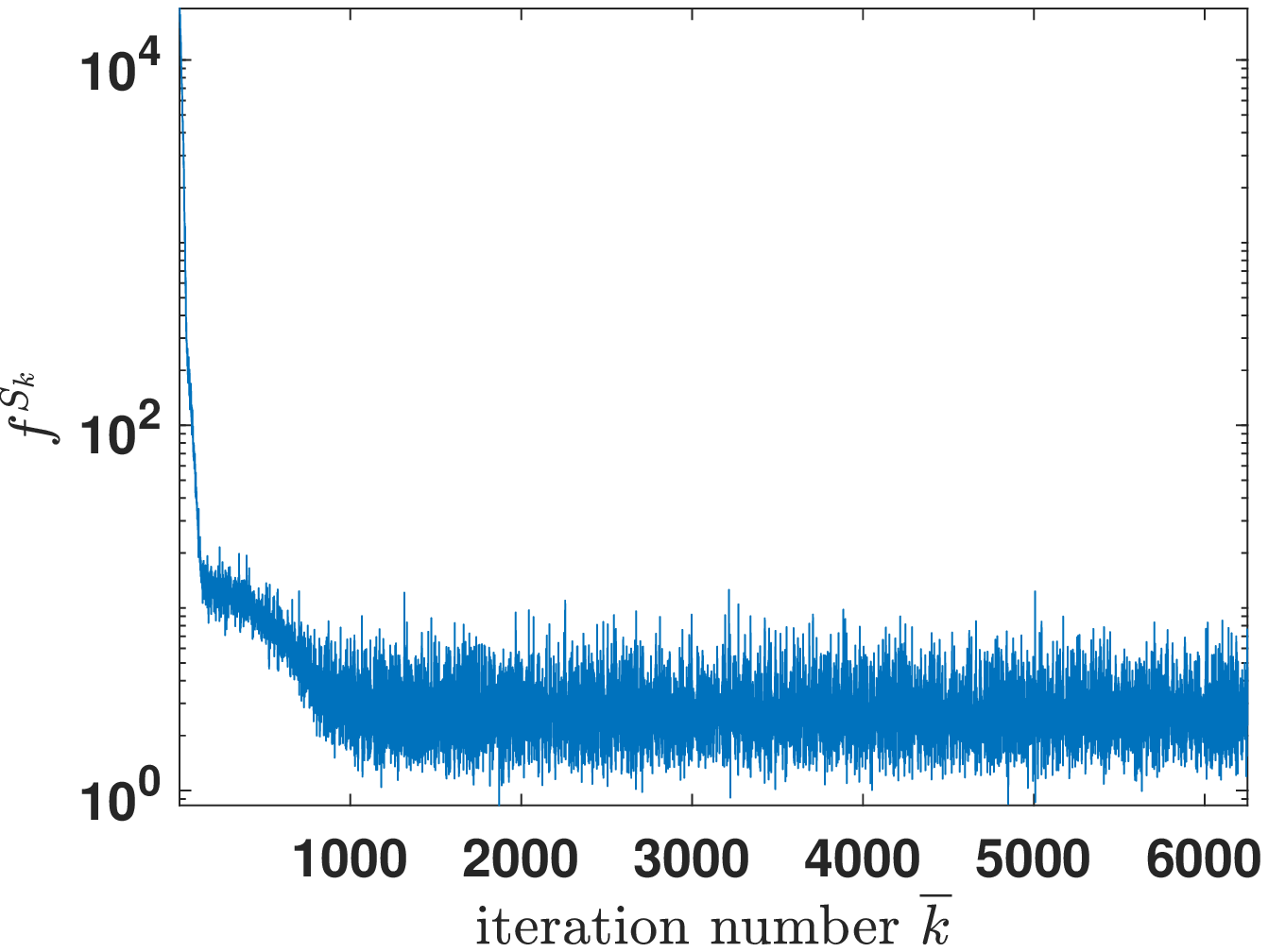}}}
\subfigure[]{\scalebox{0.4}[0.4]{\includegraphics{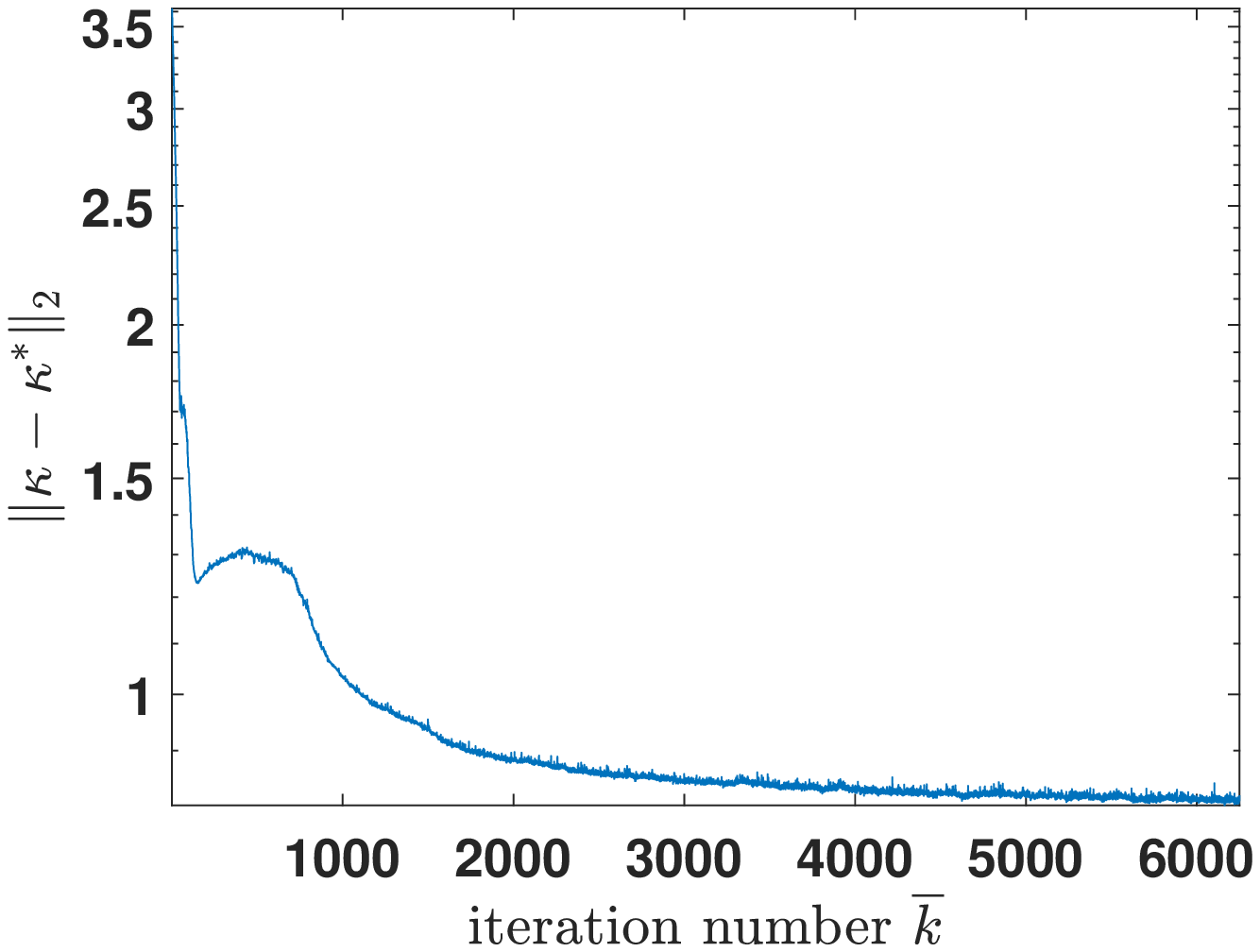}}}
\caption{Example 1. Performance of the inversion algorithm. (a) Mini-batch misfit function $f^{S_k}$; (b) $l^2$-difference between the recovered susceptibility $\kappa$ and the true susceptibility $\kappa^*$: $\|\kappa-\kappa^*\|_2$. The misfit function $f^{S_k}$ is oscillating partly due to the nonlinearity of level-set inversion and partly due to the noisy nature of stochastic gradient descent. On the other hand, the performance of the $l^2$-difference $\|\kappa-\kappa^*\|_2$ implies a smooth convergence.}
\label{Fig4}
\end{figure}

As shown in Figure \ref{Fig2}\,(a), the true model of the volume susceptibility consists of 3 magnetic sources, where the cube is centered at $(0.75,0.2,-0.3)$ with side length $0.15$, and the sphere is centered at $(0.8, 0.75, -0.2)$ with radius $0.1$. We assume an inducing field with inclination and declination $(I^0,D^0)=(75^{\circ},25^{\circ})$. Figure \ref{Fig2}\,(b) plots the magnetic modulus data on the measurement surface. Figure \ref{Fig3} provides the inversion results. Figure \ref{Fig3}\,(a) shows the recovered solution using the proposed stochastic gradient descent algorithm with partitioned-truncated SVD; Figure \ref{Fig3}\,(b) plots the data discrepancy $d-d^*$; Figures \ref{Fig3}\,(c) and \ref{Fig3}\,(d) plot the cross-sections along $x=0.25$ and $x=0.75$, respectively, where the dashed line indicates the true model and the solid line indicates the recovered solution. The solution adequately reproduces the magnetic modulus data, although we employ the partitioned-truncated SVD to approximate matrix multiplications and we never use all the 10,000 measurement data at once in the inversion algorithm. The recovered susceptibility structure matches well with the true model, where both the shape and location of the magnetic sources are successfully recovered. In Figure \ref{Fig4},  we provide more information to illustrate the performance of the inversion algorithm. Figure \ref{Fig4}\,(a) presents the evolution of the mini-batch misfit function $f^{S_k}$ as defined in equation (\ref{eqn20}); Figure \ref{Fig4}\,(b) plots $l^2$-difference between the recovered susceptibility $\kappa$ and the true susceptibility $\kappa^*$. The misfit function $f^{S_k}$ is oscillating partly due to the nonlinearity of level-set inversion and partly due to the noisy nature of stochastic gradient descent. On the other hand, the performance of the $l^2$-difference $\|\kappa-\kappa^*\|_2$ implies a smooth convergence.

\subsubsection{Example 2}
\begin{figure}
\centering
\subfigure[]{\scalebox{0.4}[0.4]{\includegraphics{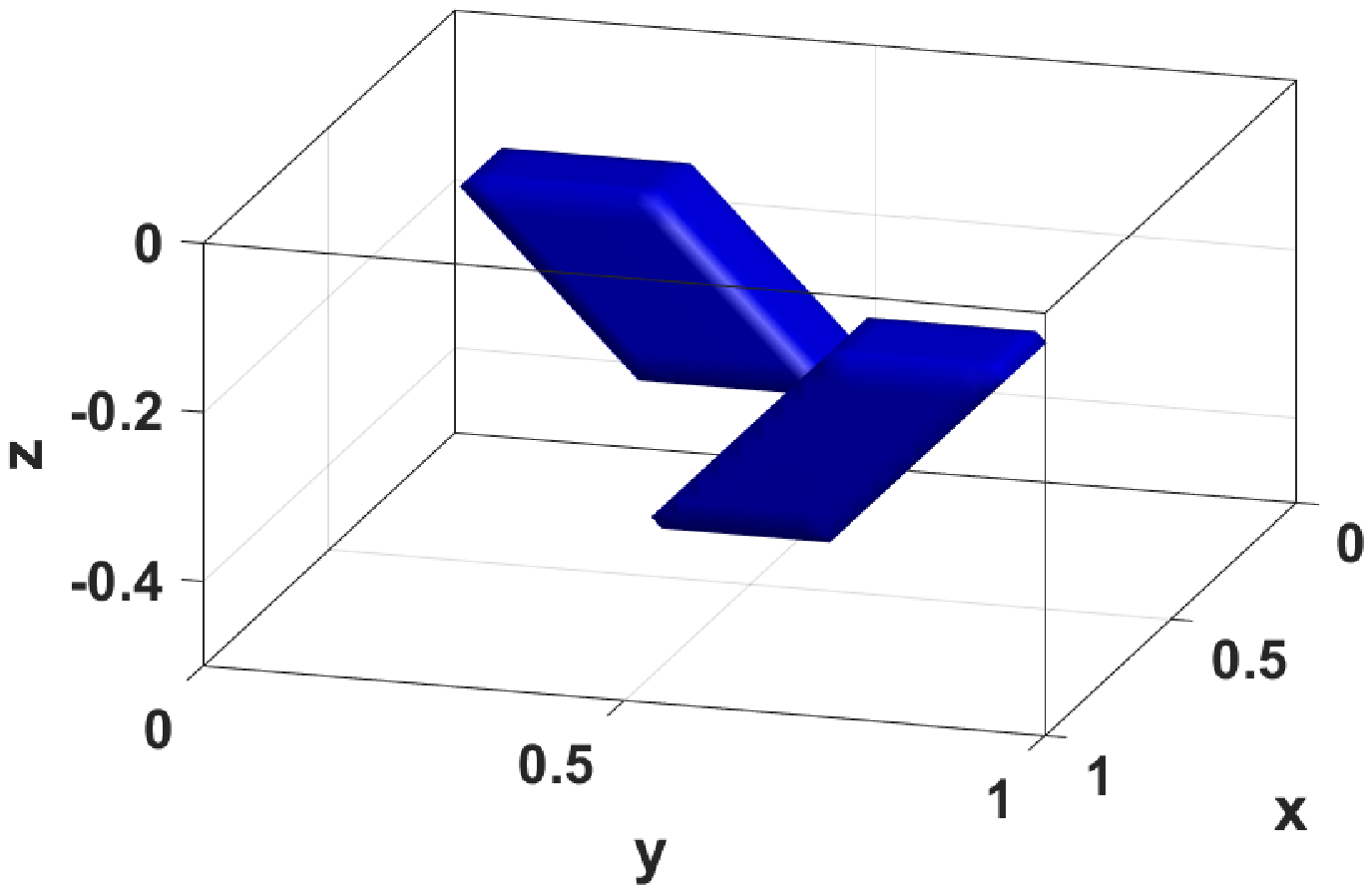}}}
\subfigure[]{\scalebox{0.4}[0.4]{\includegraphics{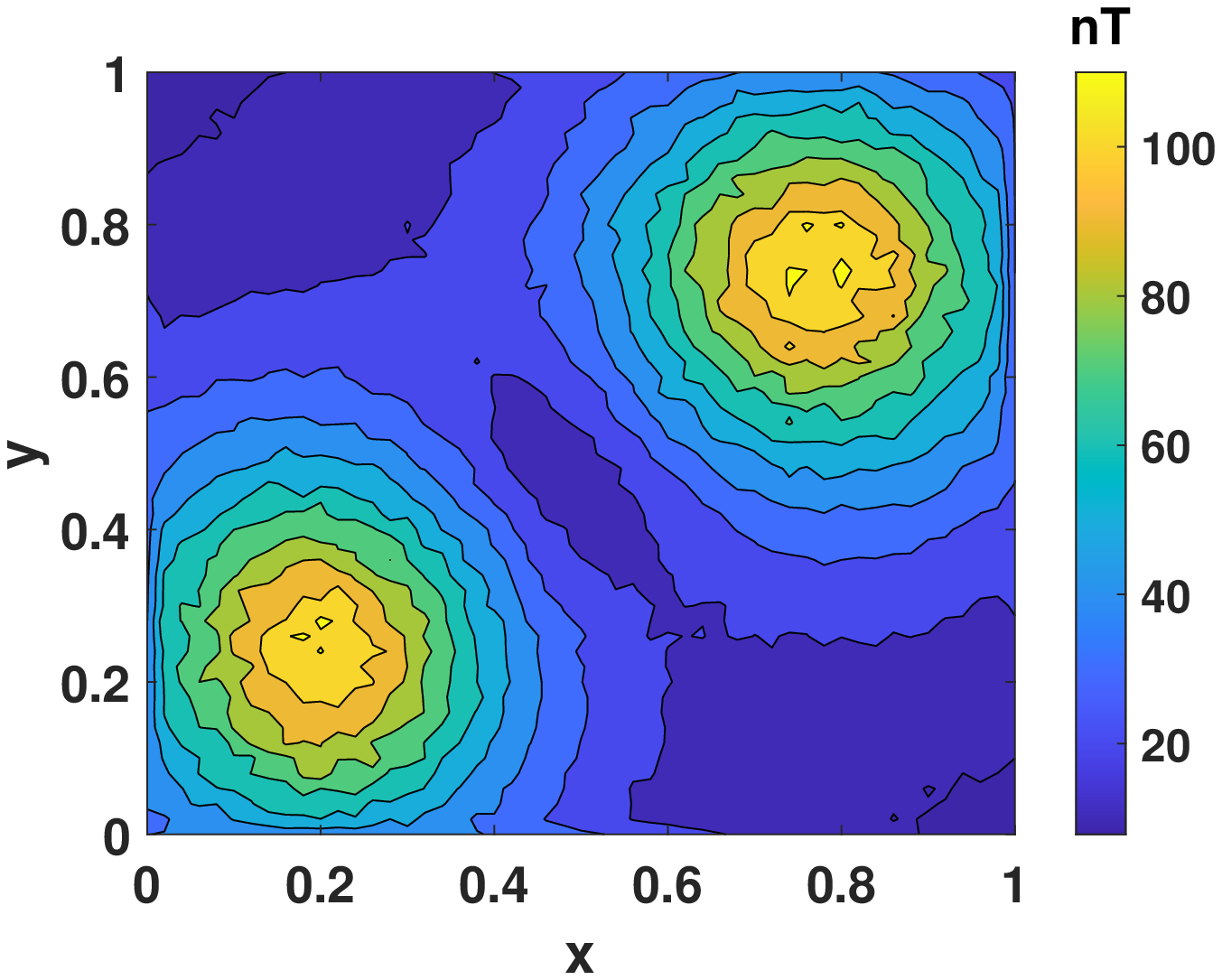}}}
\caption{Example 2. (a) True model; (b) magnetic modulus data with Gaussian noises.}
\label{Fig5}
\end{figure}

\begin{figure}
\centering
\subfigure[]{\scalebox{0.4}[0.4]{\includegraphics{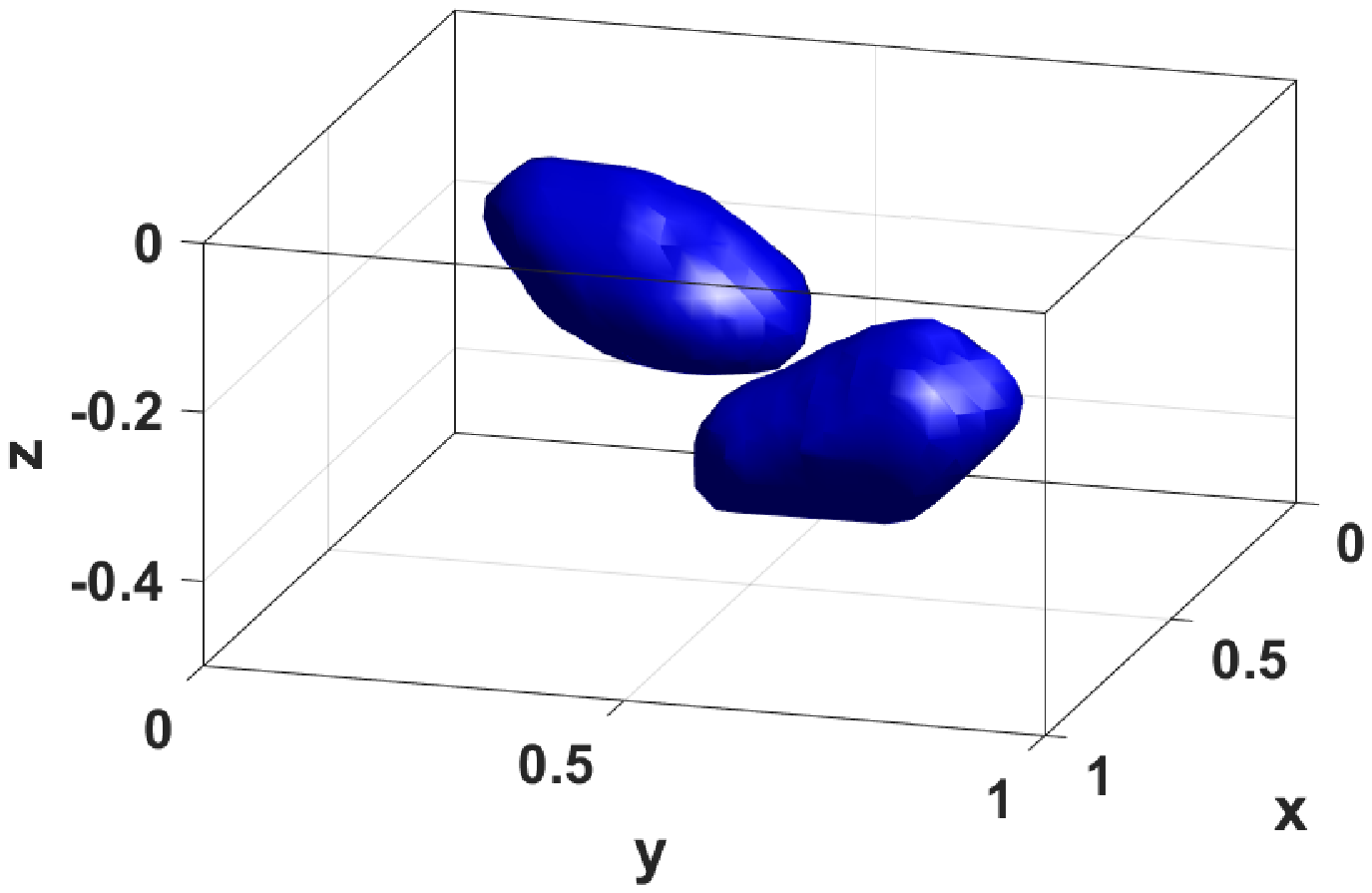}}}
\subfigure[]{\scalebox{0.4}[0.4]{\includegraphics{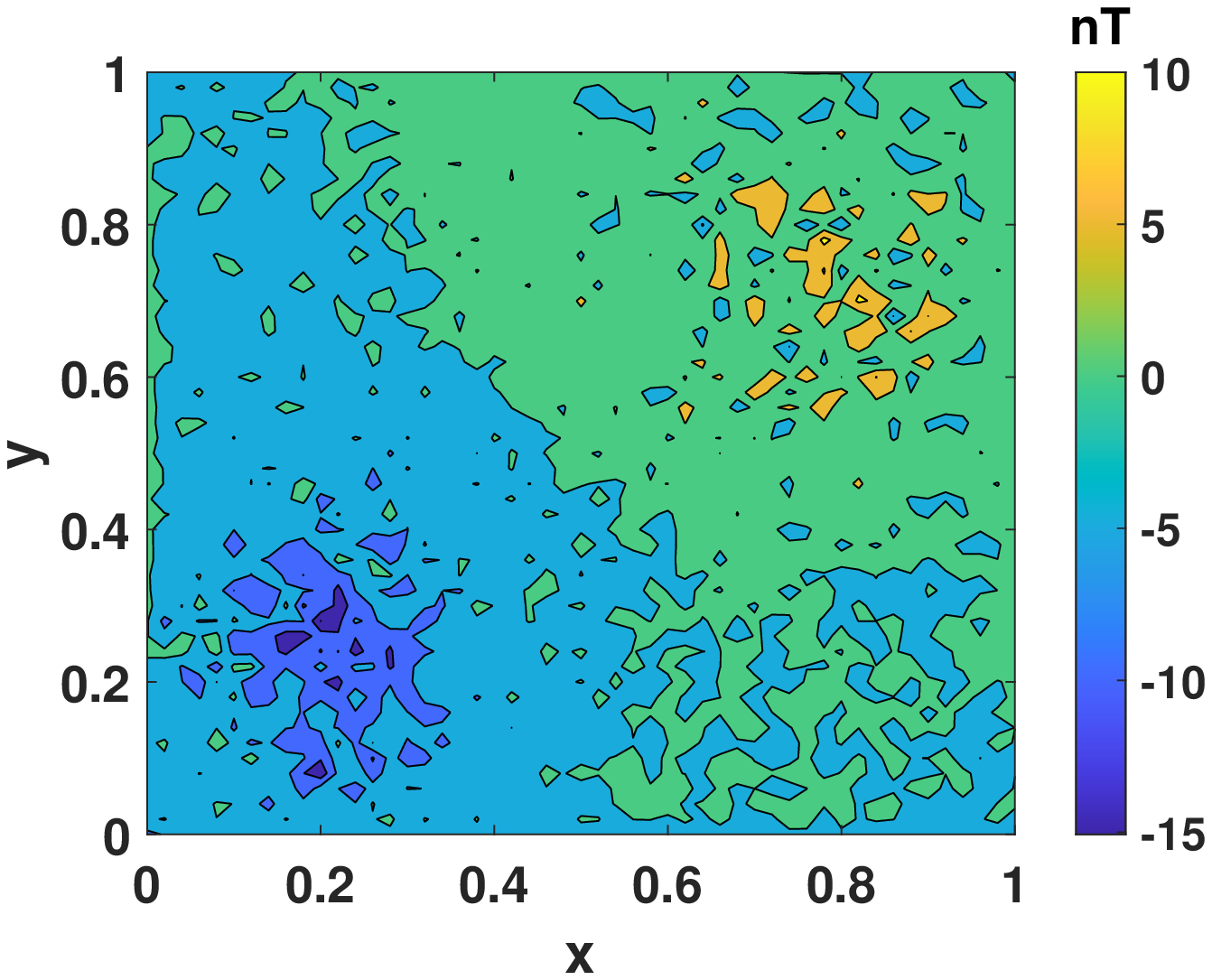}}}
\subfigure[]{\scalebox{0.4}[0.4]{\includegraphics{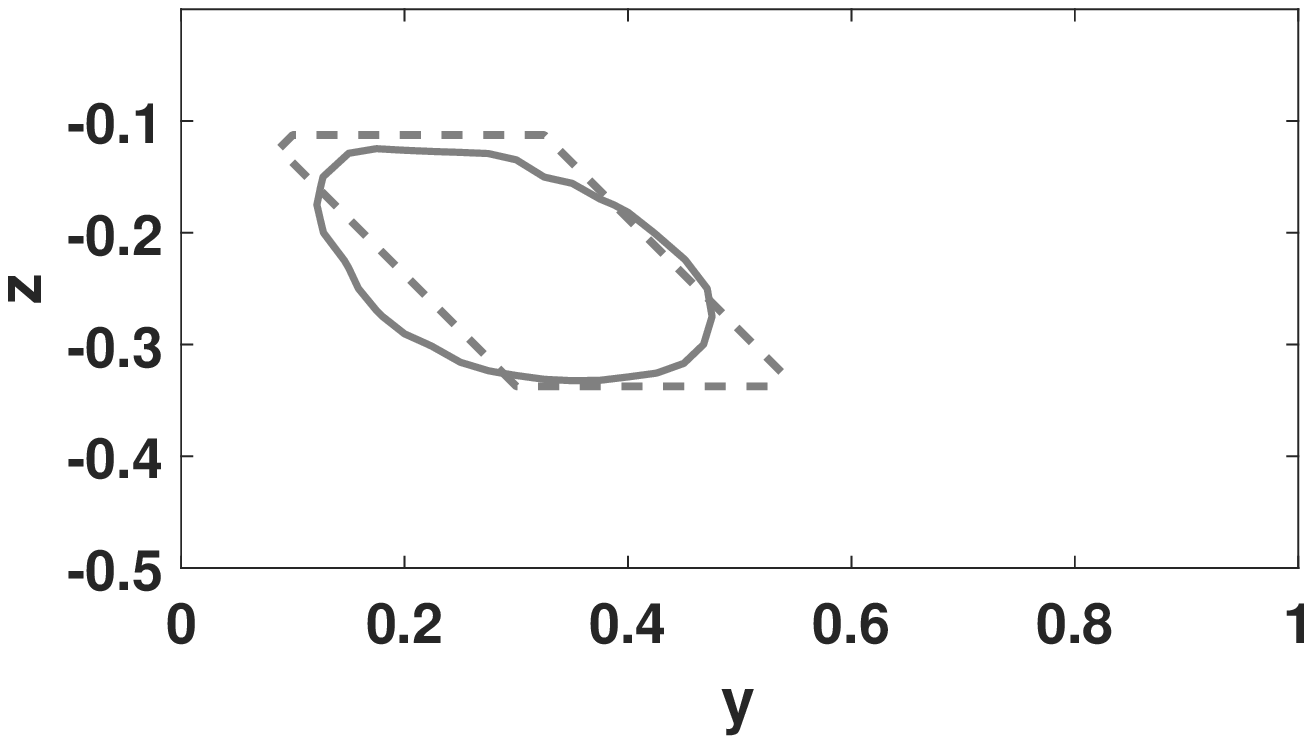}}}
\subfigure[]{\scalebox{0.4}[0.4]{\includegraphics{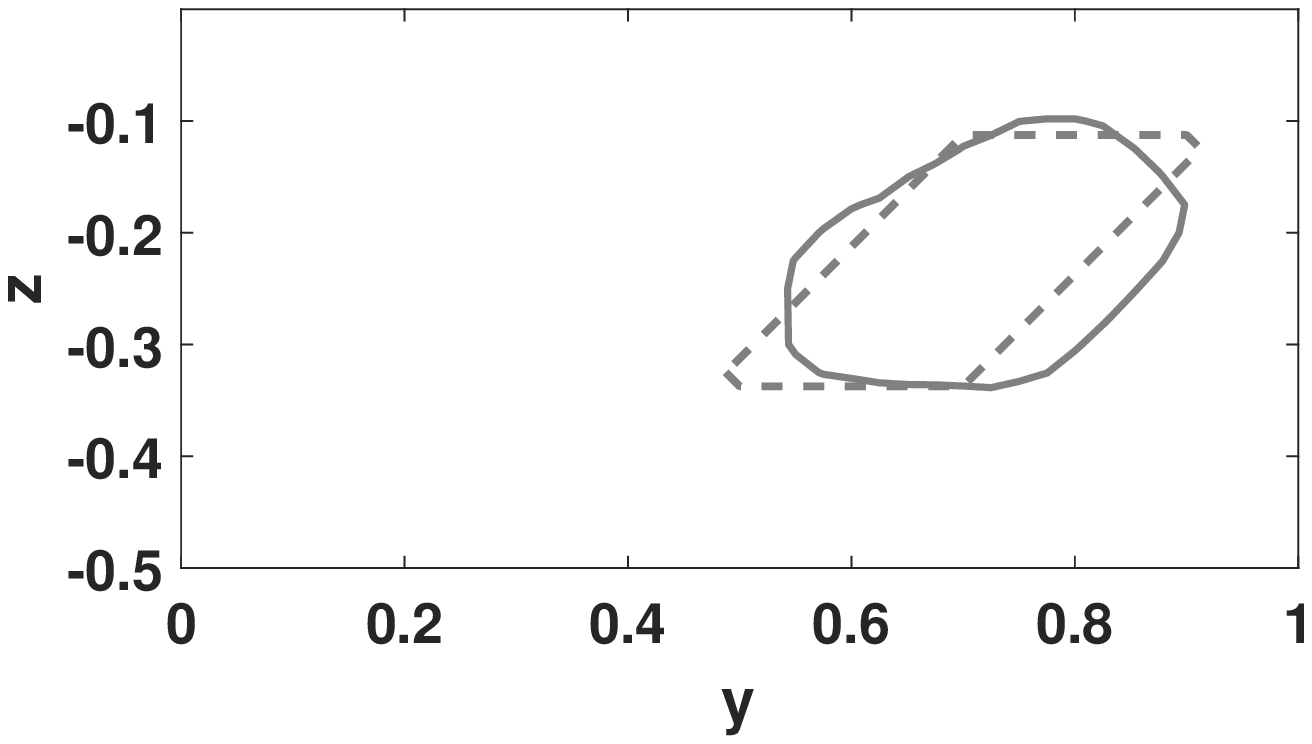}}}
\caption{Example 2. Inversion results. (a) recovered solution; (b) data discrepancy $d-d^*$; (c)-(d) cross-sections along $x=0.25$ and $x=0.75$, respectively, where the dashed line indicates the true model and the solid line indicates the recovered solution.}
\label{Fig6}
\end{figure}

\begin{figure}
\centering
\subfigure[]{\scalebox{0.4}[0.4]{\includegraphics{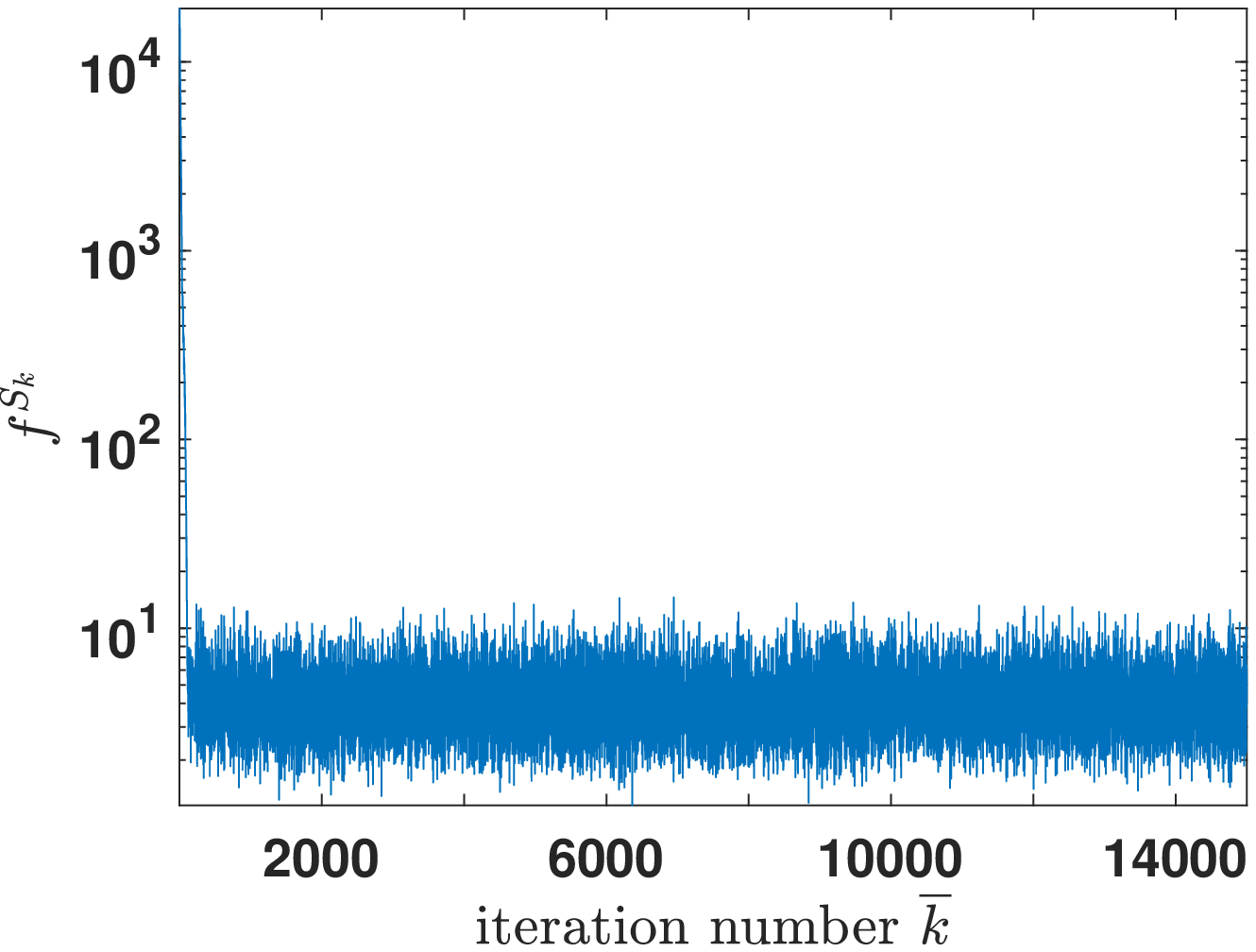}}}
\subfigure[]{\scalebox{0.4}[0.4]{\includegraphics{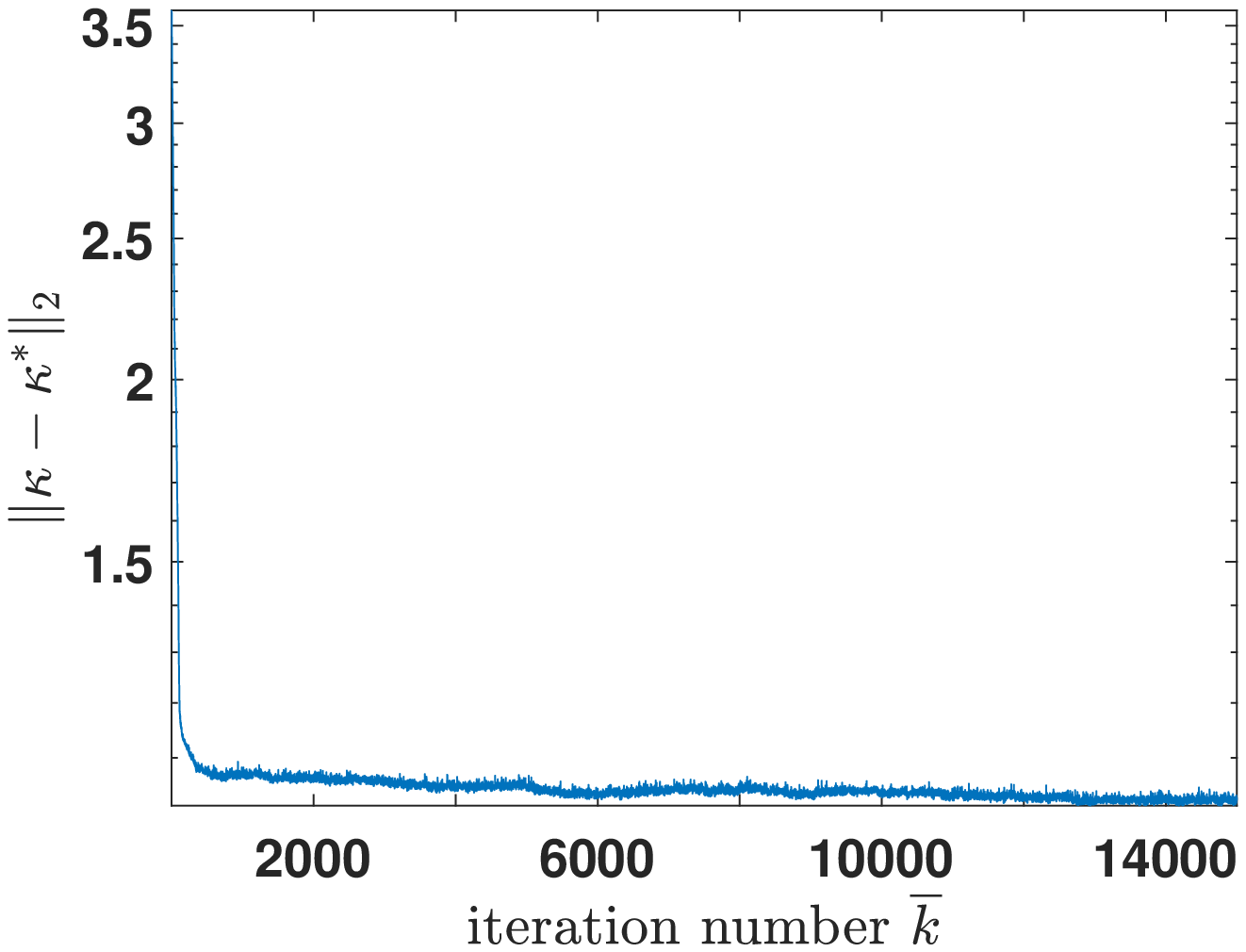}}}
\caption{Example 2. Performance of the inversion algorithm. (a) Mini-batch misfit function $f^{S_k}$; (b) $\|\kappa-\kappa^*\|_2$.}
\label{Fig7}
\end{figure}

Figure \ref{Fig5}\,(a) shows the true model of the volume susceptibility, where two dipping prisms are distributed in the computational domain. The inducing field has the inclination and declination $(I^0,D^0)=(75^{\circ},25^{\circ})$. Figure \ref{Fig5}\,(b) plots the magnetic modulus data on the measurement surface. Figure \ref{Fig6} provides the inversion results, where Figure \ref{Fig6}\,(a) shows the recovered solution, Figure \ref{Fig6}\,(b) plots the data discrepancy $d-d^*$, and Figures \ref{Fig6}\,(c),\,\ref{Fig6}\,(d) plot the cross-sections along $x=0.25$ and $x=0.75$, respectively. The magnetic modulus data are adequately reproduced. The recovered solution correctly captures the shape and location of the magnetic prisms, where the depth resolution looks perfect in the displayed cross-sections. In Figure \ref{Fig7}, we present the mini-batch misfit function $f^{S_k}$ and the $l^2$-difference $\|\kappa-\kappa^*\|_2$. The performance of $\|\kappa-\kappa^*\|_2$ shows a smooth convergence although the mini-batch misfit function $f^{S_k}$ is oscillating.

\subsubsection{Example 3}
\begin{figure}
\centering
\subfigure[]{\scalebox{0.4}[0.4]{\includegraphics{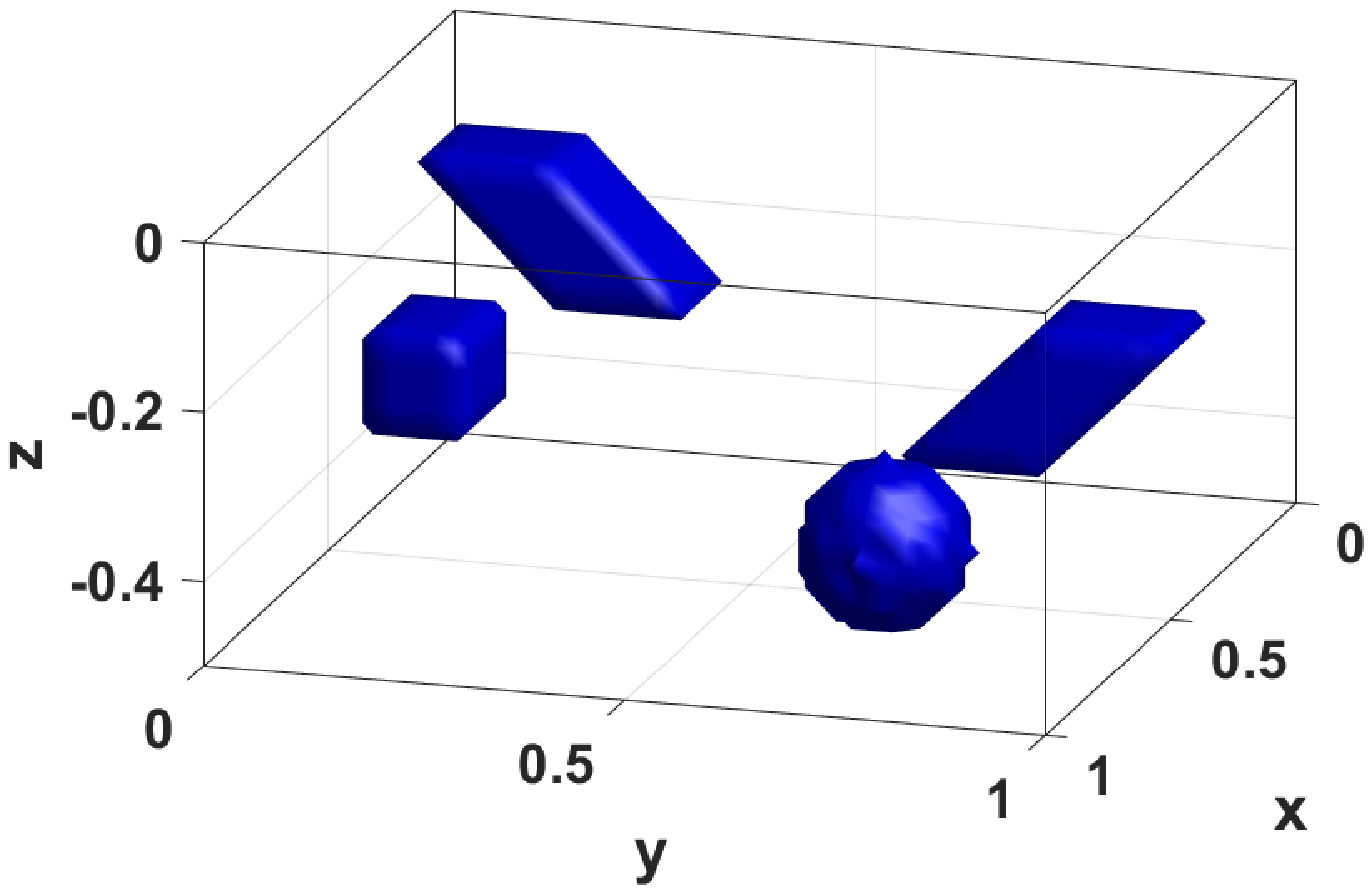}}}
\subfigure[]{\scalebox{0.4}[0.4]{\includegraphics{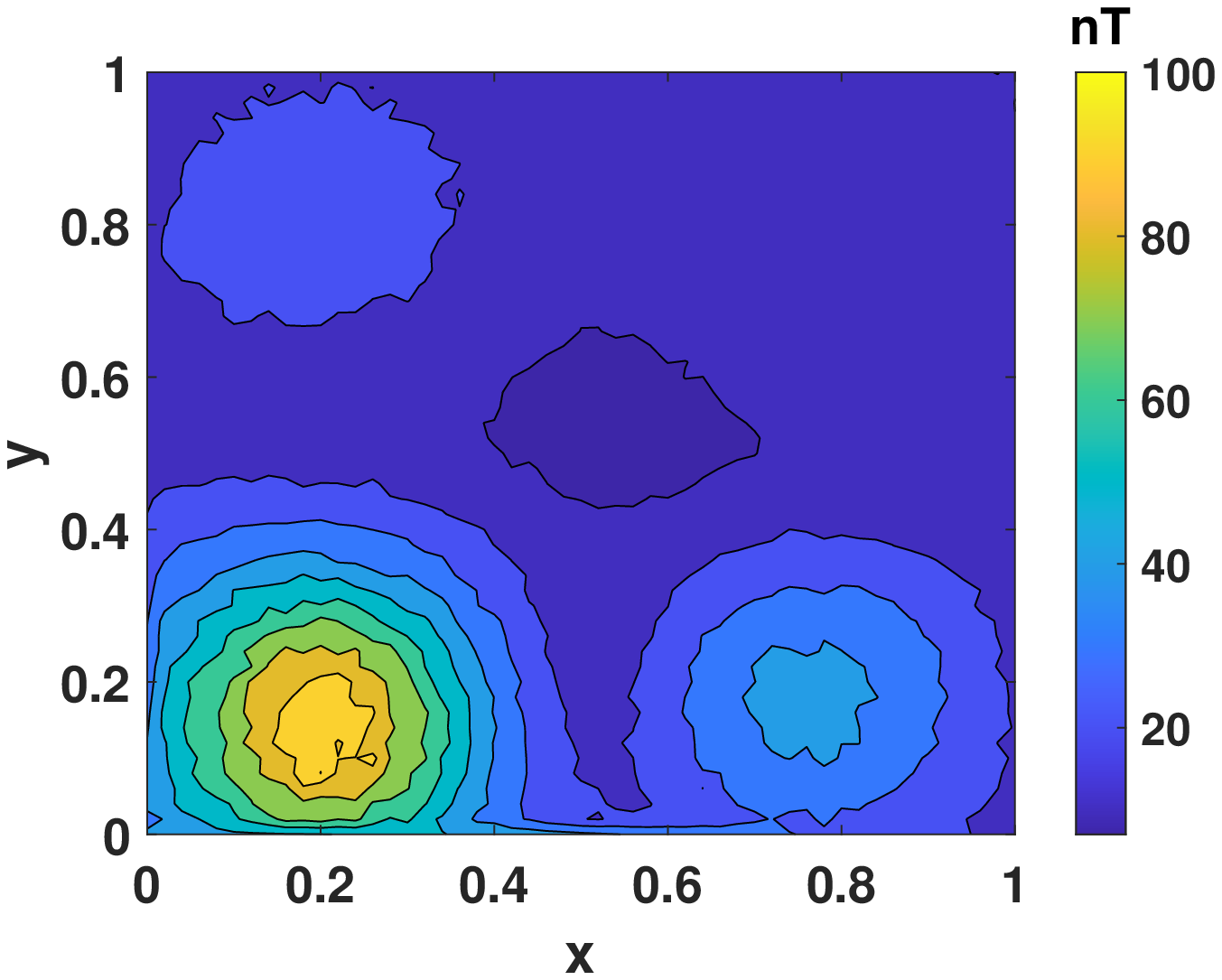}}}
\caption{Example 3. (a) True model; (b) magnetic modulus data with Gaussian noises.}
\label{Fig8}
\end{figure}

\begin{figure}
\centering
\subfigure[]{\scalebox{0.4}[0.4]{\includegraphics{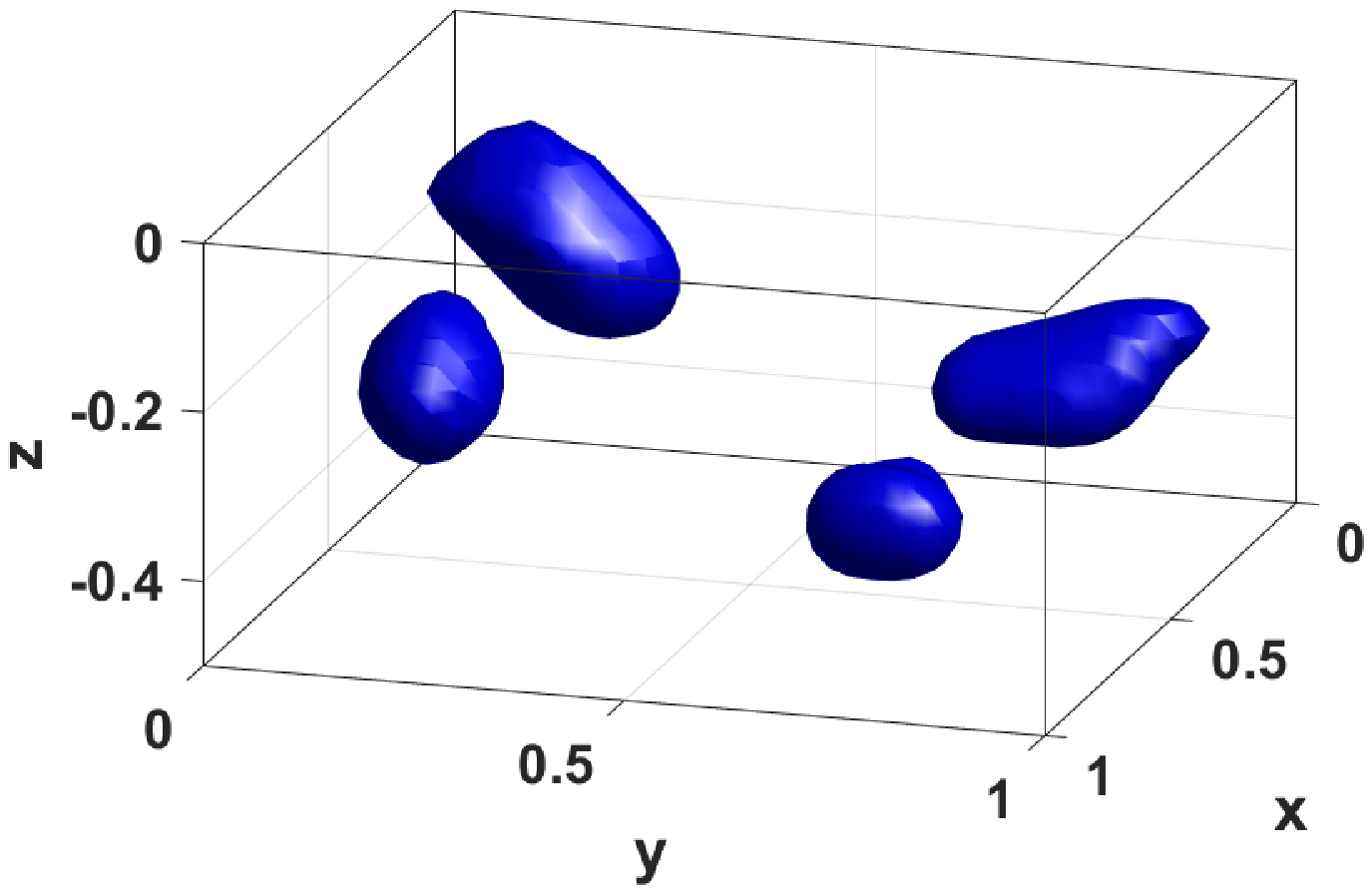}}}
\subfigure[]{\scalebox{0.4}[0.4]{\includegraphics{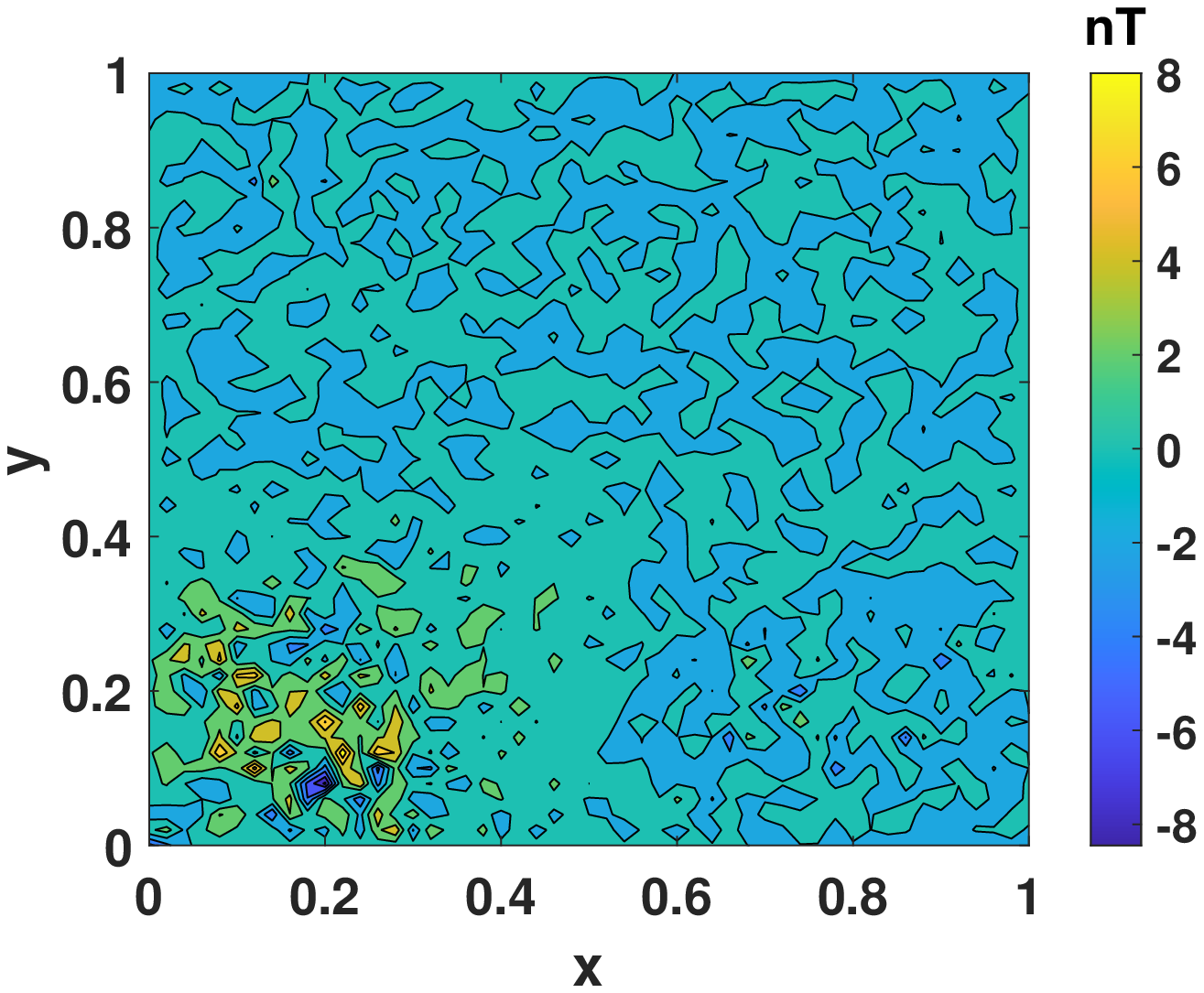}}}
\subfigure[]{\scalebox{0.4}[0.4]{\includegraphics{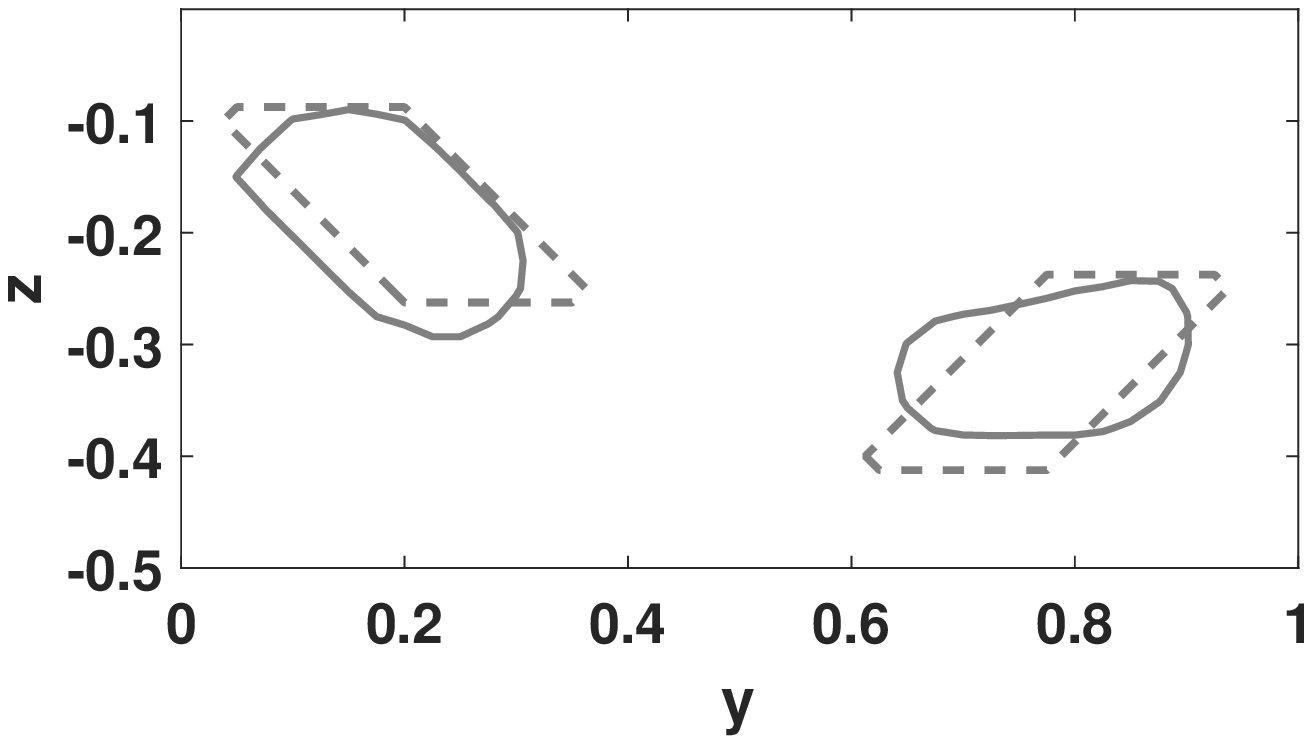}}}
\subfigure[]{\scalebox{0.4}[0.4]{\includegraphics{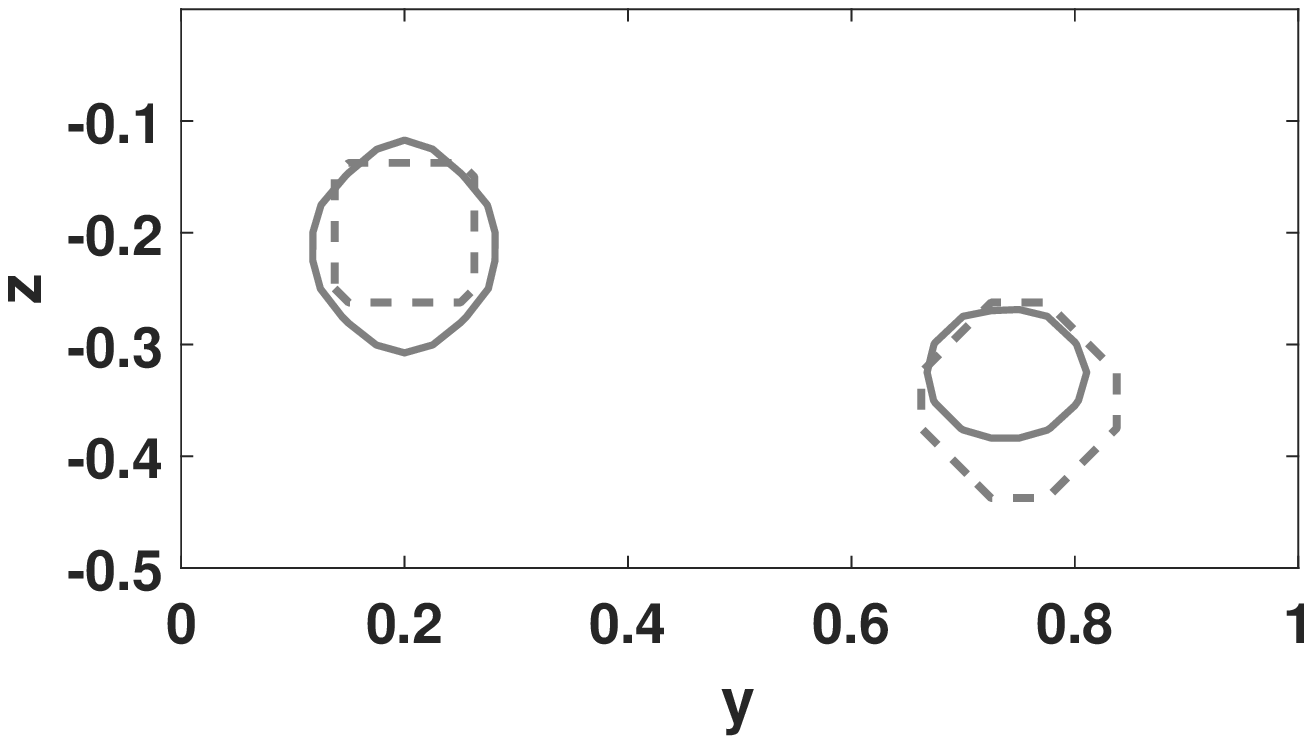}}}
\caption{Example 3. Inversion results. (a) recovered solution; (b) data discrepancy $d-d^*$; (c)-(d) cross-sections along $x=0.25$ and $x=0.75$, respectively, where the dashed line indicates the true model and the solid line indicates the recovered solution.}
\label{Fig9}
\end{figure}

\begin{figure}
\centering
\subfigure[]{\scalebox{0.4}[0.4]{\includegraphics{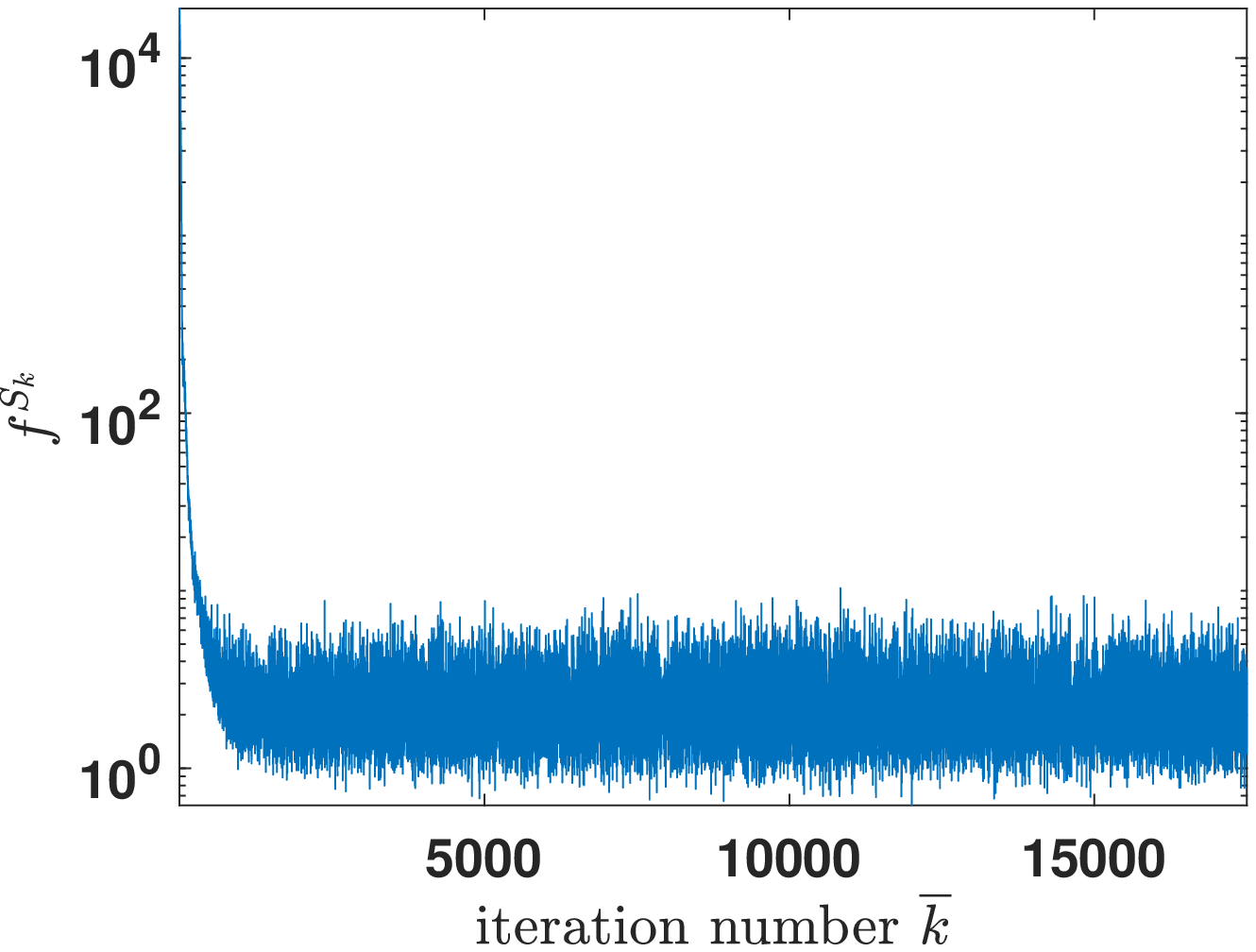}}}
\subfigure[]{\scalebox{0.4}[0.4]{\includegraphics{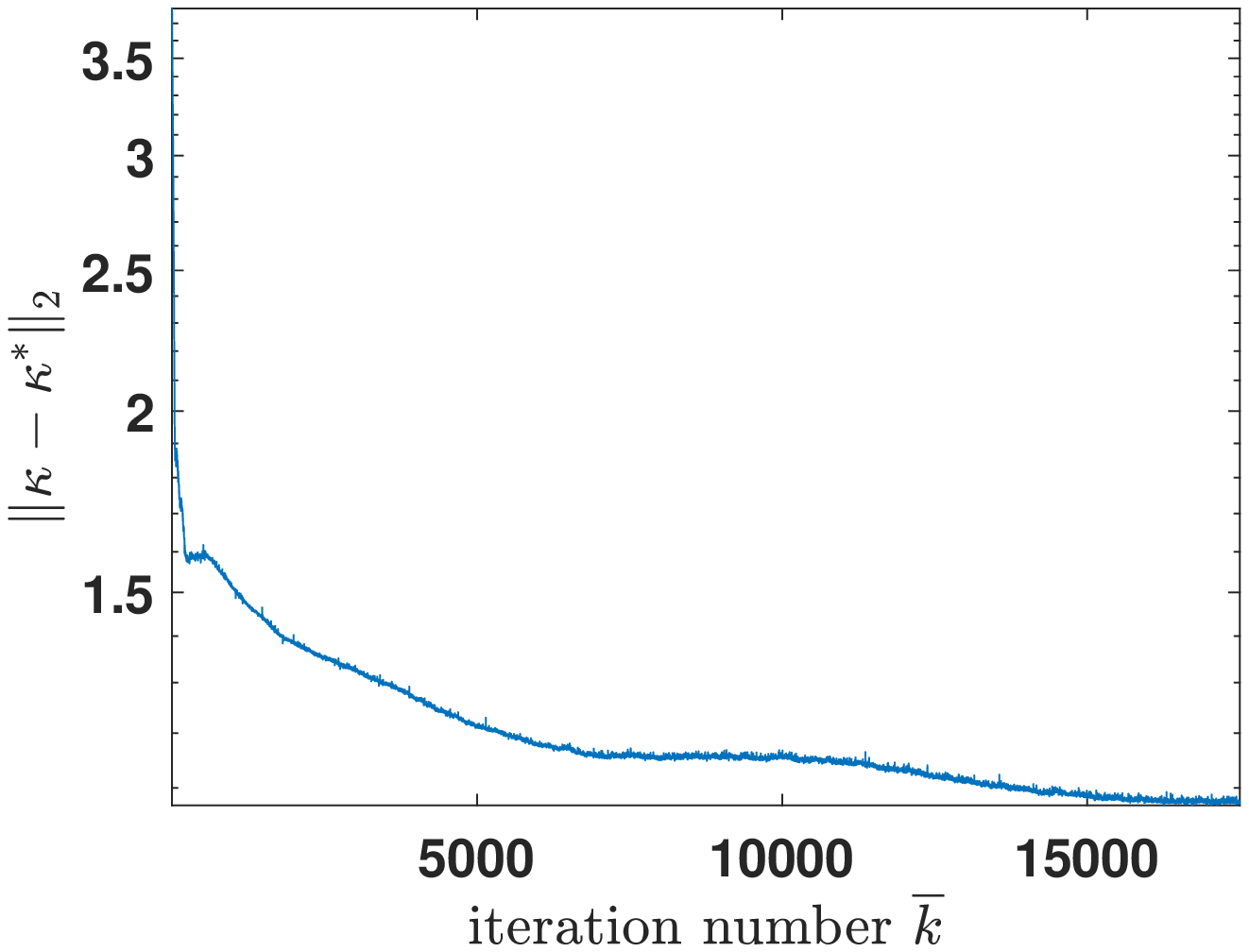}}}
\caption{Example 3. Performance of the inversion algorithm. (a) Mini-batch misfit function $f^{S_k}$; (b) $\|\kappa-\kappa^*\|_2$.}
\label{Fig10}
\end{figure}

As shown in Figure \ref{Fig8}\,(a), we consider a susceptibility model including 4 distinct magnetic sources. The cube is centered at $(0.75,0.2,-0.2)$ with side length $0.15$, the sphere is centered at $(0.8,0.75,-0.35)$ with radius $0.1$, and the two dipping prisms with different depths are enclosed in the region $\{(x,y,z) \mid x\in(0.1,0.3)\}$. We assume an inducing field with inclination and declination $(I^0,D^0)=(75^{\circ},25^{\circ})$. Figure \ref{Fig8}\,(b) shows the magnetic modulus data on the measurement surface. The inversion results are shown in Figure \ref{Fig9}, where Figure \ref{Fig9}\,(a) shows the recovered solution, Figure \ref{Fig9}\,(b) plots the data discrepancy, and Figures \ref{Fig9}\,(c),\,\ref{Fig9}\,(d) plot the cross-sections along $x=0.25$ and $x=0.75$, respectively. In addition, we present the mini-batch misfit function $f^{S_k}$ and the $l^2$-difference $\|\kappa-\kappa^*\|_2$ in Figure \ref{Fig10}. We conclude that the solution successfully recovers the susceptibility model and adequately reproduces the noisy magnetic modulus data. The performance of $\|\kappa-\kappa^*\|_2$ shows a smooth convergence although the mini-batch misfit function $f^{S_k}$ is oscillating.

\begin{figure}
\centering
\subfigure[]{\scalebox{0.4}[0.4]{\includegraphics{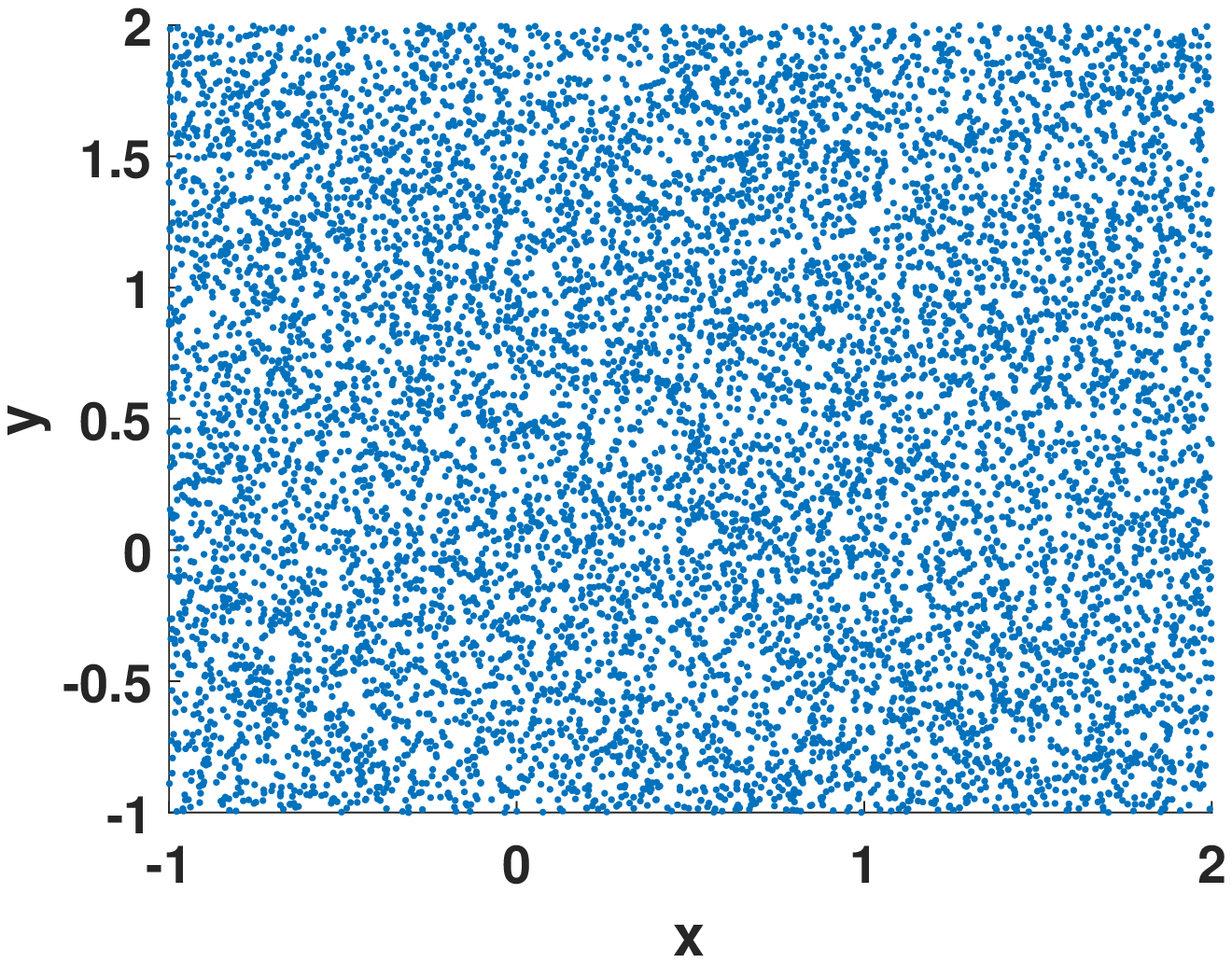}}}
\subfigure[]{\scalebox{0.4}[0.4]{\includegraphics{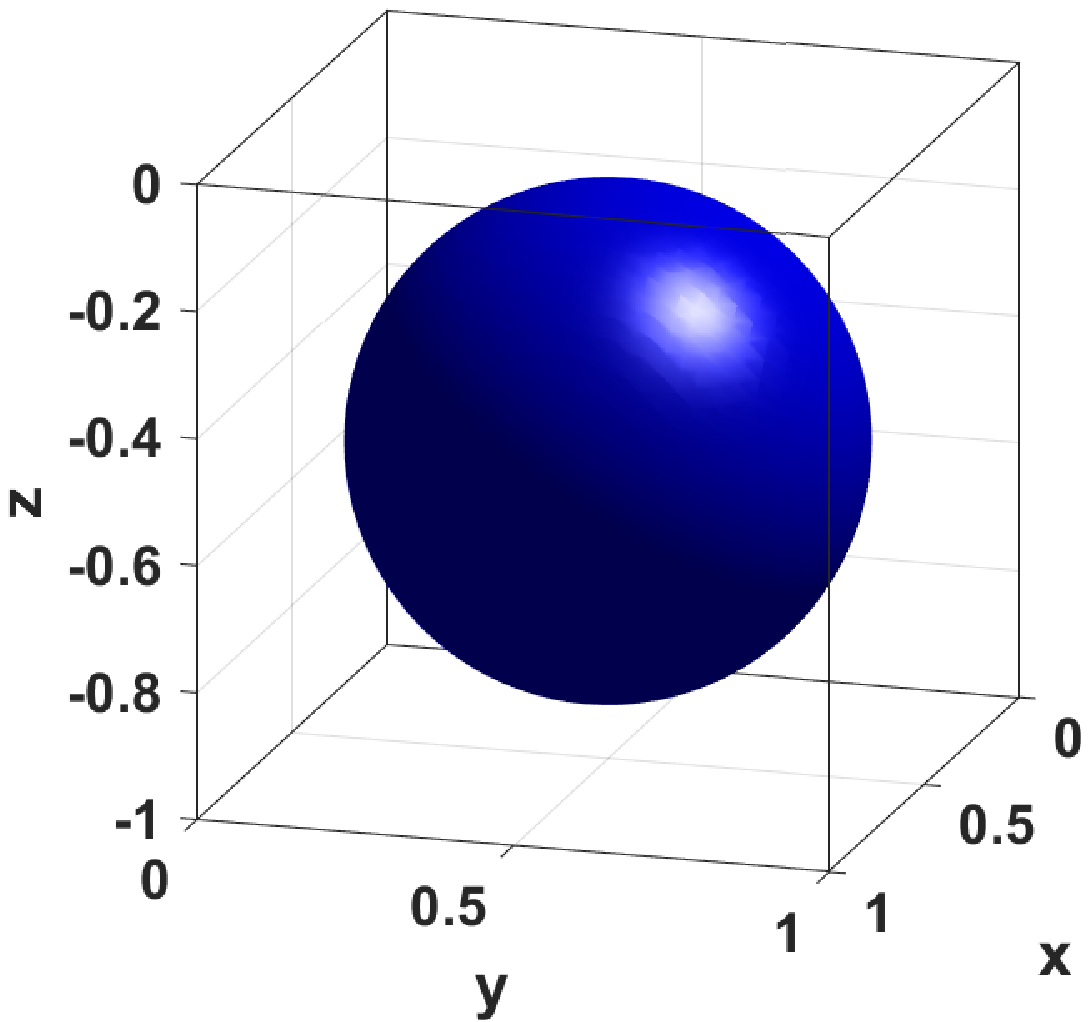}}}
\caption{Setup of examples in a cubic domain. (a) Measurement points along $\Gamma=[-1\ 2]\times[-1\ 2]\times\{z=0.1\}\,km$; (b) initial guess.}
\label{Fig11}
\end{figure}

\subsection{Examples in a cubic domain}
In the following examples, we implement the inversion algorithm in a cubic domain with larger depth. Depth resolution is generally a difficult task in magnetic inverse problem. Since the measurement surface is above the computational domain and the magnetic integral kernel decays rapidly as the distance between source and measurement increases, the responses of deep structures are likely to be suppressed by those of shallow structures and noise contaminations. As a result, depth resolution is an important criterion to evaluate the performance of a magnetic inversion algorithm. In the following examples, we include susceptibility distributions with larger variances in depth, which aims to test the capability of depth resolution for the proposed inversion algorithm.

The computational domain is $\Omega=[0\ 1]\times[0\ 1]\times[-1\ 0]\,km$, which is uniformly discretized into $41\times41\times41$ mesh grids. The measurement boundary is taken as $\Gamma=[-1\ 2]\times[-1\ 2]\times\{z=0.1\}\,km$, along which 10,000 measurement points are randomly distributed; Figure \ref{Fig11}\,(a) plots the measurement points along $\Gamma$. Again, we set the mini-batch size to be $b=200$ for the stochastic gradient descent. The initial guess of the level-set function $\phi$ is taken as:
\[
\phi_{initial}=0.4-\sqrt{(x-0.5)^2+(y-0.5)^2+(z-0.5)^2}\,,
\]
and so the zero level-set of the initial structure is a sphere as shown in Figure \ref{Fig11}\,(b). 

\subsubsection{Example 4}
\begin{figure}
\centering
\subfigure[]{\scalebox{0.4}[0.4]{\includegraphics{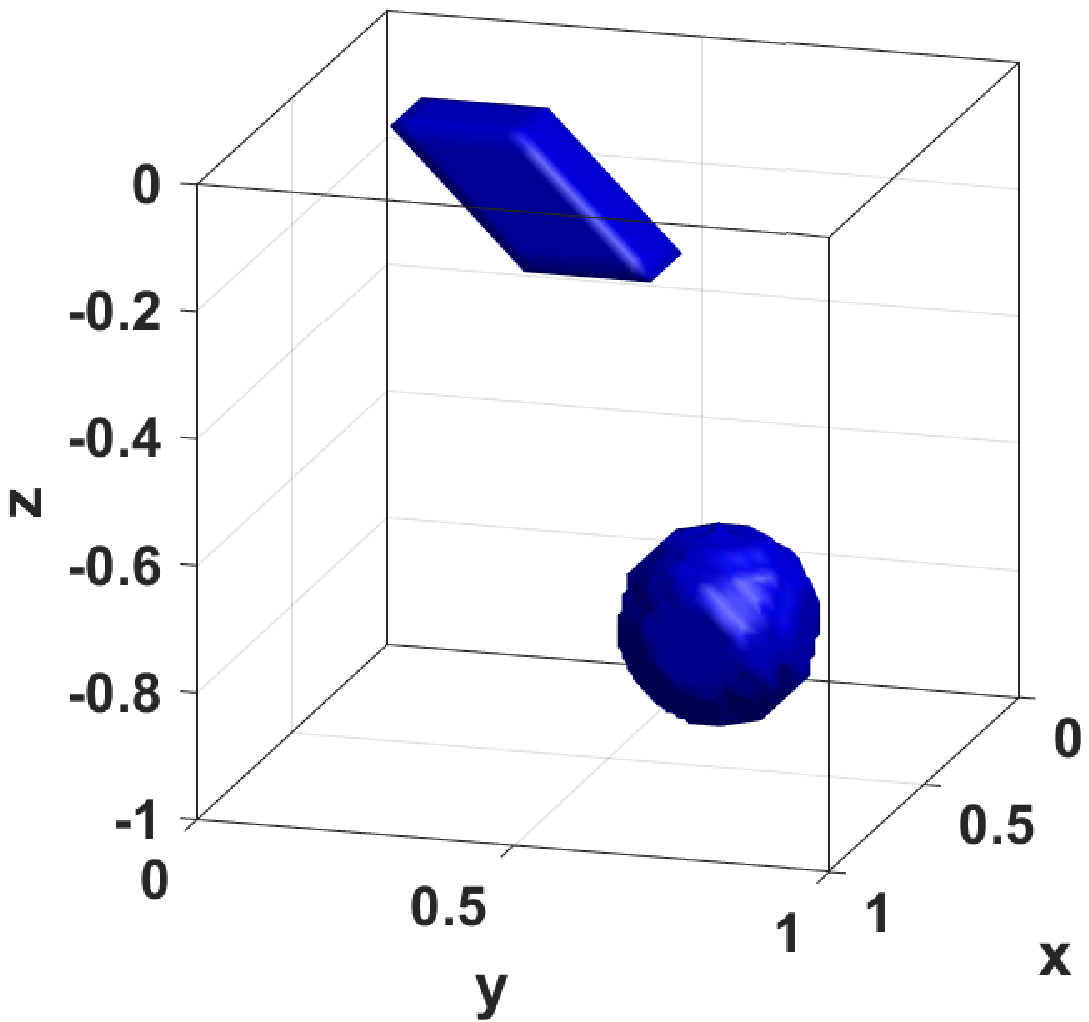}}}
\subfigure[]{\scalebox{0.4}[0.4]{\includegraphics{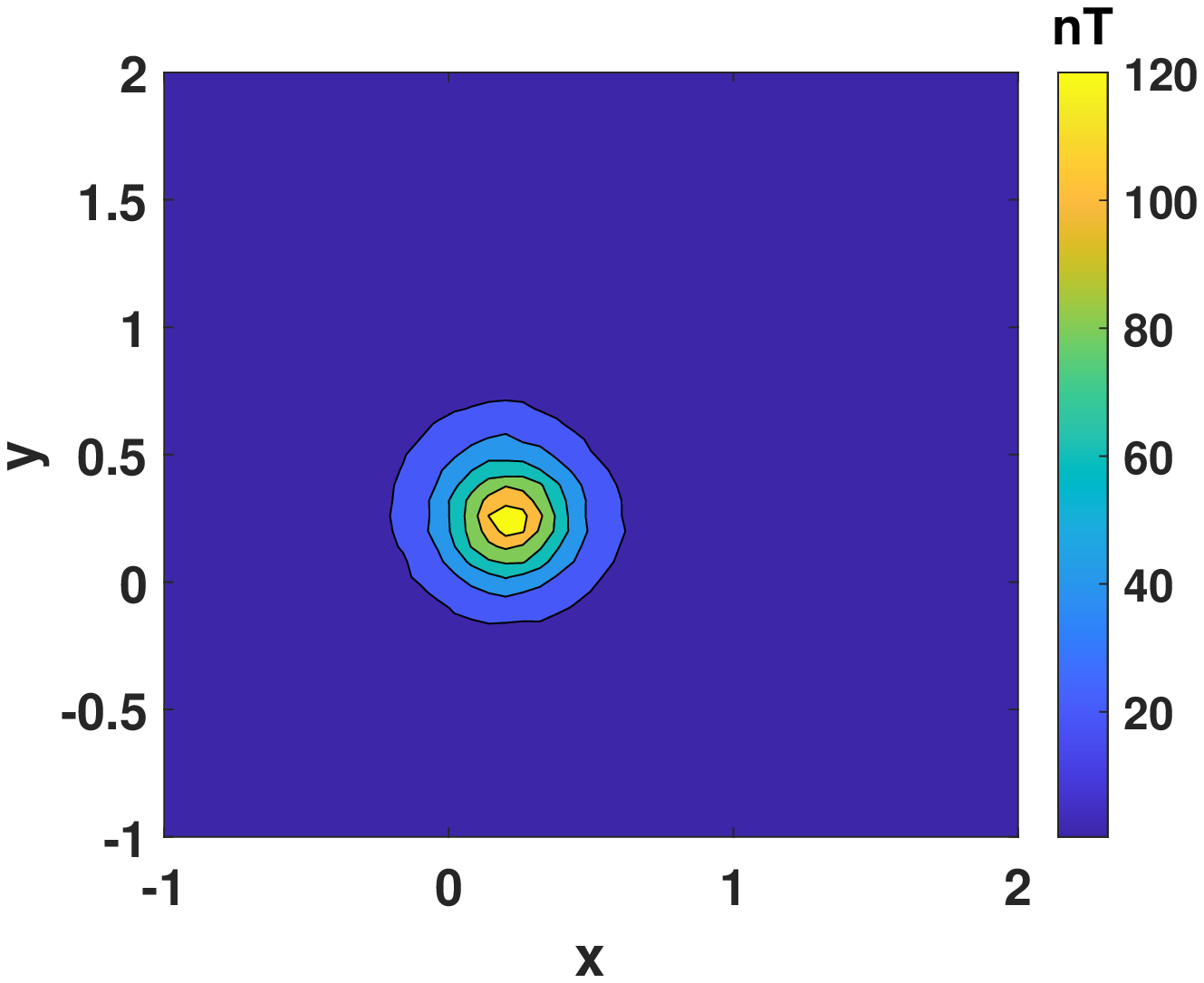}}}
\caption{Example 4. (a) True model; (b) magnetic modulus data with Gaussian noises.}
\label{Fig12}
\end{figure}

\begin{figure}
\centering
\subfigure[]{\scalebox{0.4}[0.4]{\includegraphics{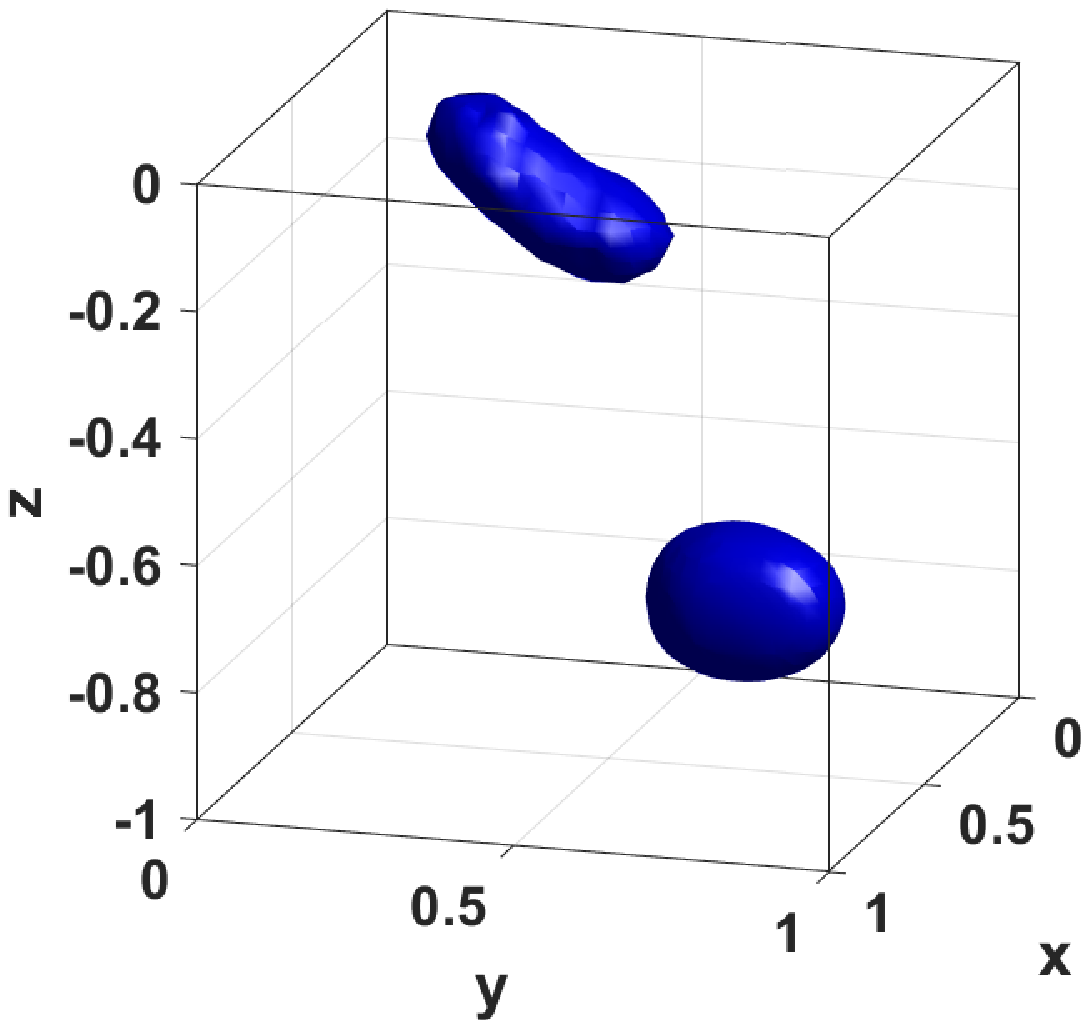}}}
\subfigure[]{\scalebox{0.4}[0.4]{\includegraphics{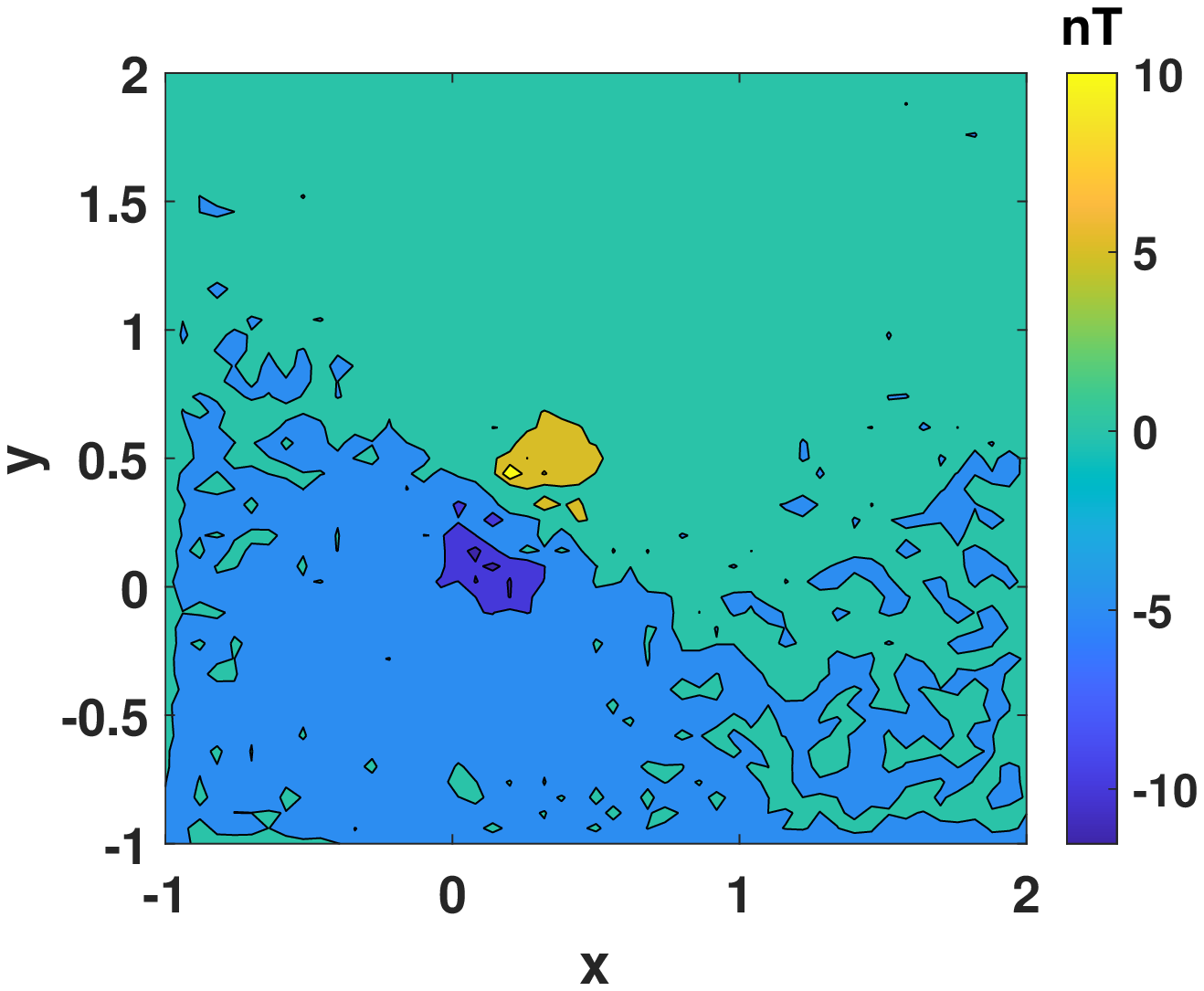}}}\\
\subfigure[]{\scalebox{0.4}[0.4]{\includegraphics{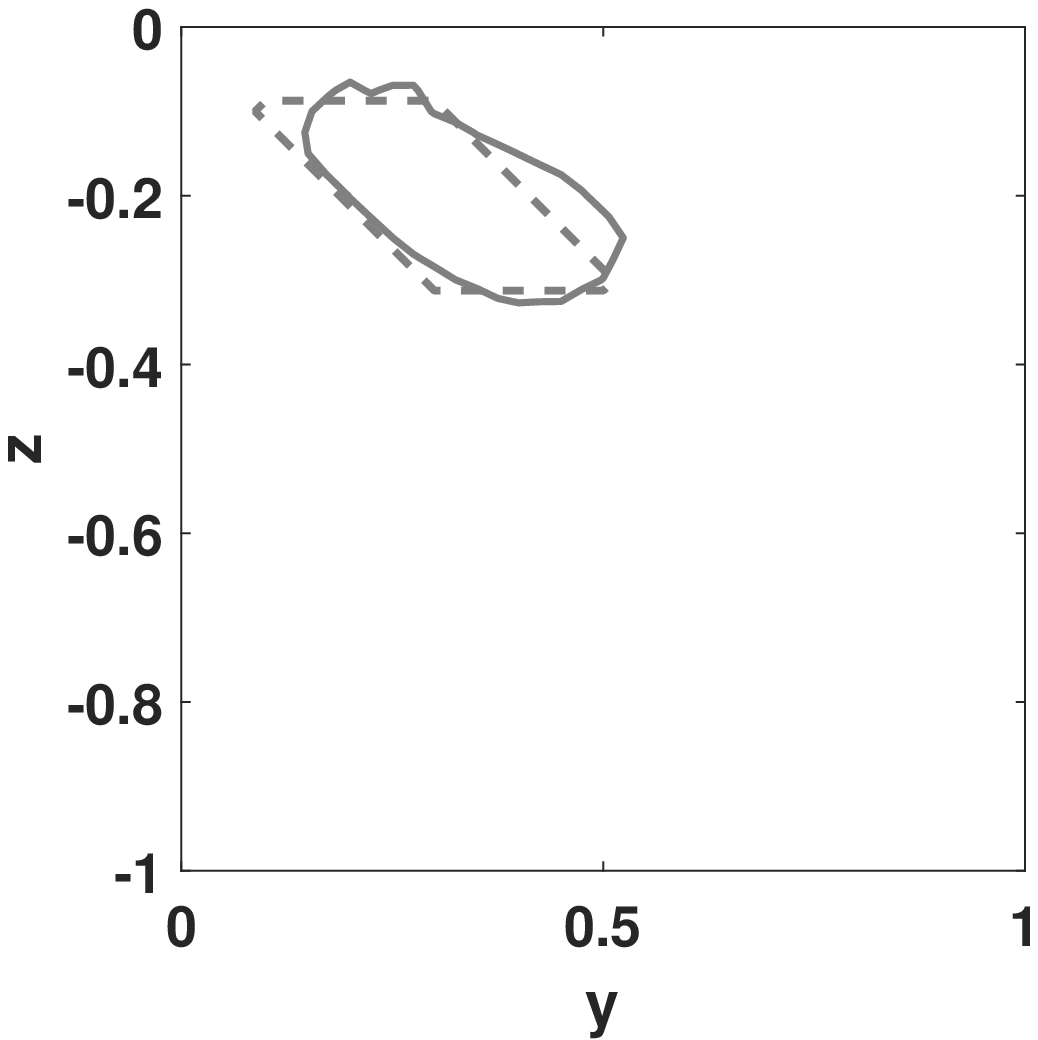}}}
\subfigure[]{\scalebox{0.4}[0.4]{\includegraphics{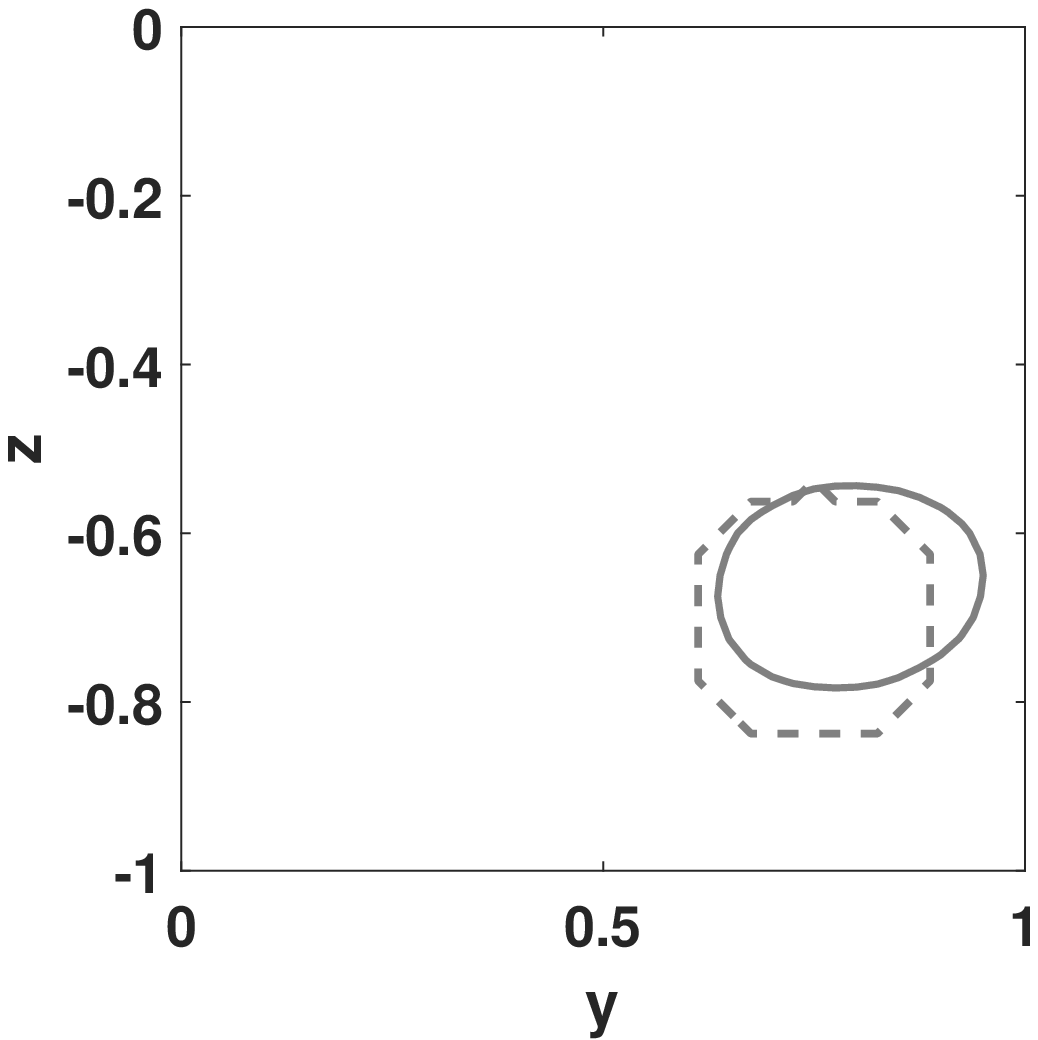}}}
\caption{Example 4. Inversion results. (a) recovered solution; (b) data discrepancy $d-d^*$; (c)-(d) cross-sections along $x=0.25$ and $x=0.75$, respectively, where the dashed line indicates the true model and the solid line indicates the recovered solution.}
\label{Fig13}
\end{figure}

\begin{figure}
\centering
\subfigure[]{\scalebox{0.4}[0.4]{\includegraphics{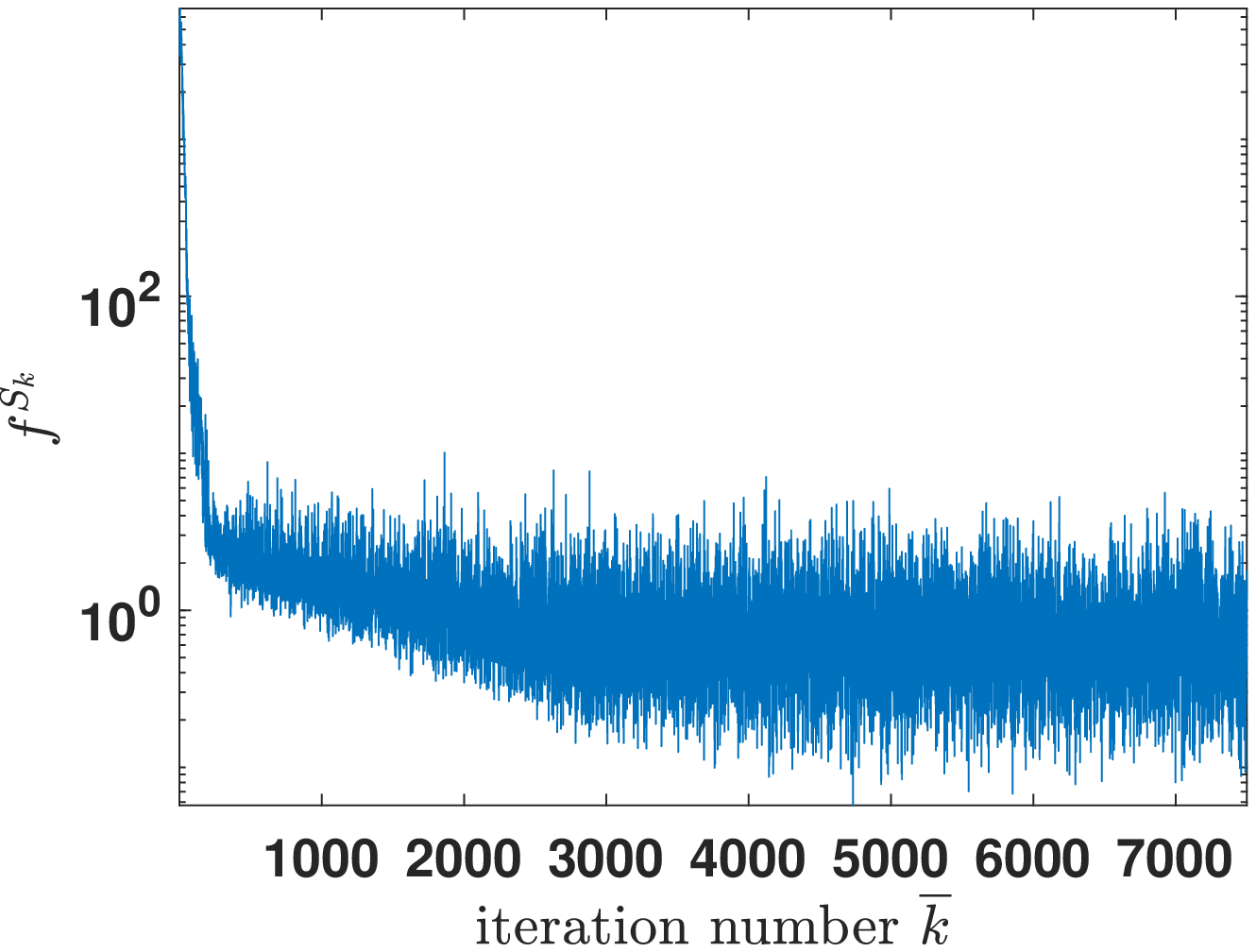}}}
\subfigure[]{\scalebox{0.4}[0.4]{\includegraphics{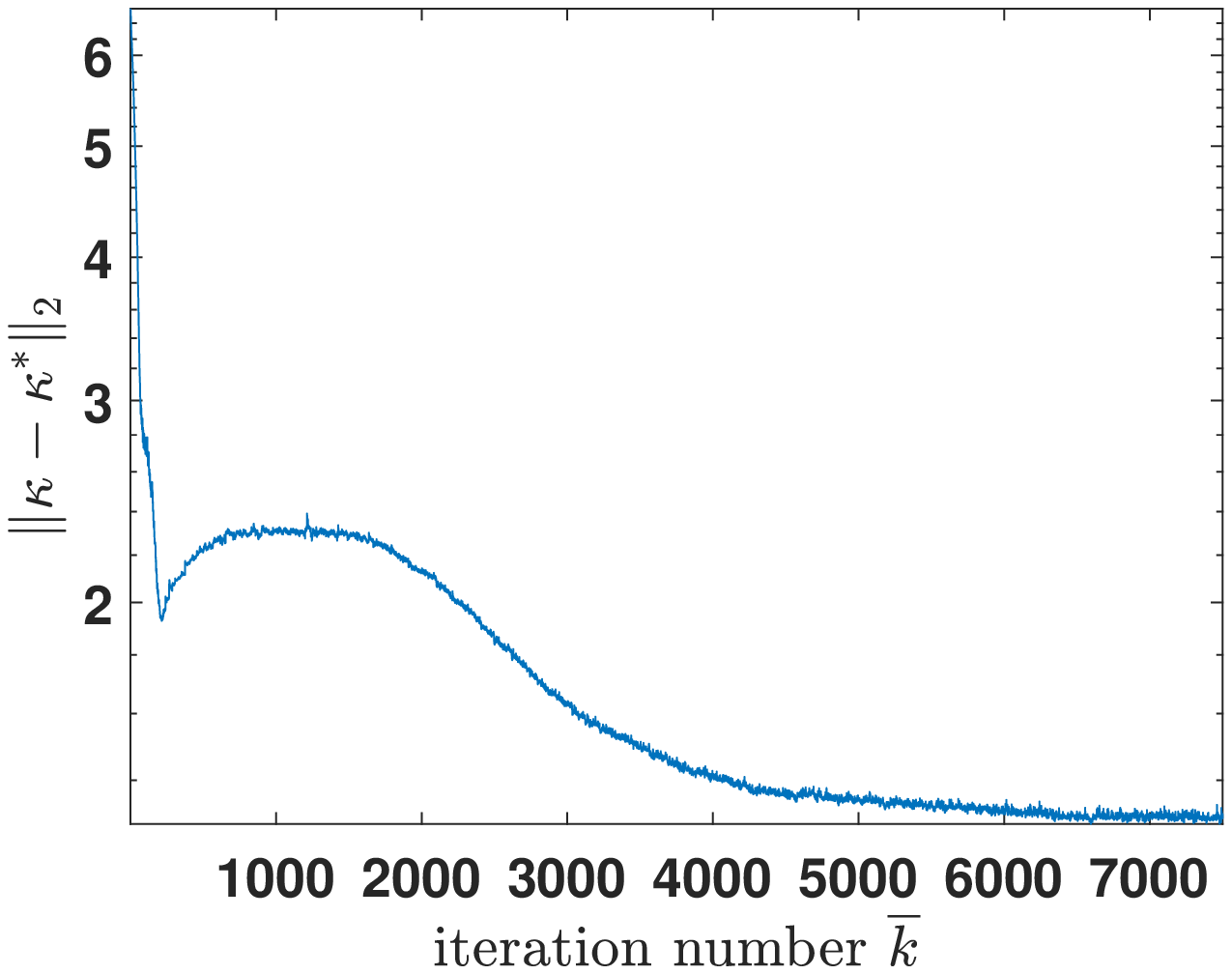}}}
\caption{Example 4. Performance of the inversion algorithm. (a) Mini-batch misfit function $f^{S_k}$; (b) $\|\kappa-\kappa^*\|_2$.}
\label{Fig14}
\end{figure}

Figure \ref{Fig12}\,(a) shows the susceptibility model, where two magnetic sources are located at different depths of the computational domain. The dipping prism is located in the shallow region, and the sphere is centered at $(0.75,0.75,-0.7)$ with radius $0.1$. We assume an inducing field with inclination and declination $(I^0,D^0)=(90^{\circ},0^{\circ})$. Figure \ref{Fig12}\,(b) shows the magnetic modulus data on the measurement surface. This is a typical example where the responses of deep structures are inapparent in the data profile. Figure \ref{Fig13} provides the inversion results using the mini-batch stochastic gradient descent algorithm with partitioned-truncated SVD. Figure \ref{Fig13}\,(a) shows the recovered solution; Figure \ref{Fig13}\,(b) plots the data discrepancy $d-d^*$; Figure \ref{Fig13}\,(c) plots the cross-section of the solution along $x=0.25$; Figure \ref{Fig13}\,(d) plots the cross-section along $x=0.75$. The performances of the mini-batch misfit function $f^{S_k}$ and the $l^2$-difference $\|\kappa-\kappa^*\|_2$ are presented in Figure \ref{Fig14}\,(a) and Figure \ref{Fig14}\,(b), respectively.

The solution successfully recovers the susceptibility model and adequately reproduces the magnetic modulus data. The depth resolution is amazing, given that the susceptibility model has large variance in depth and the spherical source in deep region has invisible response in the data profile. This test example illustrates that the inversion algorithm does not loose depth resolution in the solution, although the mini-batch stochastic gradient descent only employs part of data at each iteration and the partitioned-truncated SVD provides inaccurate approximation for matrix multiplications. In fact, since the stochastic gradient descent has the capability of escaping saddle points \cite{gehuajinyua15} and local minima \cite{kleliyua18}, we believe that it can contribute to improving depth resolution in the magnetic inverse problem.

\subsubsection{Example 5}
\begin{figure}
\centering
\subfigure[]{\scalebox{0.4}[0.4]{\includegraphics{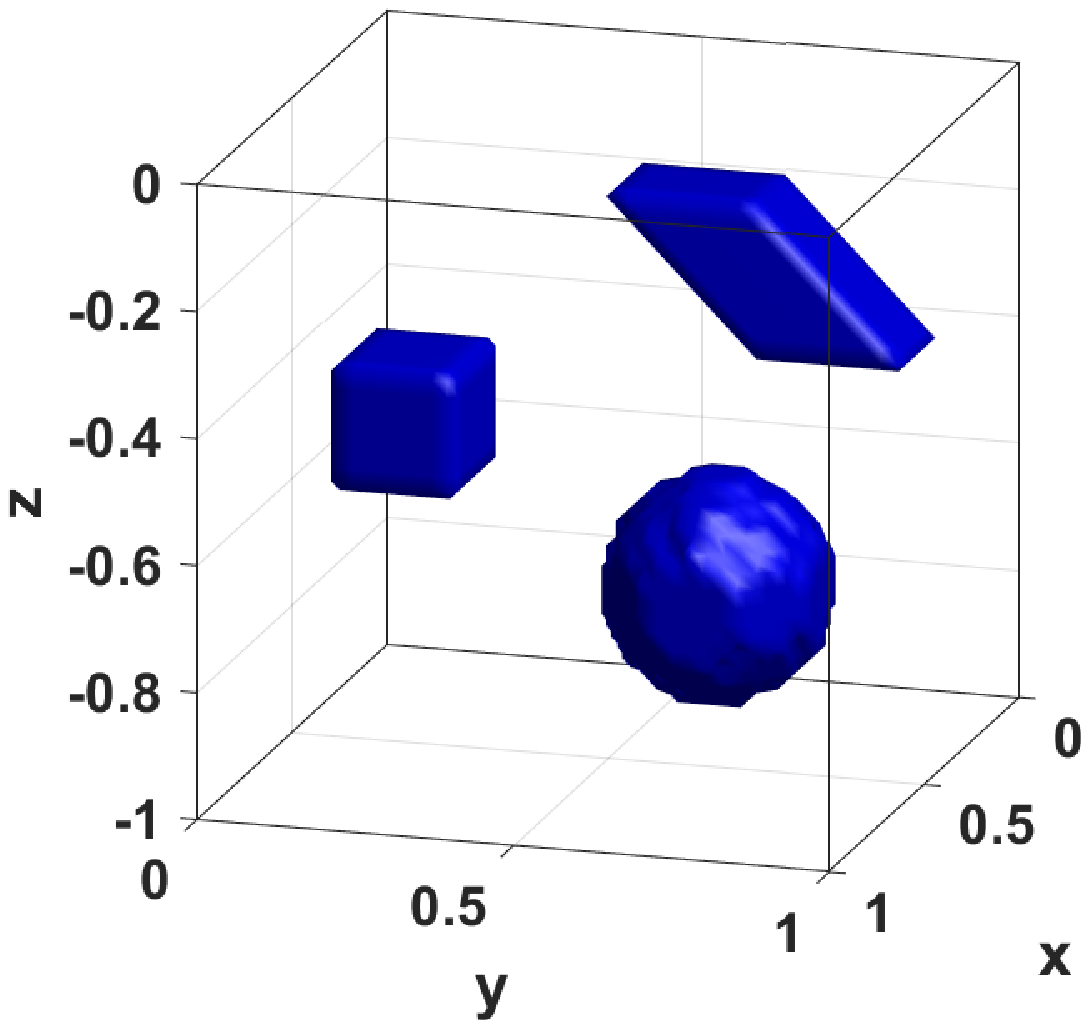}}}
\subfigure[]{\scalebox{0.4}[0.4]{\includegraphics{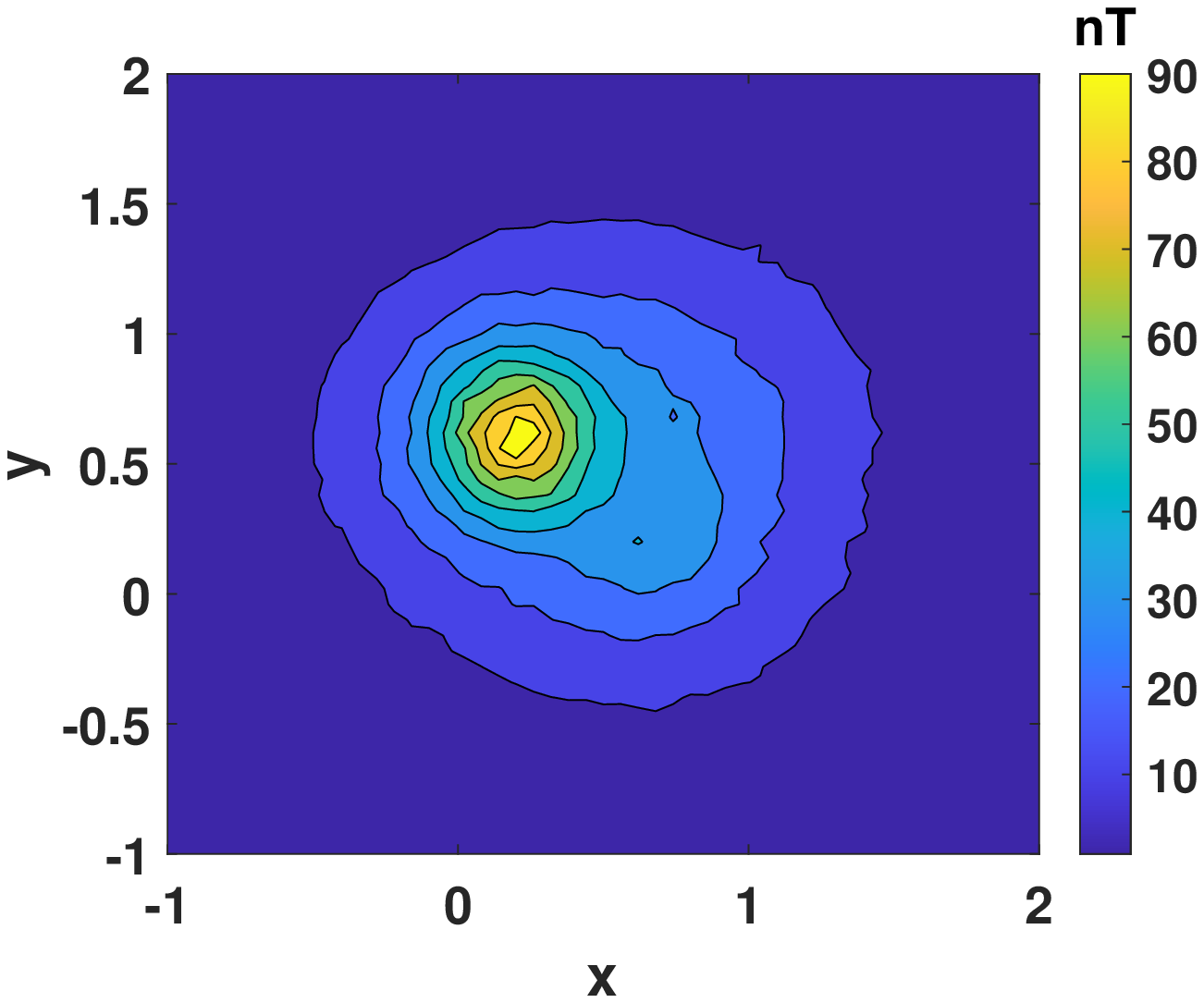}}}
\caption{Example 5. (a) True model; (b) magnetic modulus data with Gaussian noises.}
\label{Fig15}
\end{figure}

\begin{figure}
\centering
\subfigure[]{\scalebox{0.4}[0.4]{\includegraphics{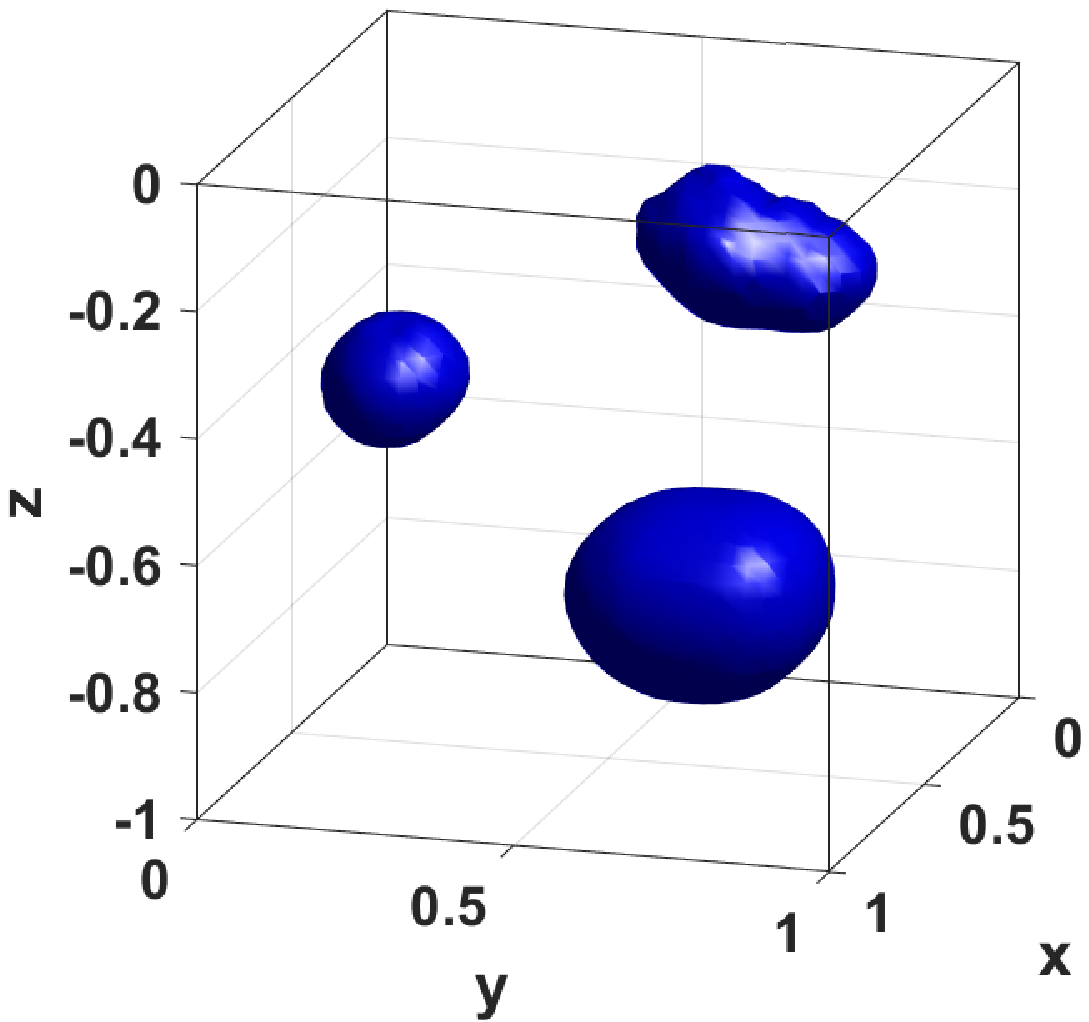}}}
\subfigure[]{\scalebox{0.4}[0.4]{\includegraphics{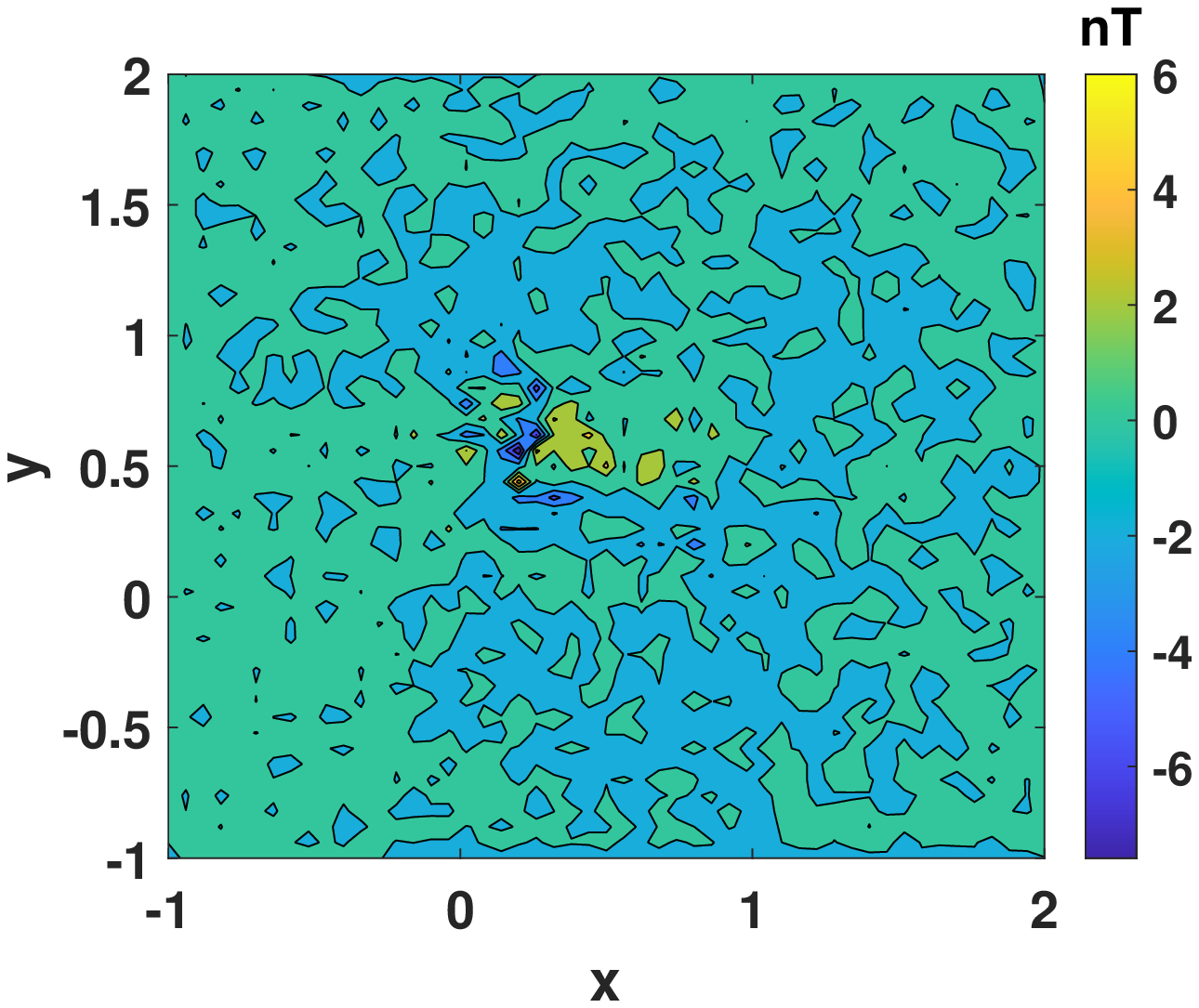}}}\\
\subfigure[]{\scalebox{0.4}[0.4]{\includegraphics{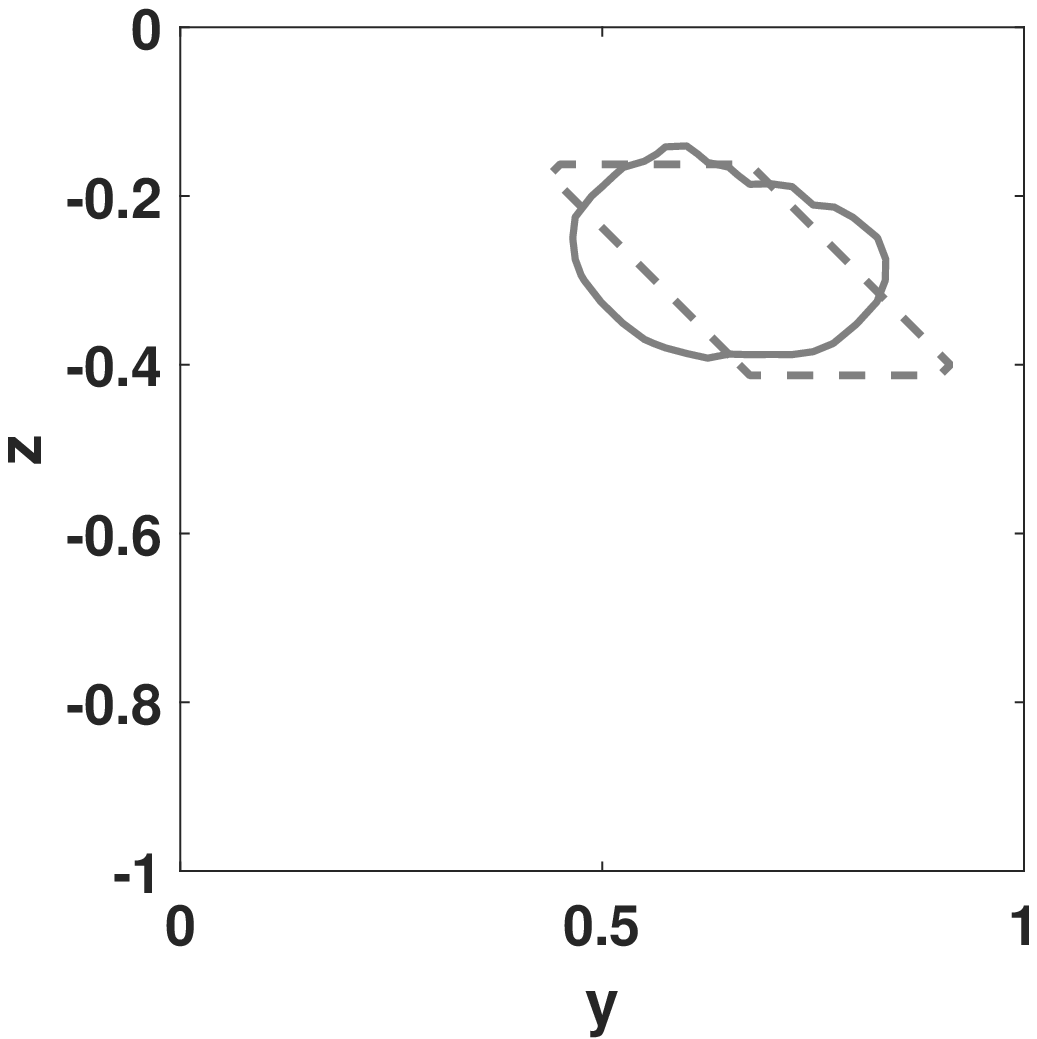}}}
\subfigure[]{\scalebox{0.4}[0.4]{\includegraphics{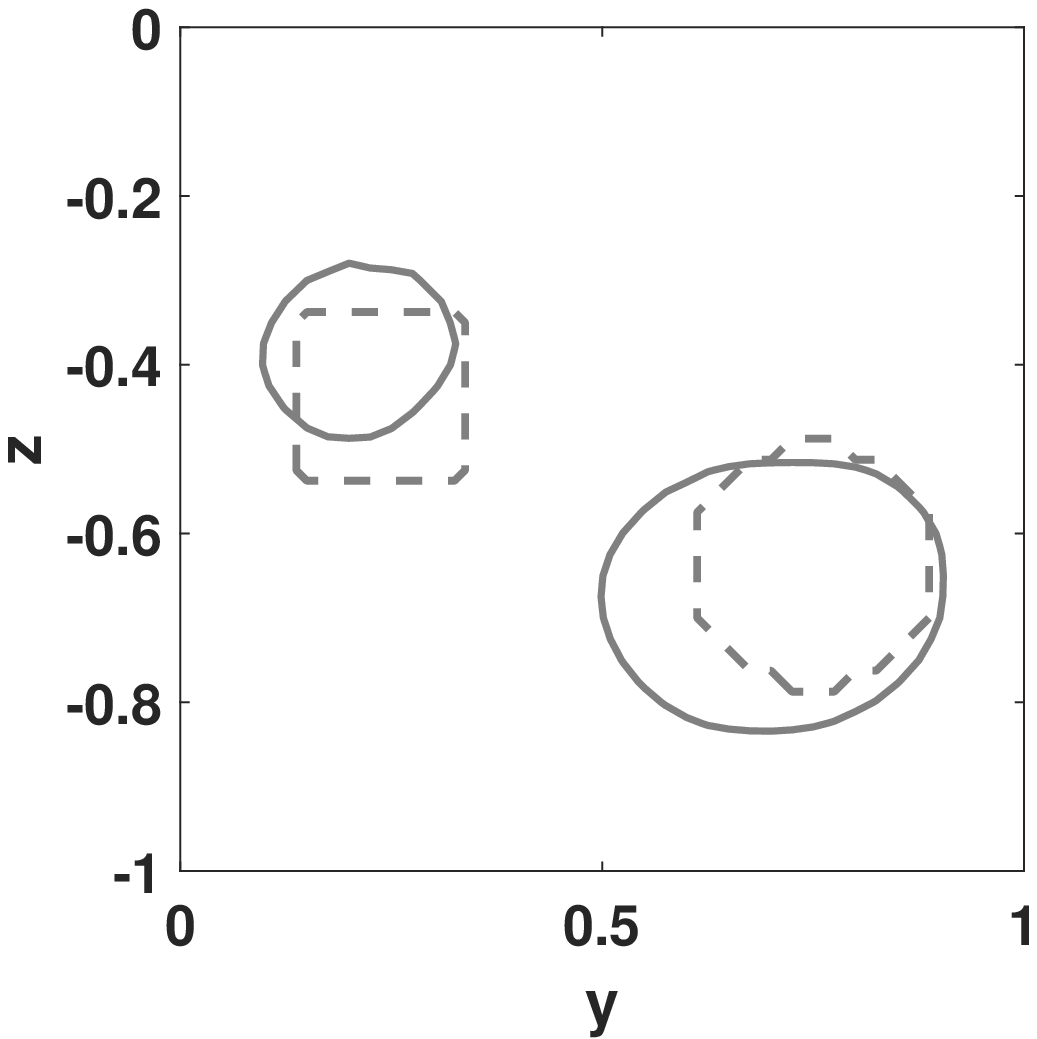}}}
\caption{Example 5. Inversion results. (a) recovered solution; (b) data discrepancy $d-d^*$; (c)-(d) cross-sections along $x=0.20$ and $x=0.65$, respectively, where the dashed line indicates the true model and the solid line indicates the recovered solution.}
\label{Fig16}
\end{figure}

\begin{figure}
\centering
\subfigure[]{\scalebox{0.4}[0.4]{\includegraphics{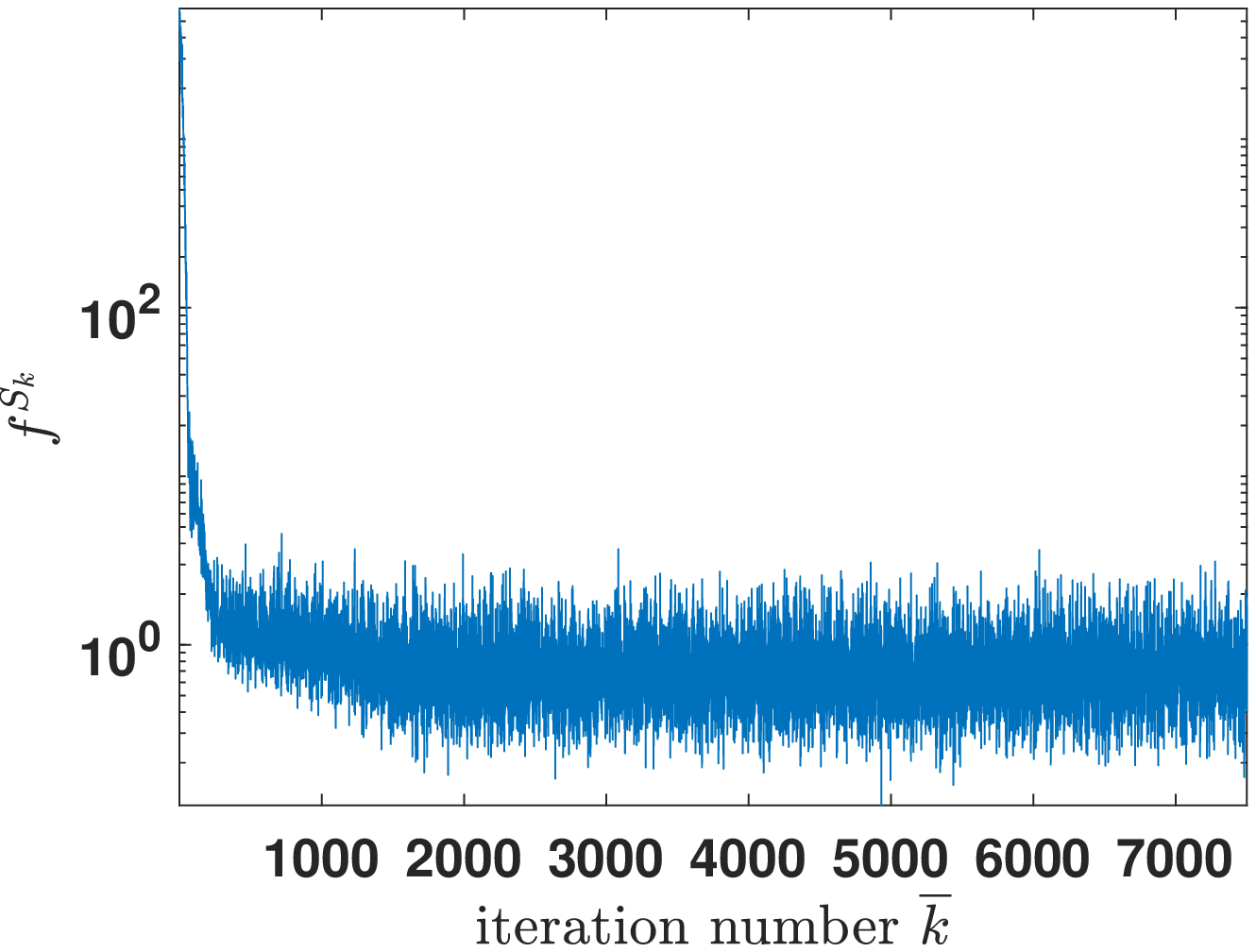}}}
\subfigure[]{\scalebox{0.4}[0.4]{\includegraphics{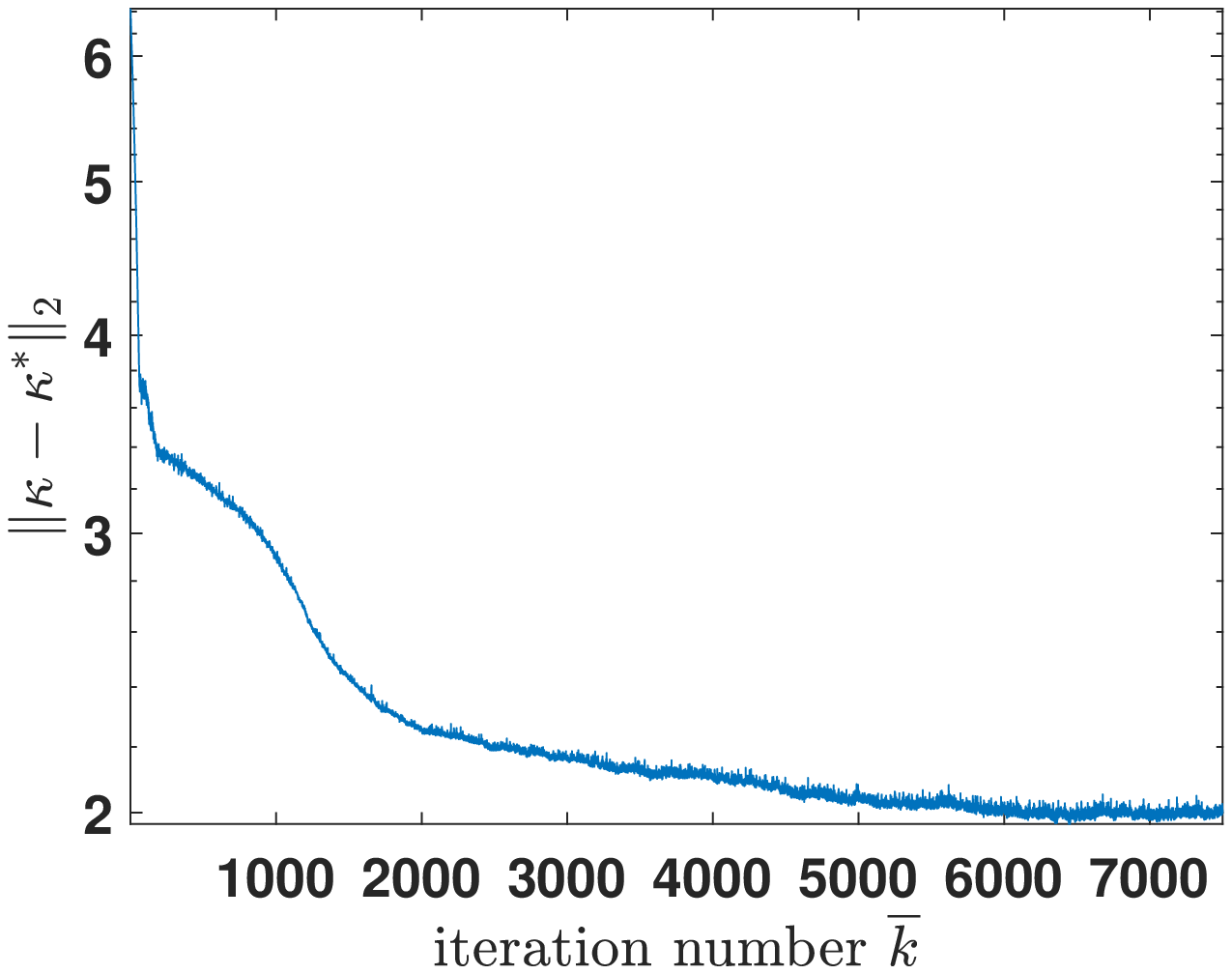}}}
\caption{Example 5. Performance of the inversion algorithm. (a) Mini-batch misfit function $f^{S_k}$; (b) $\|\kappa-\kappa^*\|_2$.}
\label{Fig17}
\end{figure}

In this example, we consider a susceptibility model including 3 magnetic sources with different depths. The true model is shown in Figure \ref{Fig15}\,(a). The measurement magnetic modulus data with Gaussian noises are shown in Figure \ref{Fig15}\,(b), where we assume an inducing field with inclination and declination $(I^0,D^0)=(75^{\circ},25^{\circ})$. The inversion results are shown in Figure \ref{Fig16}, where Figure \ref{Fig16}\,(a) shows the recovered solution, Figure \ref{Fig16}\,(b) plots the data discrepancy, and Figures \ref{Fig16}\,(c),\,\ref{Fig16}\,(d) plot two cross-sections of the solution. The performances of the mini-batch misfit function $f^{S_k}$ and the $l^2$-difference $\|\kappa-\kappa^*\|_2$ are presented in Figure \ref{Fig17}. We conclude that the inversion algorithm successfully recovers the susceptibility model and adequately reproduces the measurement data. The depth resolution is reasonable if not perfect, given that the susceptibility model has distinct variance in depth, and the task of recovering depth is generally difficult. With the techniques of stochastic gradient descent and partitioned-truncated SVD, the proposed inversion algorithm is able to deal with large-scale measurement data efficiently in the magnetic inverse problem.

\section{Conclusion and discussion}
We have proposed a mini-batch stochastic gradient descent approach with partitioned-truncated SVD for large-scale inverse problems of magnetic modulus data. We employ a level-set formulation for susceptibility model and recover the volume susceptibility distribution from nonlinear magnetic modulus data by solving the level-set function. To deal with massive amount of measurement data, we propose a stochastic gradient descent approach to solve the level-set optimization problem with large scale. When evaluating the stochastic gradients, we consider a mini-batch strategy to reduce variances, and we use the without-replacement sampling with random reshuffling to explore full data set at every epoch. Realizing that the iteration process can be formulated as a Hamilton-Jacobi equation, we propose a step-size rule for the stochastic gradient descent according to the Courant-Friedrichs-Lewy condition of the evolution PDE. To further improve the computational efficiency, we take advantage of the decaying property of the magnetic integral kernel, and propose a partitioned-truncated SVD for matrix multiplications in the context of stochastic gradient descent.

Numerical examples are included to illustrate the efficacy of the proposed method. The solutions successfully recover the susceptibility models and adequately reproduce the magnetic modulus data. The inversion algorithm is able to reconstruct susceptibility distributions with large variances in depth. As an important criterion to evaluate the performance of the magnetic inversion algorithm, the depth resolution is very good in all the test examples. We believe that the stochastic gradient approach can help to improve depth resolution due to its capability of escaping saddle points and local minima. With the techniques of stochastic gradient descent and partitioned-truncated SVD, the proposed method enables us to efficiently process large-scale magnetic data which was infeasible due to the restriction of hardware resource.

The inverse problem of magnetic modulus data can be formulated as the general form
\begin{equation} \label{general1}
d=\sigma_2\circ L\circ \sigma_1(\phi)\,:\quad \mathrm{reconstruct}\ \phi\ \mathrm{from}\ d\,,
\end{equation}
where $\sigma_1$ denotes the nonlinear level-set formulation as shown in (\ref{eqn5}), $L$ denotes the linear integral as shown in (\ref{eqn1}), and $\sigma_2$ denotes the nonlinear stacking as shown in (\ref{eqn3}); after discretization, $L$ is the linear operator of matrix multiplications. The forward operator in formula (\ref{general1}) is similar to the architecture of deep neural network:
\begin{equation} \label{general2}
\Psi=\sigma_N\circ L_N\circ\cdots\sigma_{i}\circ L_{i}\circ \cdots \circ\sigma_1\circ L_1\,,
\end{equation}
where $\sigma_i$'s denote nonlinear activating functions, and $L_i$'s denote linear operators including affine mappings, convolutions, etc. In this work, we have proposed an efficient stochastic gradient descent approach with partitioned-truncated SVD for the inverse problem of (\ref{general1}). The algorithm can be generalized to solve the inverse problem of the deep neural network as shown in (\ref{general2}), provided that the linear operators $L_i$ have some sort of low-rank property. The inverse problem should be stated as follows: given a well-trained deep neural network $\Psi=\sigma_N\circ L_N\circ\cdots\circ\sigma_1\circ L_1$ and its output response $d=\Psi(\phi)$, reconstruct the input signal $\phi$. We will do further explorations on this subject.

\section*{Acknowledgments}
Wenbin Li is supported by NSFC (grant no. 41804096), Natural Science Foundation of Guangdong Province (grant no. 2018A030313341), and Natural Science Foundation of Shenzhen (grant no. JCYJ20190806144005645).

\newpage
\bibliographystyle{plain}

\end{document}